%% file: CPR_SBP_Error__arXiv.tex
\documentclass[a4paper, 10pt]{scrartcl}
\pdfoutput=1

\title{Error Boundedness of \revPP{Discontinuous Galerkin Methods} with Variable Coefficients}

\author{Philipp \"Offner and Hendrik Ranocha}
\date{November 12, 2018}

\newcommand{\revPP}[1]{\textcolor{black}{#1}} 
\newcommand{\revsec}[1]{\textcolor{black}{#1}}

\include{includes_arXiv}

\begin{document}

\maketitle

\begin{abstract}

\input{0_Abstract}
\end{abstract}

\tableofcontents

\input{1_Introduction}

\input{2_Model_problem}

\input{3_Numerical_Methods}

\input{4_Error_bound_semidiscrete_Lobatto}

\input{5_Error_bound_semidiscrete_Gauss}

\input{6_Numerics}

\input{7_Generalisation}

\input{8_Summary}

\input{Appendix}

\section*{Acknowledgements}

This first author was supported by SVF project ``Solving advection dominated 
problems with high order schemes with polygonal meshes: application to compressible
and incompressible flow problems''
and the second author was supported by the German Research Foundation (DFG, 
Deutsche Forschungsgemeinschaft) under Grant SO~363/14-1.

\bibliographystyle{spmpsci}
\bibliography{literature}

\end{document}

%% file: includes_arXiv.tex
\usepackage[utf8]{luainputenc}
\usepackage[english]{babel}
\usepackage{csquotes}

\usepackage[a4paper, textheight=622pt, textwidth=468pt, centering]{geometry}

\usepackage[plainpages=false,pdfpagelabels,hidelinks]{hyperref}
\makeatletter
\hypersetup{pdfauthor={\@author}}
\hypersetup{pdftitle={\@title}}
\makeatother

\usepackage{cite}

\usepackage{amsmath}
\allowdisplaybreaks
\usepackage{amssymb}
\usepackage{commath}
\usepackage{mathtools}

\usepackage{siunitx}

\usepackage{amsthm}
\theoremstyle{plain}
  \newtheorem{theorem}{Theorem}[section]
  
  \theoremstyle{definition}
  
\newtheorem{remark}[theorem]{Remark}
\newtheorem{Result}[theorem]{Result}
  
\usepackage{color}
\usepackage{graphicx}
\usepackage[small]{caption}
\usepackage{subcaption}

\usepackage{pgfplots}
\pgfplotsset{compat=1.11}
\usepgfplotslibrary{external} 
\begingroup\expandafter\expandafter\expandafter\endgroup
\expandafter\ifx\csname pdfsuppresswarningpagegroup\endcsname\relax
\else
\fi

\usepackage{booktabs}
\usepackage{rotating}
\usepackage{multirow}

\usepackage{calc}
\usepackage{xparse}

\renewcommand{\vec}[1]{\underline{#1}}
\NewDocumentCommand{\mat}{mo}{%
  \IfValueTF{#2}{%
    \underline{\underline{#1}}{#2}
  }{%
    \underline{\underline{#1}}\,
  }%
}
\newcommand{\diag}[1]{\operatorname{diag}\left(#1\right)}

\renewcommand{\div}{\operatorname{div}}

\renewcommand{\d}{\operatorname{d}}

\renewcommand{\L}{\mathbf{L}}
\renewcommand{\P}{\mathbb{P}}
\newcommand{\bU}{\mathbf{U}}
\newcommand{\bF}{\mathbf{F}}

\newcommand{\Ep}{\mathbb{E}}
\newcommand{\Eps}{\mathbf{E}}

\newcommand{\artanh}{\operatorname{artanh}}

\newcommand{\fnum}{f^{\mathrm{num}}}
\newcommand{\vecfnum}{\vec{f}^{\mathrm{num}}}
\newcommand{\vecfnumk}{\vec{f}^{\mathrm{num},k}}

\newcommand{\umin}{u_{-}}
\newcommand{\uadd}{u_{+}}

\newcommand{\amin}{a_{-}}
\newcommand{\aadd}{a_{+}}

\newcommand{\vecepsilon}{\vec{\epsilon}_1^{\mathrm{num},k}}

\newcommand{\Ip}{\mathbb{I}}

\renewcommand{\epsilon}{\varepsilon}
\renewcommand{\phi}{\varphi}

\newcommand{\hrho}{\hat{\rho}}
\newcommand{\hu}{\hat{u}}
\newcommand{\vect}[1]{\begin{pmatrix}#1\end{pmatrix}}

\renewcommand{\bf}{\mathbf}
\renewcommand{\r}{\right}
\renewcommand{\l}{\left}

\newcommand{\est}[1]{\left\langle#1\right\rangle}

\newsavebox{\DelimiterBox}
\newlength{\DelimiterHeight}
\newlength{\DelimiterDepth}
\newsavebox{\ArgumentBox}
\newlength{\ArgumentHeight}
\newlength{\ArgumentDepth}
\newlength{\ResizedDelimiterHeight}

\newcommand{\jump}[1]{%
  \savebox{\ArgumentBox}{$ \displaystyle #1 $}%
  \settoheight{\ArgumentHeight}{\usebox{\ArgumentBox}}%
  \settodepth{\ArgumentDepth}{\usebox{\ArgumentBox}}%
  \savebox{\DelimiterBox}{$ [\!\![ $}%
  \settoheight{\DelimiterHeight}{\usebox{\DelimiterBox}}%
  \settodepth{\DelimiterDepth}{\usebox{\DelimiterBox}}%
  \setlength{\ResizedDelimiterHeight}{\maxof{1.2\ArgumentHeight}{\DelimiterHeight}}
  \!
  \resizebox{\width}{\ResizedDelimiterHeight}{ [\![ }
  \mkern-6.5mu
  #1
  \mkern-6.5mu
  \resizebox{\width}{\ResizedDelimiterHeight}{ ]\!] }
  \!\!
}

%% file: 0_Abstract.tex
For practical applications, the long time behaviour of the error of numerical solutions
to time-dependent partial differential equations is very important. Here, we investigate
this topic in the context of hyperbolic conservation laws and flux reconstruction schemes,
\revPP{focusing on the schemes in the discontinuous Galerkin spectral element framework}.
For linear problems with constant coefficients, it is well-known in the literature that
the choice of the numerical flux (e.g. central or upwind) and the selection of the polynomial
basis (e.g. Gauß-Legendre or Gauß-Lobatto-Legendre) affects both the growth rate and the
asymptotic value of the error.

Here, we extend these investigations of the long time error to variable coefficients using
both Gauß-Lobatto-Legendre and Gauß-Legendre nodes as well as several numerical fluxes.
We derive conditions guaranteeing that the errors are still bounded in time.
Furthermore, we analyse the error behaviour under these conditions and demonstrate in
several numerical tests similarities to the case of constant coefficients. However, if
these conditions are violated, the error shows a completely different behaviour. Indeed,
by applying central numerical fluxes, the error increases without upper bound while
upwind numerical fluxes can still result in uniformly bounded numerical errors. An
explanation for this phenomenon is given, confirming our analytical investigations.

%% file: 1_Introduction.tex
\section{Introduction}
\label{sec:introduction}

The investigation of the error of numerical solutions to hyperbolic conservation
laws has received much interest in the literature
\cite{hesthaven2002nodal,cohen2006spatial,nordstrom2003high,abarbanel2000error,koley2009higher,nordstrom2007error,kopriva2017error,oeffner2018error,zhang2004error}.
In some of these papers, a linear error growth (or nearly linear growth) in
time is observed, while the numerical error is bounded uniformly in time for others.
In \cite{nordstrom2007error}, the author explains under what conditions the error
is or is not bounded in time \revsec{if a linear problem with constant coefficients is considered}. Using finite difference approximations with
summation-by-parts (SBP) operators and simultaneous approximation terms (SATs),
the error behaviour depends on the choice of boundary procedure of the problem.
If one catches the waves in cavities or with periodic boundary conditions,
linear growth is observed like in \cite{hesthaven2002nodal},
whereas for inflow-outflow problems one obtains uniform boundedness in time.
In other words, if the boundary approach has sufficient dissipation,
the error is bounded. It does not depend on the internal discretisation.

This investigation is extended to the discontinuous Galerkin spectral
element methods (DGSEM) in \cite{kopriva2017error}
and to Flux Reconstruction (FR) schemes in \cite{oeffner2018error}.
Different from \cite{nordstrom2007error}, using DG or FR methods, the internal
approximation has an influence on the behaviour of the error, since there are
additional parameters. The choices of numerical fluxes (upwind and central)
and polynomial bases (Gauß-Lobatto-Legendre or Gauß-Legendre) have an impact on
the magnitude of the error and the speed at which the asymptotic error is reached.

In all of these works \cite{nordstrom2007error,kopriva2017error,oeffner2018error},
the model problem under consideration is a linear advection equation with constant
coefficients. In this paper, we extend these investigations by considering variable
coefficients. The introduction of variable coefficients leads to stability issues
and problems in the discretisation of the numerical fluxes as described
in \cite{ranocha2018generalised}. Using split forms in the spatial discretisation
\cite{nordstrom2006conservative, fisher2013discretely},
we are able to construct an error equation in the spirit of \cite{kopriva2017error}
for our new model problem.

Furthermore, using this error equation, we formulate conditions on the variable
coefficients to guarantee that the error is still bounded uniformly in time.
Here, it will be essential that the first derivative of the variable coefficient
$a(x)$ is positive. In numerical tests, we demonstrate that if these conditions
are fulfilled, the errors behave like in the case of constant coefficients.
If these conditions are not satisfied, we have a different behaviour. If central
numerical fluxes are applied, the errors tend to infinity, whereas the errors
using upwind fluxes in the calculation may still remain bounded uniformly in time.
This matches our analysis and the conditions which we derive in the analytical
investigations in sections \ref{sec:Error_Lobatto} and \ref{sec:Error_Gauss}.

The paper is organised as follows:
In the second section, we introduce the model problem and repeat the stability
analysis from the continuous point of view.
In section~\ref{sec:Num_methods}, the main idea of SBP-FR methods and the concrete
schemes are repeated. Then, we present the different numerical fluxes under consideration
and introduce the main focus of our study, the numerical errors. We repeat some
approximation results which we need in the following sections.
For our analysis, it is essential whether or not boundary points are included
in the nodal bases. In section~\ref{sec:Error_Lobatto},  we start by considering
Gauß-Lobatto-Legendre nodes. These include the boundary points and we demonstrate
that the error is bounded uniformly in time under some conditions on the variable
coefficient $a(x)$. Afterwards, in section~\ref{sec:Error_Gauss},
we adapt the investigation from before to Gauß-Legendre nodes which do not contain the
boundary points. We get additional error terms in our error equation and focus finally
on the different discretisations of the numerical fluxes. Similar conditions are derived
like before on $a(x)$ to guarantee that the error is bounded in time.
We confirm our investigation by numerical experiments in section~\ref{sec:Num}, which
includes also a physical interpretation of the test cases under consideration.
\revsec{Furthermore, a first analytical study about the error inequalities is given if one
of the conditions on $a$ is not fulfilled.  In section \ref{sec:Generalisation}, we generalize our
investigation to systems (linearized Euler equations and magnetic induction equation) and demonstrate problems
which arise in these cases. We give an outlook for further research.}
Finally, we summarise and discuss our results.

%% file: 2_Model_problem.tex
\section{Model Problem and Continuous Setting}\label{sec:Model_problem}

The problem under consideration is the following linear advection equation
\begin{equation}\label{eq:Model_problem}
 \begin{aligned}
   \partial_t u(t,x)+\partial_x (a(x)u(t,x))&=0, &t>0,\; x\in (x_L, x_R),\\
   u(t,x_L)&=g_L(t), & t\geq 0, \\
   u(0,x)&=u_0(x), & x\in [x_L,x_R],
 \end{aligned}
\end{equation}
with variable speed $a(x)>0$ and compatible initial and boundary conditions $u_0, \; g_L$.
Furthermore, the initial and boundary values are chosen in such a way that  $u(t,x)\in H^m(x_L,x_R) $ for
$m>1$ and that its norm $||u(t)||_{H^m}$ is bounded uniformly in time.
This condition is physically meaningful, e.g. for problems with sinusoidal boundary inputs.
However, we will also present in section \ref{sec:Num} an example where this condition is violated and our whole analysis will break down.

The impact of the boundary condition and the variable coefficient $a$ on the solution
is essential and will shortly be repeated from \cite{ranocha2018generalised,
nordstrom2017conservation}.
The energy of the solution $u$ of  the initial boundary value problem \eqref{eq:Model_problem}
is measured by the classical $\L^2$-norm
$||u||^2= \int_{x_L}^{x_R} u^2\d x$.
Focusing on the weak formulation of the advection equation \eqref{eq:Model_problem},
a test function $\phi \in C^1[x_L,x_R]$ is multiplied  and  integrated over the domain
\begin{equation}\label{eq:Continuous}
 \int_{x_L}^{x_R} (\partial_t u) \phi \d x +\int_{x_L}^{x_R} (\partial_x (au)) \phi \d x =0.
\end{equation}
Setting $\phi=u$, application of the product rule and integration-by-parts yields
\begin{align*}
\frac{\d}{\d t} ||u||^2=&2\int_{x_L}^{x_R} u \partial_t u \d x=-2 \int_{x_L}^{x_R}  u \partial_x (au) \d x\\
=&- \int_{x_L}^{x_R} \l( u \partial_x(au)+au \partial_x u +u^2 \partial_x a\right) \d x
= -au^2|_{x_L}^{x_R} - \int_{x_L}^{x_R} u^2 \partial_x a \d x\\
=&a(x_L) g_L^2-a(x_R)u(x_R)^2-\int_{x_L}^{x_R} u^2 \partial_x a \d x.
\end{align*}
Integration in time over an interval $[0,T]$ leads to
\begin{equation}\label{eq:stability}
\begin{aligned}
  ||u(T)||^2-||u_0||^2=&a(x_L) \int_0^T g_L^2(t) \d t -a(x_R) \int_0^T u^2(t,x_R) \d t \\
  &-
 \int_{x_L}^{x_R} \l(\int_0^T u^2(t,x) \d t \r) \partial_x a(x) \d x.
\end{aligned}
\end{equation}
Here, the change of energy at time $T$ can be expressed by the energy
added at the left side through the boundary
condition minus the energy lost through the right side, and an energy term
considering the variation of the coefficient $a$.
If $\partial_x a$ is bounded, the energy is also bounded \revsec{for a fixed time interval}.
It can be found in  \cite[Section 2]{nordstrom2017conservation} that the energy fulfils
\begin{equation*}
\begin{aligned}
 ||u(t)||^2 \leq &\exp \l(t||\partial_x a||_{\L^\infty} \r)
  \cdot \Bigg( ||u_0||^2 \\ &+
 \int_0^t \exp \l(-\tau||\partial_x a||_{\L^\infty} \r)
  \l( a(x_L) g_L(\tau)^2-a(x_R)u(\tau,x_R)^2\r)  \d \tau \Bigg).
\end{aligned}
\end{equation*}
The numerical scheme has to be constructed such that the approximation imitates this behaviour.
Special focus has to be given on an adequate discretisation of the flux function $f$, which depends
on the space coordinate $x$ via the variable coefficients $a(x)$.
The numerical fluxes have to be adjusted. We will specify this in section \ref{subsec:Numerical_fluxes}.

%% file: 3_Numerical_Methods.tex
\section{Flux Reconstruction with Summation-by-parts Operators and Numerical Fluxes} \label{sec:Num_methods}

In the first part of this section, we shortly repeat the main ideas of Flux Reconstruction (FR),
also known as Correction Procedure via Reconstruction, using Summation-by-parts Operators (SBP).
A more detailed introduction to this topic can be found in the articles \cite{ranocha2016summation,
ranocha2017extended} and references therein.

\subsection{Flux Reconstruction using Summation-by-parts Operators}\label{sec:SBP}

We consider a one-dimensional  scalar conservation law
\begin{equation}
\label{eq:scalar-CL}
  \partial_t u(t,x) + \partial_x f(t,x,u(t,x) ) = 0, \qquad t>0, x\in (x_0,x_K),
\end{equation}
equipped with appropriate initial and boundary conditions. The domain
 $(x_0,x_K)$ is split into $K$ non-overlapping elements $[x_0,x_K]= [ x_0,x_1]\bigcup\cdots \bigcup [x_{K-1},x_K]$.
 The FR method is a semidiscretisation applying a polynomial approximation using
 a nodal basis on each element.
 Therefore, each  interval $[x_{i-1},x_i]$ is transferred onto a standard element, which is in our case simply $[-1,1]$.
 All calculations are conducted within this reference element.
 Let $\P^N$ be the space of polynomials of degree $ \leq N$, $-1 \leq \zeta_i \leq 1$ ($i\in 0,\cdots, N)$
 the interpolation points and $\Ip^N:\L^2(-1,1)\cap C(-1,1)\to \P^N(-1,1)$ be
 the interpolation operator and $P^m_{N-1}$ be the orthogonal projection of $u$ onto $\P^{N-1}$ with respect to the inner product of the Sobolev space $H^m((-1,1))$.
 The solution is approximated by a
 polynomial $U\in \P^N$ and the basic formulation of a nodal Lagrange basis\footnote{Modal bases are also possible \cite{ranocha2017extended},
 but we won't consider these in
 this paper.} is employed.
 Instead of working with $U$ one may also express the numerical solution as the
 vector $\vec{u}$ with coefficients $\vec{u}_i = U(\zeta_i), i \in \set{0, \dots, N}$.
 All the relevant information are stored in these coefficients and
 one may write
\begin{equation}\label{eq:Approx}
 u(\zeta) \approx U(\zeta) = \sum\limits_{i=0}^N \vec{u}_i l_i(\zeta),
\end{equation}
where $l_i(\zeta)$ is the $i$-th Lagrange interpolation polynomial that satisfies $l_i(\zeta_j)=\delta_{ij}$.
In finite difference (FD) methods, it is natural to work with the coefficients only and since
we are working with SBP operators with origins lying in the FD community \cite{kreiss1974finite},
we utilise the coefficients.
Finally, the flux
$f(u)$ is also approximated by
a polynomial, where the coefficients are given by
$\vec{f}_i = f \left( \vec{u}_i \right) = f \left( U(\zeta_i) \right)$.

Now, with respect to the chosen basis (interpolation points),
(an approximation of) the derivative is represented by the
matrix $\mat{D}$. Moreover, a discrete scalar product is represented by the symmetric and positive definite
mass/norm\footnote{Both names are used. In the DG community \cite{gassner2013skew},
the matrix is called mass matrix, whereas the name norm matrix is common for FD methods.} matrix $\mat{M}$, approximating the usual $L^2$ scalar product, i.e.
\begin{equation}\label{eq:approx_matrix}
 \mat{D}\vec{u}\approx \vec{\partial_x u} \text{ and }
(\vec{u},\vec{v})_N:= \vec{u}^T\mat{M}\vec{v} \approx \int_{x_i}^{x_{i+1}} u v\d x.
\end{equation}
Using Lagrange polynomials, we get  $D_{kj}=l_j'(\zeta_k)$ and $\mat{M}=\diag{\omega_0,\dots,\omega_N}$, where $\omega_j$ are the quadrature weights associated with the nodes $\zeta_j$. For Gauß-Legendre nodes, $\omega_j=\int_{-1}^1 l_j(x)^2\d x$. For other quadrature nodes such as Gauß-Lobatto-Legendre nodes, the mass matrix is in general not exact.

SBP operators are constructed in such way that they mimic integration-by-parts on a discrete
level, as described in the review articles  \cite{svard2014review, fernandez2014review}
and references cited therein.
Until now, we have expressions/approximations for the derivative as well as for the integration. Hence, only the evaluation
on the boundary is missing.
Here, we have to introduce two different operators. First, the restriction operator, which is
represented by the matrix $\mat{R}$, approximates the interpolation of a function to the boundary
points $\{x_{i-1},x_i\}$. Second, the diagonal boundary matrix $\mat{B}=\diag{-1,1}$ gives the difference of boundary
values. It is
\begin{equation*}
 \mat{R} \vec{u}\approx \begin{pmatrix}
                         u(x_{i-1})\\
                         u(x_i)
                        \end{pmatrix} \text{ and } (u_L,u_R)\cdot \mat{B}\cdot\begin{pmatrix}
                         v_L\\
                         v_R
                        \end{pmatrix}=u_Rv_R-u_Lv_L.
\end{equation*}
Finally, all operators are introduced and they
have to fulfil the SBP property
\begin{equation}
\label{eq:SBP}
  \mat{M} \mat{D} + \mat{D}[^T] \mat{M}
  = \mat{R}[^T] \mat{B} \mat{R},
\end{equation}
in order to mimic integration-by-parts on a discrete level
\begin{equation}
  \vec{u}^T \mat{M} \mat{D} \vec{v} + \vec{u}^T \mat{D}[^T] \mat{M} \vec{v}
  \approx
  \int_{x_{i-1}}^{x_i} u \, (\partial_x v) + \int_{x_{i-1}}^{x_i} (\partial_x u) \, v
  = u \, v \big|_{x_{i-1}}^{x_i}
  \approx
  \vec{u}^T \mat{R}[^T] \mat{B} \mat{R} \vec{v}.
\end{equation}
Here, we investigate the long time error behaviour of linear problems with variable coefficients.
To represent these coefficients in our semidiscretisation, multiplication operators
are necessary. If the function $U$ is represented by $\vec{u}$, the discrete operator
approximating the linear operator $v\longmapsto vU$ is represented by the matrix $\mat{u}$,
mapping $\vec{v}$ to $\mat{u}\vec{v}$.
In a nodal basis\footnote{For a modal basis see \cite{ranocha2017extended}.},
the standard multiplication operators consider pointwise multiplication.
This means that $\mat{u}$ is diagonal with $\mat{u}=\diag{\vec{u}}$ and
$(\mat{u} \vec{v})_i = \vec{u}_i \vec{v}_i$.

One central point in our investigation in sections \ref{sec:Error_Lobatto}
and \ref{sec:Error_Gauss}
will be whether the boundary points are included in the set of interpolation nodes
(section \ref{sec:Error_Lobatto}) or not (section \ref{sec:Error_Gauss}).
This is  an essential point in this paper and also in others
\cite{ranocha2016summation, ranocha2017extended, ranocha2017shallow,
ranocha2018generalised, manzanero2017insights,
nordstrom2017conservation}.
If the boundary points $\{x_{i-1},\; x_i\}$ are included, the
restriction operators are simply
\begin{equation*}
 \mat{R}=\begin{pmatrix}
          1 &0 &\cdots &0 &0\\
          0 &0& \cdots &0 &1\\
         \end{pmatrix},\qquad
\mat{R}[^T]\mat{B}\mat{R}=\diag{-1,0,\dots,0,1}.
\end{equation*}
Thus, restriction to the boundary and multiplication commute, i.e.
\begin{equation}\label{eq:multi_commute}
\begin{pmatrix} \vec{u}_0 \, \vec{v}_0 \\ \vec{u}_N \, \vec{v}_N \end{pmatrix}
 =\l( \mat{R}\vec{u} \r)\cdot \l(\mat{R}\vec{v} \r)=\mat{R}\mat{u}\vec{v}
 =\mat{R}\mat{v}\vec{u}.
\end{equation}
At the continuous level this property is fulfilled and so we want this property
also in our semidiscretisation. However, if the boundary nodes are excluded, restriction
and multiplication will not commute in general. Therefore, some corrections have to be applied
\cite{ranocha2016summation, ranocha2017extended, ranocha2017shallow, ranocha2017comparison}.
It is common to use some linear combination/splitting of the
terms $\l(\mat{R}\vec{u}\r)\cdot\l( \mat{R}\vec{v}\r)$ and $\mat{R}\mat{v}\vec{u}$
to mimic \eqref{eq:multi_commute} at a discrete level.
We have to mention that the construction of these correction terms
can be very difficult (e.g. \cite{ranocha2017shallow} and \cite[Section~4.5]{ranocha2018thesis}) and for some equations like
Euler for example, it is still an open problem if such correction terms
exist \cite{ranocha2017comparison}.

Now, the general aspects of SBP operators are introduced
and we can focus on our FR approach.
Contrary to DG methods, we do not apply a variational formulation (i.e. weak form) of \eqref{eq:scalar-CL}.
Instead, the differential form is used, corresponding to a strong form DG method. To describe the semidiscretisation all operators
are introduced. We apply the
 discrete derivative matrix $\mat{D}$ to $\vec{f}$.
The divergence is $\mat{D} \vec{f}$.
 Since the solutions will probably have discontinuities
across elements, we will have this in the discrete flux, too.
In order to avoid this problem, a numerical flux $\vecfnum$ is introduced
which computes a common flux at the boundary using values from both neighbouring elements.
The main idea of the FR schemes is that the numerical flux at the boundaries
will be corrected by
functions in such manner that information of two neighbouring elements
interact and basic properties like conservation hold also in the  semidiscretisation.
Therefore, we add a correction term using a correction matrix $\mat{C}$ at the boundary nodes.
This gives \emph{Flux Reconstruction} its name.
Hence,  a simple FR (or correction procedure via reconstruction, CPR) method for
\eqref{eq:scalar-CL} with boundary nodes included reads
\begin{equation}
\label{eq:SBP CPR}
  \partial_t \vec{u}
  = - \mat{D} \vec{f}
    - \mat{C}\left( \vecfnum - \mat{R} \vec{f} \right).
\end{equation}
A general choice of the correction matrix $\mat{C}$ recovers the linearly stable
flux reconstruction methods of \cite{vincent2011newclass, vincent2015extended},
as described by \cite{ranocha2016summation}.
The canonical choice for the correction matrix is
\begin{equation}\label{eq:C}
  \mat{C} := \mat{M}[^{-1}] \mat{R}[^T] \mat{B}.
\end{equation}
It is a generalisation of simultaneous approximation terms (SATs) used in finite
difference methods  \cite{carpenter1999stable} and corresponds
to a strong form of the discontinuous Galerkin method \cite{kopriva2010quadrature}.
In this paper we concentrate  on the correction term using \eqref{eq:C}. However,
a generalisation to the schemes of Vincent et al. \cite{vincent2011newclass}
is possible and can be done as in \cite{ranocha2016summation, oeffner2018error}.
\revPP{However, further problems emerge concerning the interchangeability of coefficients in the broken norms and one has to be
careful.}

\subsection{Numerical Fluxes}\label{subsec:Numerical_fluxes}
Special focus has to be given on an adequate discretisation of the flux function $f$, which depends
on the space coordinate $x$ via the variable coefficients $a(x)$. The numerical fluxes
have to be adjusted. The numerical fluxes under consideration will be
\begin{align}
  \label{eq:Edge_based_cen_flux}
  \text{Edge based central flux}\qquad
  \fnum(\umin,\uadd) &= a(x)\frac{\umin+\uadd}{2},
  \\
  \label{eq:split_based_cen_flux}
  \text{Split central flux}\qquad
  \fnum(\umin,\uadd) &= \frac{\amin\umin+\aadd\uadd}{2},
  \\
  \label{eq:unsplit_based_cen_flux}
  \text{Unsplit central flux}\qquad
  \fnum(\umin,\uadd) &= \frac{(au)_-+ (au)_+}{2},
  \\
  \label{eq:Edge_based_up_flux}
  \text{Edge based upwind flux}\qquad
  \fnum(\umin,\uadd) &= a(x)\umin,
  \\
  \label{eq:split_based_up_flux}
  \text{Split upwind flux}\qquad
  \fnum(\umin,\uadd) &= \amin\umin,
  \\
  \label{eq:unsplit_based_up_flux}
  \text{Unsplit upwind flux}\qquad
  \fnum(\umin,\uadd) &= (au)_-.
\end{align}
If boundary nodes are included  and the coefficients $\vec{a}$ of the discrete version of the function $a$ are
obtained by evaluating $a$ at these nodes, \eqref{eq:Edge_based_cen_flux},
\eqref{eq:split_based_cen_flux}, \eqref{eq:unsplit_based_cen_flux}
 and \eqref{eq:Edge_based_up_flux}, \eqref{eq:split_based_up_flux}, \eqref{eq:unsplit_based_up_flux}  are identical like in the next section
 \ref{sec:Error_Lobatto}.
 From the stability analysis in \cite{ranocha2018generalised}, we know that using the unsplit fluxes \eqref{eq:unsplit_based_cen_flux},
 \eqref{eq:unsplit_based_up_flux} may result in stability issues. Furthermore, applying the other fluxes and  to guarantee stability, we
 need that the interpolation speeds have to be exact. In this cases, the edge based (\eqref{eq:Edge_based_cen_flux}, \eqref{eq:Edge_based_up_flux})  and the split numerical fluxes (\eqref{eq:Edge_based_up_flux}, \eqref{eq:split_based_up_flux})
 are equivalent. This exactness can be achieved by evaluating the speed $a(x)$ at $N+1$ Gauss-Lobatto points\footnote{We assume here a nodal basis using $N+1$ points to represent polynomials of degree $\leq N$. } and
then the unique interpolating polynomial can be evaluated at the nodes used in the basis not including the boundary.  We will consider this later in detail in section \ref{sec:Error_Gauss}.

 As it was described in section \ref{sec:SBP}, all calculations are done in a standard element $[-1,1]$.
 Therefore, a transformation of every element $e^k= [x_{k},x_{k-1}]$ to this standard element is necessary.
 Equation \eqref{eq:Continuous}  is transformed to
 \begin{equation}\label{eq:Continuous2}
  \frac{\Delta x_k}{2} \est {\partial_t u, \phi^k} + \est{\partial_{\xi} (a^k u),\phi^k}=0,
 \end{equation}
where $ \est{\cdot,\cdot}$ is the $\L^2$-scalar product, $\phi^k$ is a test function in the $k$-th element
and  the factor $\frac{\Delta x_k}{2} = \frac{x_k - x_{k-1}}{2}$ comes from the transformation.
Applying the product rule and integration-by-parts to \eqref{eq:Continuous2} yields
\begin{equation}\label{eq:Continuous3}
   \frac{\Delta x_k}{2} \est {\partial_t u, \phi^k} +
   \frac{1}{2}\l( \est{\partial_{\xi} (a^k u),\phi^k} +\est{ a^k \partial_\xi u,\phi^k}
   +\est{ u \partial_\xi a^k,\phi^k}  \r)   =0,
   \end{equation}
   \begin{equation} \label{eq:Continuous4}
    \begin{aligned}
     \frac{\Delta x_k}{2} \est {\partial_t u, \phi^k} +
   \frac{1}{2}\Bigg( a^k u \phi^k|_{-1}^1- \est{a^k u, \partial_\xi \phi^k}& \\
   +\est{ a^k \partial_\xi u,\phi^k} +\est{ u \partial_\xi a^k,\phi^k}  \Bigg)&   =0.
    \end{aligned}
   \end{equation}

Formulation \eqref{eq:Continuous4} will be used to construct the error equations.

\subsection{Numerical Errors and Approximation Results}
The error in every element is  given by
$E^k:=u^k(t,x(\xi))-U^k(t,\xi)$, where $u$ represents the solution in the $k$-th element and $U$ is the spatial approximation.
Using the interpolation operator and adding zero to the error, $E^k$ can be split in two parts:
\begin{equation}\label{eq:error_basis}
 E^k=\underbrace{(\Ip^N(u^k)-U^k)}_{=:\epsilon^k_1\in \P^N}+\underbrace{(u^k-\Ip^N(u^k))}_{=:\epsilon_p^k }.
\end{equation}
With the triangle inequality, one may bound this by
\begin{equation}\label{eq:error-esti}
 ||E^k||_N\leq ||\epsilon_1^k||_N+||\epsilon_p^k||_N,
\end{equation}
where $||\cdot||_N$ is the discrete norm induced by the discrete scalar product \eqref{eq:approx_matrix}.
$\epsilon_p^k$ is the interpolation error,
which is the sum of the series truncation error and the aliasing error.
As it was already described in
\cite{canuto2006spectral, funaro2008polynomial, hesthaven2008filtering, offner2013spectral,
offner2015zweidimensionale}, its continuous norm
converges spectrally fast for the different bases under consideration.
It is
\begin{equation*}
 |u|_{H^{m;N}(-1,1)}:= \l(\sum\limits_{j=\min(m,N+1)}^m ||u^{(j)}||_{L^2(-1,1)}^2\r)^\frac{1}{2}
\end{equation*}
the seminorm of the Sobolev space $H^m(-1,1)$. For Gauß-Lobatto-Legendre/Gauß-Legendre points,
 \begin{equation}\label{eq:Gauss-esti}
||u-\Ip(u)||_{L^2(-1,1)} \leq C N^{-m}|u|_{H^{m;N}(-1,1)},
 \end{equation}
where $C$ depends on $m$.
In view of our investigation, one needs to consider the
interpolation error not only in the standard interval $[-1,1]$, but in each element $e^k$.
Therefore, the estimation \eqref{eq:Gauss-esti} will be transform  to every element\footnote{A more detailed analysis can be found in
 \cite{bernardi1989properties,bernardi1992polynomial}.}.
 With a combination of \cite[Theorem 6.6.1]{funaro2008polynomial} and
 \cite[Section 5.4.4]{canuto2006spectral}, for Gauß-Legendre nodes
 \begin{equation}\label{eq:Gauss-esti_element}
||\epsilon^k_p||_{H^n(e^k)} \leq  C  \l(\Delta x_k\r)^{n-\min\{m,N\}+
\frac{1}{2}} N^{n-m+\frac{1}{2}}|u|_{H^{m;N}(e^k)}
 \end{equation}
 for $n=0,1$. For Gauss-Lobatto-Legendre nodes, delete $\frac{1}{2}$ on the right side of
 \eqref{eq:Gauss-esti_element}.
A finite dimensional normed vector space is considered and
all norms are equivalent there.  This allows to bound
the discrete norm in terms of the continuous ones and implies that $||\epsilon_p^k||_N$
in \eqref{eq:error-esti} decays spectrally fast in all cases of consideration.
In other words, $\epsilon_1^k$ has to be investigated in detail.
This error describes the difference of the interpolation of $u$ and
the spatial approximation $U$.

Therefore, we have to consider the numerical schemes under consideration.
The semidiscretisation of \eqref{eq:Model_problem} is given by the following form\revPP{:}
\begin{equation}\label{eq:semi_disc_Gauss}
\begin{aligned}
  \partial_t \vec{u} =& -\frac{1}{2} \mat{D} \mat{a} \vec{u}- \frac{1}{2}\mat{a} \mat{D} \vec{u}
 -\frac{1}{2}\mat{u} \mat{D} \vec{a} \\
 &- \mat{M}[^{-1}]\mat{R}[^T] \mat{B}
 \l( \vecfnum-\frac{1}{2}\mat{R} \mat{a}\vec{u}-\frac{1}{2}
 \l(\mat{R} \vec{a}\r) \cdot \l(\mat{R} \vec{u}\r)  \r),
\end{aligned}
\end{equation}
where analogously to the continuous setting a split formulation has been applied.
The last term is due to the fact that for Gauss-Legendre nodes the restriction operators $\mat{R}$ do not commute with the multiplication operators.
Therefore, corrections have to be used. If boundary nodes are included, multiplication and restriction commute and we can simplify \eqref{eq:semi_disc_Gauss}
to
\begin{equation}\label{eq:semi_disc}
 \partial_t \vec{u} +\frac{1}{2} \l(\mat{D} \mat{a} \vec{u}+ \mat{a} \mat{D} \vec{u}
 +\mat{u} \mat{D} \vec{a}\r)  + \mat{M}[^{-1}]\mat{R}[^T] \mat{B}\l( \vecfnum-\mat{R} \mat{a}\vec{u} \r)=0.
\end{equation}
In \eqref{eq:semi_disc_Gauss} and  \eqref{eq:semi_disc},
the terms 2 -- 4 approximate the split form \\
$\frac{1}{2} \l(\partial_x (au)+a (\partial_x u)+u(\partial_x a \r) $ of the flux derivative $\partial_x (au)$
of \eqref{eq:Model_problem}.
Since the semidiscretisation  is used in every element $e^k$, one obtains for every element the following form:
\begin{equation}\label{eq:semi_disc_el}
\begin{aligned}
  \frac{\Delta x_k}{2}
 \partial_t \vec{u} +\frac{1}{2} \l(\mat{D} \mat{a} \vec{u}+ \mat{a} \mat{D} \vec{u}
 +\mat{u} \mat{D} \vec{a}\r)&  \\
 + \mat{M}[^{-1}]\mat{R}[^T] \mat{B} \l( \vecfnum-\frac{1}{2}\mat{R} \mat{a}\vec{u}-\frac{1}{2} \l(\mat{R} \vec{a}\r)
 \cdot \l(\mat{R} \vec{u}\r) \r)&=0.
\end{aligned}
\end{equation}
Using a Galerkin approach, $\vec{\phi}^{k,T} \mat{M}$
is multiplied to \eqref{eq:semi_disc_el}, resulting due to the
SBP property \eqref{eq:SBP} in
\begin{equation}\label{eq:Semi_disc_gal}
\begin{aligned}
 \frac{\Delta x_k}{2}\vec{\phi}^{k,T}\mat{M}
 \partial_t \vec{u} +\frac{1}{2} \vec{\phi}^{k,T} \mat{M}\l(\mat{D} \mat{a} \vec{u}+ \mat{a} \mat{D} \vec{u}
 +\mat{u} \mat{D} \vec{a}\r)& \\
 + \vec{\phi}^{k,T}\mat{R}[^T] \mat{B}
 \left(  \vecfnum-\frac{1}{2}\mat{R} \mat{a}\vec{u}-\frac{1}{2} \l(\mat{R} \vec{a}\r)
 \cdot \l(\mat{R} \vec{u}\r)  \right)& =0,\\
 \frac{\Delta x_k}{2}\vec{\phi}^{k,T}\mat{M}
 \partial_t \vec{u} +\frac{1}{2} \vec{\phi}^{k,T}\l(\mat{R}[^T] \mat{B} \mat{R} - \mat{D}[^T] \mat{M}
  \r)\mat{a} \vec{u}
  + \frac{1}{2} \vec{\phi}^{k,T} \mat{M} \mat{a} \mat{D} \vec{u} & \\
 + \frac{1}{2} \vec{\phi}^{k,T} \mat{M} \mat{u} \mat{D} \vec{a}+
 \vec{\phi}^{k,T}\mat{R}[^T] \mat{B}\left(  \vecfnum-\frac{1}{2}\mat{R}
 \mat{a}\vec{u}-\frac{1}{2} \l(\mat{R} \vec{a}\r) \cdot \l(\mat{R} \vec{u}\r)  \right)&=0.
\end{aligned}
\end{equation}
The diagonal multiplication operators are self-adjoint with respect to $\mat{M}$, i.e. $\mat{M} \mat{a} =\mat{a} \mat{M}$,
and $\mat{M} \mat{u} =\mat{u}\mat{M}$. Thus, \eqref{eq:Semi_disc_gal}  is
\begin{multline}\label{eq:semi_disc_ver_Gauss}
   \frac{\Delta x_k}{2}\vec{\phi}^{k,T} \mat{M} \partial_t \vec{u}-\frac{1}{2}
   \vec{\phi}^{k,T} \mat{D}[^T] \mat{M}\mat{a} \vec{u}+
 \frac{1}{2} \vec{\phi}^{k,T}  \mat{a} \mat{M} \mat{D} \vec{u}+ \frac{1}{2}
 \vec{\phi}^{k,T}  \mat{u} \mat{M} \mat{D} \vec{a}\\
 + \vec{\phi}^{k,T} \mat{R}[^T]\mat{B} \left(  \vecfnum-\frac{1}{2} \l(\mat{R}
 \vec{a}\r) \cdot \l(\mat{R} \vec{u}\r)  \right)=0\revPP{,}
\end{multline}
or with boundary nodes included
\begin{equation}\label{eq:semi_disc_ver}
\begin{aligned}
   \frac{\Delta x_k}{2}\vec{\phi}^{k,T}\mat{M}
 \partial_t \vec{u} - \frac{1}{2} \vec{\phi}^{k,T}\mat{D}[^T] \mat{M}
  \mat{a} \vec{u}+ \frac{1}{2} \vec{\phi}^{k,T} \mat{a} \mat{M}  \mat{D} \vec{u}& \\
 + \frac{1}{2} \vec{\phi}^{k,T}  \mat{u} \mat{M}\mat{D} \vec{a}  +
 \vec{\phi}^{k,T}\mat{R}[^T] \mat{B} \l( \vecfnum- \frac{1}{2}\mat{R} \mat{a}\vec{u} \r)&=0.
\end{aligned}
\end{equation}
The error equations will be derived using both semidiscretisations. Before starting
with the Gauß-Lobatto-Legendre case in the next section \ref{sec:Error_Lobatto}, we
shortly repeat for clarification again the notation which will be used in this paper
in Table~\ref{tab:notation}.

\renewcommand{\arraystretch}{1.2}
\begin{table}[!ht]
\caption{Summary of the notations used in this article.}
\label{tab:notation}
\centering
\begin{tabular}{l|l}
  \multicolumn{1}{c}{\textbf{\textit{Notation}}} & \textbf{\textit{Interpretation}}\\ \hline
  $u$ & is the solution of \eqref{eq:Model_problem}.\\
  $U$ & is the spatial approximation of $u$ given by \eqref{eq:Approx}. \\
  $\vec{u}$ & \multirow{2}{*}{\shortstack[l]{are the coefficients of $U$, evaluated at the \\ \quad interpolation/quadrature nodes.}}\\
  & \\
  $\mat{D}$ & is the discrete derivative matrix. \\
  $\mat{R}$ & \multirow{2}{*}{\shortstack[l]{is the restriction operator performing interpolation \\ \quad to the boundary.}}\\
  & \\
  $\mat{M}$ & is the mass/norm matrix.\\
  $\est{\cdot,\cdot}$ & is the usual $L^2$ scalar product.\\
  $||\cdot|| $ & is the norm induced by the $L^2$ scalar product.\\
  $(\cdot,\cdot)_N$ & is the discrete scalar product given by \eqref{eq:approx_matrix}.\\
  $||\cdot||_N$ & is the norm induced by the discrete scalar product from above.\\
  $\Ip^N$  &  is the interpolation operator. \\
  $P^m_{N-1}(u)$ & \multirow{2}{*}{\shortstack[l]{is the orthogonal projection of $u$ onto $\P^{N-1}(-1,1)$ using \\ \quad the inner product of $H^m(-1,1)$ (\revPP{See also appendix \ref{sec:Appendix}}). } } \\
  & \\
  $E^k= u^k-U^k$ & is the total error in the $k$-th element.\\
  $\epsilon_1^k:=\Ip^N(u^k)-U^k$ & \multirow{2}{*}{\shortstack[l]{is the difference between interpolation and spatial \\ \quad approximation in the $k$-th element.}} \\
  & \\
  $\epsilon_p^k=u^k-\Ip^N(u^k)$ & is the interpolation error which decays spectrally fast.
\end{tabular}
\end{table}
\renewcommand{\arraystretch}{1}

%% file: 4_Error_bound_semidiscrete_Lobatto.tex
\section{Error Behaviour using Gauß-Lobatto Nodes}\label{sec:Error_Lobatto}

In this section, Gauß-Lobatto-Legendre nodes will be used in the discretisation,
resulting in diagonal norm SBP operators including the boundary nodes. In this
case, multiplication and restriction to the boundary commute and the interpolated
speed $a(x)$ is automatically continuous.
\revsec{Before starting our investigation in this section,
we will briefly summarize our final results for both cases (Gau\ss-Lobatto-Legendre and Gau\ss-Legendre).}

\begin{Result}
\revsec{$\eta(t)$ is a factor which depends on $\epsilon_1$, the values of $a$
and $a'$. If there exits a positive constant $\delta$, such that the mean value of
 $\eta(t)$ can be bounded from below, then there exists a constant $C$ such
that the errors $\epsilon_1^k(t)$ of \eqref{eq:error_basis} satisfy the inequality
\begin{equation*}
||\epsilon_1(t)||_N\leq \frac{1-\exp(-\delta t)}{\delta} C,
\end{equation*}
in the discrete norm $||\cdot||_N$. The total error is bounded in time.}
\end{Result}

\revsec{
In the following, we will derive the exact conditions when the above inequality
is fulfilled and specify in detail what factors play a key role in the definition of $\eta$ and $\delta$.}
\revPP{
We  outline the steps of our analysis.
All steps of the investigation in sections \ref{sec:Error_Lobatto}
and \ref{sec:Error_Gauss} are almost analogous except
that in step 5 we have to consider the different flux functions
\eqref{eq:Edge_based_cen_flux} - \eqref{eq:unsplit_based_up_flux}
in our investigation\footnote{We have an additional error term in section \ref{sec:Error_Gauss}, but this does not change
the major steps of the study.}.
The main stepts are the following:
\begin{itemize}
 \item[1.] We derive an error equation for $\epsilon_1^k$ of \eqref{eq:error_basis} by inserting the error $E^k$ into
 the continuous equation \eqref{eq:Continuous4} for every element.
 \item[2.] By adding zero in a suitable way, we are able to split the equations into
 a continuous and a discrete part.
 \item[3.] We add both parts for every element and get the error behaviour for the total domain.
 \item[4.] We estimate the continuous terms and get an inequality for the error $\epsilon_1$ in the discrete norms.
 \item[5.] We split the terms with the numerical fluxes. In the Gauss-Legendre case (section \ref{sec:Error_Gauss}),
 we have to be careful with respect to
 the used implementation of the numerical fluxes.
 \item[6.] We estimate the long time error behaviour under some assumptions.
\end{itemize}
}
In the following, the error equation for $\epsilon_1^k=\Ip^N(u^k)-U^k$ will be
derived. Starting by considering Gauß-Lobatto-Legendre nodes in our
semidiscretisation and putting
 $u=\Ip^N(u^k)+\epsilon_p^k$ into \eqref{eq:Continuous4} yields
\begin{align*}
  &\frac{\Delta x_k}{2} \est {\partial_t \Ip^N(u^k), \phi^k} +
   \frac{1}{2}\Bigg( a^k \Ip^N(u^k) \phi^k|_{-1}^1- \est{a^k \Ip^N(u^k), \partial_\xi \phi^k} \\
   &+
   \est{ a^k \partial_\xi \Ip^N(u^k),\phi^k}
   +\est{ \Ip^N(u^k) \partial_\xi a^k,\phi^k}  \Bigg)   \\
   =&- \frac{\Delta x_k}{2} \est {\partial_t \epsilon_p^k, \phi^k} +
   \frac{1}{2} \est{a^k \epsilon_p^k, \partial_\xi \phi^k} -\frac{1}{2}
   \est{ a^k \partial_\xi \epsilon_p^k,\phi^k}
   -\frac{1}{2}\est{ \epsilon_p^k \partial_\xi a^k,\phi^k},
\end{align*}
where $\phi^k \in \P^N$ is a polynomial test function.
For Gauß-Lobatto-Legendre nodes, $a^k\epsilon_p^k=0$ at the endpoints,  since the interpolant
is equal to the solution there. Thus, $a^k\epsilon_p^k\phi^k\Big|_{-1}^1=0$.
Using integration-by-parts for  $\est{ a^k\epsilon_p^k, \partial_\xi \phi^k }$ yields
\begin{equation}\label{eq:continuous_error}
 \begin{aligned}
   &\frac{\Delta x_k}{2} \est {\partial_t \Ip^N(u^k), \phi^k} +
   \frac{1}{2}\Bigg( a^k \Ip^N(u^k) \phi^k|_{-1}^1- \est{a^k \Ip^N(u^k), \partial_\xi \phi^k} \\
   &+
   \est{ a^k \partial_\xi \Ip^N(u^k),\phi^k}
   +\est{ \Ip^N(u^k) \partial_\xi a^k,\phi^k}  \Bigg)   \\
   =&- \frac{\Delta x_k}{2} \est {\partial_t \epsilon_p^k, \phi^k} -
   \frac{1}{2} \est{\partial_\xi(a^k \epsilon_p^k),  \phi^k} -\frac{1}{2}
   \est{ a^k \partial_\xi \epsilon_p^k,\phi^k}
   -\frac{1}{2}\est{ \epsilon_p^k \partial_\xi a^k,\phi^k}.
 \end{aligned}
\end{equation}
\revPP{We have to transfer the continuous scalar product from \eqref{eq:continuous_error}
to the discrete ones. Therefore, we are following the ideas from \cite{oeffner2018error},
add zero to the above equation and rearrange these terms.
We will explain this for the first term on the left side of \eqref{eq:continuous_error}
in detail.
The third to fifth terms on the left side are handled analogously and details can be found in the
appendix \ref{sec:Appendix}.
Applying the interpolation operator together with discrete norms results in
\begin{align}
\est{\partial_t \Ip^N(u^k),\phi^k  }=&
\l( \partial_t \vec{\Ip^N(u^k)}, \vec{\phi}^k \r)_N \nonumber \\ &
+\l\{ \est{ \partial_t \Ip^N(u^k), \phi^k }-
  \l( \partial_t \vec{\Ip^N(u^k)}, \vec{\phi}^k \r)_N  \r\}, \label{eq:interpolation}
\end{align}
Now, we are introducing in the factor $Q$ in the above equation which measures
the projection error of a polynomial of degree $N$ to a polynomial of degree $N-1$.
We can rewrite \eqref{eq:interpolation} as
\begin{align*}
 \est{\partial_t \Ip^N(u^k),\phi^k  }=
\l(\vec{\Ip^N(u^k)}, \vec{\phi}^k \r)_N
+\l\{ \est{ Q(u^k), \phi^k }-
  \l(  \vec{Q(u^k)}, \vec{\phi}^k \r)_N  \r\}
\end{align*}
where $Q(u^k):= \partial _t \l( \Ip^N(u^k)-P^m_{N-1} \l( \Ip^N(u^k)\r)\r)$ and
$P^m_{N-1}$ is the orthogonal projection\footnote{More
details can be found in the appendix \eqref{sec:Appendix}.}
of $u$ onto $\P^{N-1}$ using the inner product of $H^m(e^k)$.
We get similar factors ($Q_1-Q_3$) for the other three terms.
Since $u$ and $a$ are bounded, also all of theses values have to be bounded.
}
Finally, the values of the interpolation polynomial at the boundaries of the element ($-1$ and $1$)
can be approximated by a limitation process from
the left side $\Ip^N(u^k)^{-}$ and right side $\Ip^N(u^k)^{+}$.
To simplify the notation, let
\begin{multline}
 \vecfnumk \l(\Ip^N(u^k)^{-},\Ip^N(u^k)^{+}\r)
 \\:=
 \l(\fnum \l(\Ip_R^N(u^{k-1}),\Ip_L^N(u^{k})\r), \fnum \l(\Ip_R^N(u^{k}),\Ip_L^N(u^{k+1})\r)      \r)^T.
\end{multline}
\revPP{ For boundary points included, the interpolation is continuous
(because the exact solution $u$ is continuous)
and all numerical fluxes are exactly the products of the interpolation and the
coefficient values. One obtains
\begin{equation}\label{eq:in_flux}
 \begin{aligned}
   \frac{1}{2} a^k\Ip^N(u^k)\phi^k\bigg|_{-1}^1=& \vec{\phi}^{k,T} \mat{R}[^T] \mat{B}
 \l(\vecfnumk \l(\Ip^N(u^k)^{-},\Ip^N(u^k)^{+}\r)-\frac{1}{2 }\mat{R}\mat{a^k}\vec{u} \r)
 .
 \end{aligned}
\end{equation}}

\revPP{Using the above investigation} and putting \eqref{eq:interpolation}--\eqref{eq:in_flux} in
\eqref{eq:continuous_error} results in
\begin{equation}\label{eq:Delta}
\begin{aligned}
 &\quad
 \frac{\Delta x_k}{2} \l( \partial_t \vec{\Ip^N(u^k)}, \vec{\phi}^k \r)_N
 + \vec{\phi}^{k,T} \mat{R}[^T] \mat{B}
 \l(\vecfnumk \l(\Ip^N(u^k)^{-},\Ip^N(u^k)^{+}\r)-\frac{1}{2 }\mat{R}\mat{a^k}\vec{u} \r)
 \\&\quad
 -\frac{1}{2} \l(  \mat{a^k} \vec{\Ip^N(u^k)} , \partial_\xi \vec{\phi}^k\r)_N
 +\frac{1}{2} \l(  \mat{a^k}  \partial_\xi \vec{\Ip^N(u^k)}, \vec{\phi}^k   \r)_N
 +\frac{1}{2} \l(  \mat{\Ip^N(u^k)} \partial_\xi \vec{a}^k,  \vec{\phi}^k   \r)_N
 \\&=
 +\frac{\Delta x_k}{2} \est{T^k(u), \phi^k}
 +\frac{\Delta x_k}{4} \Bigg( (\vec{Q(u^k)}, \vec{\phi}^k)_N
 -(\vec{Q_1(u^k)}, \partial_x \vec{\phi}^k)_N  \\&\quad
 +(\mat{a^k}\vec{Q_2(u^k)}, \vec{\phi}^k)_N
 +(\mat{Q_3(u^k)}\partial_x \vec{a}^k,  \vec{\phi^k}  )_N   \Bigg)
 +\frac{\Delta x_k}{4}  \est{Q_1(u^k), \partial_x \phi^k},
\end{aligned}
\end{equation}
with
\begin{equation}\label{Def_T}
\begin{aligned}
  T^k(u):= &- \Bigg\{ \partial_t \epsilon_p^k+
 \frac{1}{2} \l( \partial_x(a^k \epsilon_p^k)+\epsilon_p^k \partial_x a^k +a^k \partial_x \epsilon_p^k \r)\\
 &+ \frac{1}{2}
 \l( Q(u^k)+ a^k Q_2(u^k) +(Q_3(u^k) \partial_x a^k) \r) \Bigg\}.
\end{aligned}
\end{equation}
Here, in definition \eqref{Def_T} we have again the derivatives in $x$ since we make the term independent
from the transformation. Therefore, we have in \eqref{eq:Delta} a $\frac{\Delta x_k}{2}$
in the $T^k$ terms.

By \eqref{eq:Gauss-esti_element}, the interpolation error $\epsilon_p^k$ converges in $N$ to zero,
if $m>1$ and the Sobolev norm of the solution is uniformly bounded in time\footnote{Therefore, we need the initial and boundary conditions in the
model problem \eqref{eq:Model_problem}. }.
Equation \eqref{eq:semi_disc_ver} is subtracted form \eqref{eq:Delta}
and with $\epsilon_1^k =\Ip^N(u^k)-U^k$ one
obtains
\begin{align*}\label{eq:Delta2}
& \frac{\Delta x_k}{2} \l( \partial_t \vec{\epsilon}_1^k, \vec{\phi}^k \r)_N
  + \vec{\phi}^{k,T} \mat{R}[^T] \mat{B}
  \l(\vecfnumk \l( (\epsilon_1^k )^{-},(\epsilon_1^k )^{+}\r)-\frac{1}{2 }\mat{R}\mat{a^k}\vec{\epsilon}_1^k
  \r) \\
  -&\frac{1}{2} \l(  \mat{a^k} \vec{\epsilon}_1^k , \partial_\xi \vec{\phi}^k\r)_N
  +\frac{1}{2} \l(  \mat{a^k}  \partial_\xi \vec{\epsilon}_1^k, \vec{\phi}^k   \r)_N
  +\frac{1}{2} \l(  \mat{\epsilon}_1^k \partial_\xi \vec{a}^k,  \vec{\phi}^k   \r)_N \\
  =&+\frac{\Delta x_k}{2} \est{T^k(u), \phi^k}
   +\frac{\Delta x_k}{4}  \est{Q_1(u^k), \partial_x \phi^k}
  +\frac{\Delta x_k}{4} \Bigg( (\vec{Q(u^k)}, \vec{\phi}^k)_N \\
  &-(\vec{Q_1(u^k)}, \partial_x \vec{\phi}^k)_N
  +(\mat{a^k}\vec{Q_2(u^k)}, \vec{\phi}^k)_N
  +\l(\mat{Q_3(u^k)}\partial_x \vec{a}^k,  \vec{\phi^k}   \right)_N
\Bigg) .
\end{align*}
Putting $\phi^k=\epsilon_1^k$ results in the energy equation
\begin{equation}\label{eq:energy1}
\begin{aligned}
 &\quad \frac{\Delta x_k}{4} \frac{\d}{\d t} ||\epsilon_1^k||_N^2+
 \vec{\epsilon}_1^{k,T} \mat{R}[^T] \mat{B}
 \l(\vecfnumk \l( (\epsilon_1^k )^{-},(\epsilon_1^k )^{+}\r)-\frac{1}{2 }\mat{R}\mat{a^k}\vec{\epsilon}_1^k \r)
  \\&\quad-\frac{1}{2} \l(  \mat{a^k} \vec{\epsilon}_1^k , \partial_\xi \vec{\epsilon}_1^k \r)_N
 +\frac{1}{2} \l(  \mat{a^k}  \partial_\xi \vec{\epsilon}_1^k, \vec{\epsilon}_1^k  \r)_N
 +\frac{1}{2} \l(  \mat{\epsilon}_1^k \partial_x \vec{a}^k,  \vec{\epsilon}_1^k   \r)_N
 \\&=
 +\frac{\Delta x_k}{2} \est{T^k(u), \epsilon_1^k}
  +\frac{\Delta x_k}{4}  \est{Q_1(u^k), \partial_x \epsilon_1^k}
  \\&\quad
 +\frac{\Delta x_k}{4} \Bigg( (\vec{Q(u^k)}, \vec{\epsilon}_1^k)_N
 -(\vec{Q_1(u^k)}, \partial_x \vec{\epsilon}_1^k)_N
 +(\mat{a^k}\vec{Q_2(u^k)}, \vec{\epsilon}_1^k)_N \\&\quad  +\l(\mat{Q_3(u^k)}\partial_x \vec{a}^k,  \vec{\epsilon}_1^k   \right)_N \Bigg).
\end{aligned}
\end{equation}

\revPP{Since $\mat{M}[^T]=\mat{M}$, we get
\begin{equation}\label{eq:Deleting_a_terms}
\begin{split}
 &\frac{1}{2} \l(\mat{a^k} \vec{\epsilon}_1^k, \partial_\xi \epsilon_1^k \r)_N= \frac{1}{2}\vec{\epsilon}_1^{k,T} \mat{a^k}[^T] \mat{M}\mat{D}\epsilon_1^k,\\
  &\frac{1}{2} \l( \partial_\xi \epsilon_1^k, \mat{a^k} \vec{\epsilon}_1^k\r)_N=
  \frac{1}{2} \vec{\epsilon}_1^{k,T} \mat{D}[^T] \mat{M} \mat{a^k} \epsilon_1^k=\frac{1}{2}\vec{\epsilon}_1^{k,T} \mat{a}^{k,T} \mat{M}\mat{D}\epsilon_1^k,
\end{split}
\end{equation}
and}
one obtains in \eqref{eq:energy1}
\begin{equation}\label{eq:energy2}
\begin{aligned}
 &\quad \frac{\Delta x_k}{4}   \frac{\d}{\d t} ||\epsilon_1^k||_N^2+
 \vec{\epsilon}_1^{k,T} \mat{R}[^T] \mat{B}
 \l(\vecfnumk \l( (\epsilon_1^k )^{-},(\epsilon_1^k )^{+}\r)-\frac{1}{2 }\mat{R}\mat{a^k}\vec{\epsilon}_1^k \r)
 \\& +\frac{\Delta x_k}{4} \l(  \mat{\epsilon}_1^k  \vec{\epsilon}_1^k  ,  \partial_x \vec{a}^k \r)_N=
 \frac{\Delta x_k}{2} \est{T^k(u), \epsilon_1^k}
  +\frac{\Delta x_k}{4}  \est{Q_1(u^k), \partial_x \epsilon_1^k}
 \\&
 +\frac{\Delta x_k}{4} \Bigg( (\vec{Q(u^k)}, \vec{\epsilon}_1^k)_N
 -(\vec{Q_1(u^k)}, \partial_x \vec{\epsilon}_1^k)_N
 +(\mat{a^k}\vec{Q_2(u^k)}, \vec{\epsilon}_1^k)_N \\ &
 +\l(\mat{Q_3(u^k)}\partial_x \vec{a}^k,  \vec{\epsilon}_1^k \right)_N \Bigg).
 \end{aligned}
\end{equation}
Summing this up over all elements
\revPP{ and by defining the numerical flux of the error as
$\vecepsilon:=\vecfnumk \l( \l( \epsilon_1^k\r)^-,\l(\epsilon_1^k\r)^+\r)$, the global
energy of the error is
}
\begin{equation}\label{eq:energy3}
\begin{aligned}
   &\quad \frac{1}{2}   \frac{\d}{\d t} \sum_{k=1}^K \frac{\Delta x_k}{2}||\epsilon_1^k||_N^2+\sum_{k=1}^K
 \vec{\epsilon}_1^{k,T} \mat{R}[^T] \mat{B}
\l( \vecepsilon-\frac{1}{2 }\mat{R}\mat{a^k}\vec{\epsilon}_1^k \r)
\\&\qquad
 +\frac{1}{2} \sum_{k=1}^K \frac{\Delta x_k}{2} \l(  \mat{\epsilon}_1^k  \vec{\epsilon}_1^k  ,  \partial_x \vec{a}^k \r)_N
  \\&=
  \sum_{k=1}^K\frac{\Delta x_k}{2} \est{T^k(u), \epsilon_1^k}
  + \sum_{k=1}^K\frac{\Delta x_k}{4}  \est{Q_1(u^k), \partial_x \epsilon_1^k} - \sum_{k=1}^K\frac{\Delta x_k}{4}
  (\vec{Q_1(u^k)}, \partial_x \vec{\epsilon}_1^k)_N
  \\&\quad
 +  \sum_{k=1}^K\frac{\Delta x_k}{4} \left( (\vec{Q(u^k)}, \vec{\epsilon}_1^k)_N
 +(\mat{a^k}\vec{Q_2(u^k)}, \vec{\epsilon}_1^k)_N
 +\l(\mat{Q_3(u^k)}\partial_x \vec{a}^k,  \vec{\epsilon}_1^k \right)_N      \right).
 \end{aligned}
\end{equation}
The right-hand side of \eqref{eq:energy3} will be estimated using the Cauchy-Schwarz
inequality. \revsec{For example (the others terms are handled similarly),
\begin{align}
  \sum_{k=1}^K\frac{\Delta x_k}{2} \est{T^k(u), \epsilon_1^k} &\leq
  \sqrt{ \sum_{k=1}^K\frac{\Delta x_k}{2} ||T^k(u)||^2 }  \sqrt{ \sum_{k=1}^K\frac{\Delta x_k}{2} ||\epsilon_1^k||^2   } ,
  \label{eq:T_est}\\
   \sum_{k=1}^K\frac{\Delta x_k}{4} \l(\vec{Q_1(u^k)}, \partial_x \vec{\epsilon}_1^k\r)_N
   &\leq  \frac{1}{2} \sqrt{ \sum_{k=1}^K\frac{\Delta x_k}{2} ||\vec{Q_1(u^k)}||_N^2 }
   \sqrt{ \sum_{k=1}^K\frac{\Delta x_k}{2} ||\partial_x \vec{\epsilon}_1^k||_N^2     } , \label{Eq:est_Q_1}
\end{align}}%
Using an estimation for the differential operator $\partial_x$ and the fact that $\epsilon_1\in \P^N$, it is
$||\partial_x \vec{\epsilon}_1^k||_N^2 \leq c_1N^2   ||\vec{\epsilon}_1^k||_N^2$ with a positive constant  $c_1$.
This is due to the fact that all norms are equivalent and we can estimate with a
Markov-Bernstein type inequality,
see \cite{govil1999markov}.
 The estimation is used \revsec{for example} in \eqref{Eq:est_Q_1}. 
 An alternative approach would have been
 to use the summation-by-parts property \eqref{eq:SBP} and estimate analogously.

With the global norm over all elements and the equivalence between the continuous and discrete norms, we obtain
\begin{equation}\label{eq:R_absch}
\begin{aligned}
&\frac{1}{2}   \frac{\d}{\d t} \sum_{k=1}^K \frac{\Delta x_k}{2}||\epsilon_1^k||_N^2+\sum_{k=1}^K
 \vec{\epsilon}_1^{k,T} \mat{R}[^T] \mat{B}
\l( \vecepsilon-\frac{1}{2 }\mat{R}\mat{a^k}\vec{\epsilon}_1^k \r)
 \\ & +\frac{1}{2} \sum_{k=1}^K \frac{\Delta x_k}{2} \l(  \mat{\epsilon}_1^k  \vec{\epsilon}_1^k  ,
 \partial_x \vec{a}^k \r)_N
 \leq  \Bigg\{ c||T ||+\frac{cN}{2}||Q_1|| +
 \frac{1}{2}\Bigg( ||Q||_N+ N\tilde{c_1}||Q_1||_N \\
 &+||a Q_2||_N +||Q_3\partial_x a||_N \Bigg)   \Bigg\}
 ||\epsilon_1||_N
 \equiv \hat{\Ep}(t,N) ||\epsilon_1||_N
 \end{aligned}
\end{equation}
Applying the same approach like in \cite{kopriva2017error} and splitting the sum
into three parts (one for the left physical boundary, one for the right physical
boundary and a sum over the internal element endpoints), it is
\begin{align*}
 &\quad \sum\limits_{k=1}^K \vec{ \epsilon}_1^{k,T}
  \mat{R}[^T] \mat{B} \l(\vecepsilon -\frac{1}{2} \mat{R}\mat{a^k}\vec{\epsilon_1}^k \r)
  \\&=
 \sum\limits_{k=1}^K \vec{ \epsilon}_1^{k,T}
  \mat{R}[^T] \mat{B} \Bigg(\vecfnumk \l( \l( \epsilon_1^k\r)^-,\l(\epsilon_1^k\r)^+\r) -
  \frac{1}{2} \mat{R}\mat{a^k} \vec{\epsilon}_1^k \Bigg)
 \\&=
 - \Eps_L^1 \l( f^{\mathrm{num},1}_L -\frac{1}{2} a_L^1\Eps_L^1 \r)+
 \sum\limits_{k=2}^K \l( f^{\mathrm{num},k}_L -\frac{1}{2} a^{k-1}_R\left(\Eps_R^{k-1}+\Eps_L^{k}  \right) \r)
 \\&\qquad \l( \Eps_R^{k-1}-\Eps_L^{k}  \r)
 + \Eps_R^K \l(f^{\mathrm{num},K}_R -\frac{1}{2} a^{K}_R \Eps^K_R \r).
\end{align*}
 Here,
 $\Eps^k_i$ ($i=L,R; \; k=1,\dots, K$) represents the error $\epsilon_1^k$ at the
 the position in the elements, and (to shorten the notation)
 $f^{\mathrm{num},k}_L:=f^{\mathrm{num},k} \l(\Eps^{k-1}_R,\Eps^{k}_L \r)$, $ f^{\mathrm{num},1}_L:=f^{\mathrm{num},1} \l(0,\Eps^{1}_L \r)$
 and $f^{\mathrm{num},K}_R:= f^{\mathrm{num},1} \l(\Eps^K_R,0 \r) $.
 The external states for the physical boundary contributions are zero,
 because $\Ip^N(u^1)=g$ at the left boundary and the external state for $U^1$ is set to $g$.
 At the right boundary, where the upwind numerical flux is used,
 it doesn't matter what the external state is,
 since its coefficient in the numerical flux is zero.
 One gets for the inner element with $\jump{\Eps^k}=\Eps_R^{k-1}-\Eps_L^{k} $
 \begin{equation*}
 \begin{aligned}
  \sum\limits_{k=2}^K \l( f^{\mathrm{num},k}_L -\frac{1}{2} a^{k-1}_R\left(\Eps_R^{k-1}+\Eps_L^{k}  \right)
  \r)
  \l( \Eps_R^{k-1}-\Eps_L^{k}  \r)& = \sum\limits_{k=2}^K \frac{\sigma a^{k-1}_R}{2}
  \l(\jump{\Eps^{k} }\r)^2
  \geq0,\quad \\
 &\text{with}
 \begin{cases}
 \sigma=0 \quad \text{central flux, }\\
 \sigma=1 \quad \text{ upwind flux. }  \end{cases}
 \end{aligned}
 \end{equation*}
\revsec{For the left and right boundaries, it is finally
\begin{align*}
 &\text{left:} &-\Eps_L^1 \l( f^{\mathrm{num},1}_L -
 \frac{1}{2}  a_L^1\Eps_L^1 \r)=
 \frac{\sigma a_L^1}{2}\l(\Eps_L^1\r)^2,
 \\
 &\text{right:}&\Eps_R^K \l(f^{\mathrm{num},K}_R -
 \frac{1}{2} a_R^K\Eps^K_R \r)=
 \frac{\sigma  a_R^K}{2} \l(\Eps^K_R  \r)^2.
\end{align*}}
Therefore, the energy growth rate is bounded by
\begin{equation}\label{eq:Errorbound}
\begin{aligned}
   &\frac{1}{2} \frac{\d}{\d t} ||\epsilon_1||^2_N +
  \underbrace{\frac{\sigma}{2} \l( a^K_R\l(\Eps^K_R  \r)^2 +a_L^1\l(\Eps_L^1\r)^2  \r) +
  \frac{\sigma}{2} \sum\limits_{k=2}^K a^{k-1}_R\l(\jump{\Eps^{k} }\r)^2 }_{BTs }
  \\ &+
  \underbrace{\frac{1}{2} \sum_{k=1}^K \frac{\Delta x_k}{2} \l(  \mat{\epsilon}_1^k
  \vec{\epsilon}_1^k,  \partial_x \vec{a}^k \r)_N}_{Int_d}
  \leq \Ep(t,N)
  ||\epsilon_1||_N.
\end{aligned}
\end{equation}

It is $BTs\geq0$.
If $Int_d\geq 0$, then \eqref{eq:Errorbound} has the same form as in
 in \cite{ nordstrom2007error}
and one may estimate/bound analogously
to \cite{kopriva2017error, nordstrom2007error} the error
in time.
The $\Ep$ term depends also on $N$,
but this has no influence in the estimation here.
We rewrite \eqref{eq:Errorbound} as
\begin{equation}\label{eq:Errorbound2}
\revsec{\frac{\d}{\d t}} ||\epsilon_1||_N +
\underbrace{ \frac{BTs+ Int_d}{  ||\epsilon_1||_N^2} }_{\eta(t) } ||\epsilon_1||_N  \leq \Ep(t).
\end{equation}
Like it was described in \cite{nordstrom2007error},
it is assumed that the mean value of $\eta(t)$ over any
finite time interval is bounded by a positive constant $\delta_0$
from below. This means that
$\overline{\eta}\geq \delta_0>0$.
Under the assumption for $u$, the right hand side $\Ep(t,N)$ is also bounded
in time and one can put $\max\limits_{s\in [0,\infty)} \Ep(s,N)\leq C_1<\infty$.
Applying these facts in \eqref{eq:Errorbound2} and integrating over time,
the following inequality for the error is obtained
\begin{equation}\label{eq:Errorbound3}
 ||\epsilon_1(t)||_N \leq \frac{1-\exp(-\delta_0 t)}{\delta_0} C_1,
\end{equation}
see \cite[Lemma 2.3]{nordstrom2007error} for details.
\begin{remark}\label{remark_int}
The term $Int_d$ is a crucial factor.
 If $\partial_x \vec{a}^k >0$, one may estimate the left side of \eqref{eq:Errorbound2}
using the minimum of the discrete values of $\vec{a}$.
Then, $Int_d\geq \frac{1}{2} \min\{\partial_x \vec{a}^k\} ||\epsilon_1||_N^2>0$ and
the above assumption
on $\eta$ is inevitably fulfilled.\\
Simultaneously, the term $Int_d$ can also destroy the error boundedness
if the derivatives of $a$ are negative.
It depends then on the sum of $BTs$ and $Int_d$. The upwind fluxes
can therefore rescue the error boundedness \eqref{eq:Errorbound3} whereas applying
the central flux ($\sigma= 0)$ will contribute to an unlimited growth of the error.
We demonstrate this in some examples in section \ref{sec:Num} and \revsec{make a first analytical estimation in \ref{subsec:analytical}.}
\end{remark}

%% file: 5_Error_bound_semidiscrete_Gauss.tex
\section{Error Behaviour using Gauß-Legendre Nodes }\label{sec:Error_Gauss}

Here, Gauß-Legendre nodes are used, yielding diagonal norm SBP operators not including
the boundary nodes, contrary to Gauß-Lobatto-Legendre nodes discussed in the previous
section~\ref{sec:Error_Lobatto}. Thus, care has to be taken of several potential problems.
Firstly, the restriction to the boundary and multiplication do not commute. Secondly, the
numerical flux functions \eqref{eq:Edge_based_cen_flux}--\eqref{eq:split_based_up_flux}
are now different from each other and have to be considered separately. \\
However, even if there are more problems, there are also some reasons to consider Gauß-Legendre
nodes. Indeed, Gauß-Legendre nodes have a higher order of accuracy in the quadrature and as
investigated in \cite{oeffner2018error}, for the linear advection equation with constant coefficients
using Gauß-Legendre nodes, the error reaches always faster its asymptotic value. Moreover, this
asymptotic value is lower than the corresponding one using Gauß-Lobatto-Legendre nodes. Furthermore, the
influence of the numerical fluxes is not that essential.\\
\revsec{Using $u=\Ip^N(u^k)+\epsilon^k_p$ in \eqref{eq:Continuous4},
where the terms are rearranged similar to section \ref{sec:Error_Lobatto},
we get analogously an equation similar to
\eqref{eq:Delta} except an additional error term due to the fact that boundary terms are not included.
We obtain
\begin{equation}\label{eq:Error_Gauss_2}
 \begin{aligned}
   &\quad \frac{\Delta x_k}{2} \l( \partial_t \vec{\Ip^N (u^k)}, \vec{\phi}^k \r)_N +\vec{\phi}^{k,T} \mat{R}[^T] \mat{B} \l( \vecfnumk \l(\Ip^N(u^k)^-, \Ip^N(u^k)^+\r)-
 \frac{1}{2} \l(\mat{R}\vec{a^k} \r)  \cdot
 \l(\mat{R} \vec{u} \r) \r)
 \\&\quad
 + \underbrace{\l( \frac{1}{2}  a^k \Ip^N(u^k)\revPP{\phi^k}\Bigg|_{-1}^1 -
 \vec{\phi}^{k,T} \mat{R}[^T] \mat{B} \l( \vecfnumk \l(\Ip^N(u^k)^-, \Ip^N(u^k)^+\r)-
 \frac{1}{2} \l(\mat{R}\vec{a^k} \r) \cdot
 \l(\mat{R} \vec{u} \r) \r) \r)}_{=:\epsilon_2^k(a^k)}
 \\&\quad
 - \frac{1}{2} \l(  \mat{a^k} \vec{\Ip^N(u^k)} ,  \partial_\xi \vec{\phi}^k \r)_N
  +\frac{1}{2} \l( \partial_\xi \vec{\Ip^N(u^k)}, \mat{a^k} \vec{\phi}^k    \r)_N
  +\frac{1}{2} \l( \partial_\xi \vec{a}^k , \mat{\Ip^N(u^k)}  \vec{\phi}^k \r)_N\\
 &= \frac{\Delta x_k}{2} \est{ \hat{T}^k(u^k) , \phi^k} +\frac{\Delta x_k}{4} \est{Q_1(u^k), \partial_x \phi^k }
 \\&\quad
 +\frac{\Delta x_k}{4}  \Big\{
  \l( \vec{Q(u^k)}, \vec{\phi}^k \r)_N -
   \l( \vec{Q_1(u^k)}, \partial_x \vec{\phi}^k \r)_N +
   \revPP{ \l(  \vec{Q_2(u^k)}, \mat{a}^k\vec{\phi}^k \r)_N } \\ & \quad +
     \l( \partial_x \vec{a}^k, \mat{Q_3(u^k)}\vec{\phi}^k \r)_N \Big \}
 \end{aligned}
\end{equation}
 with
 \begin{equation*}
 \begin{aligned}
     \hat{T}^k(u^k) :=& -\Bigg\{\partial_t \epsilon_p^k+\frac{1}{2} \l( \partial_x \l( a^k \epsilon_p^k\r)
     +\epsilon_p^k \partial_x a^k \partial_x \epsilon_p^k \r)\\
 &+\frac{1}{2} \l( Q(u^k) +a^kQ_2(u^k) +Q_3(u^k) \partial_x a^k \r)
 \Bigg\}.
 \end{aligned}
 \end{equation*}}
\revsec{Following the approach from section \ref{sec:Error_Lobatto}
we get the estimate\footnote{Details of main steps can also be found in the appendix \ref{sec:Appendix}.}}
\begin{equation}\label{eq:R_absch_2}
\begin{aligned}
&\frac{1}{2}   \frac{\d}{\d t} \sum_{k=1}^K \frac{\Delta x_k}{2}||\epsilon_1^k||_N^2+\sum_{k=1}^K
 \vec{\epsilon}_1^{k,T} \mat{R}[^T] \mat{B}
\l( \vecfnumk
 \l( (\epsilon_1^k )^{-},(\epsilon_1^k )^{+}\r)-\frac{1}{2 }\l( \mat{R}\vec{a^k} \r) \cdot
 \l(\mat{R} \vec{\epsilon}_1^k \r) \r)
   \\ &
   +\underbrace{\frac{1}{2} \sum_{k=1}^K \frac{\Delta x_k}{2}
   \revPP{\l(  \partial_x \vec{a}^k, \mat{\epsilon}_1^k
  \vec{\epsilon}_1^k \r)_N} }_{Int_d}
 \leq - \underbrace{\frac{1}{2} \sum\limits_{k=1}^K \frac{\Delta x_k}{2} \epsilon_2^k(a^k)}_{:=\Theta_2} +\\
 &\underbrace{ \l\{ c_1||T ||+\frac{cN}{2}||Q_1|| +
 \frac{1}{2}\l( ||Q||_N+ N||Q_1||_N +||a Q_2||_N +||Q_3\partial_x a||_N \r)   \r\}  }_{:=  \hat{\Ep}_G (t,N) }
 ||\epsilon_1||_N .
 \end{aligned}
\end{equation}
\begin{remark}
The sum of the terms $\epsilon_2^k$ depends on $a$ and the interpolation of the flux functions.
It is given by the
formula
\begin{equation*}
\begin{aligned}
  \epsilon_2^k(a^k):=&\Bigg( \frac{1}{2}  a^k\epsilon_1^k \Ip^N(u^k)\Bigg|_{-1}^1 -
 \vec{\epsilon}_{\revPP{1}}^{k,T} \mat{R}[^T] \mat{B} \Big( \vecfnumk \l(\Ip^N(u^k)^-, \Ip^N(u^k)^+\r)\\
 &-
 \frac{1}{2} \l(\mat{R}\vec{a^k} \r)
 \l(\mat{R} \vec{u} \r) \Big) \Bigg).
\end{aligned}
 \end{equation*}
Using Gauß-Lobatto nodes and an upwind flux, these terms are zero, see section \ref{sec:Error_Lobatto}.
If the sum over all elements is positive, i.e. $\Theta_2\geq 0$, then this term decreases the upper
bound of the error $\epsilon_1$.

If $\Theta_2 <0$, then it increases the total error. The error depends on $u$, $a$ and the jumps between
interfaces. Under the assumption that $u$ is continuous, $\Theta_2$ will be bounded from below,
resulting in an upper bound on the right side. Nevertheless, this makes it hard to study the
behaviour of the total error analytically.
\end{remark}

We consider the first line of \eqref{eq:R_absch_2}, especially the term
\begin{align*}
 \sum_{k=1}^K
 \vec{\epsilon}_1^{k,T} \mat{R}[^T] \mat{B}
\l( \vecfnumk
 \l( (\epsilon_1^k )^{-},(\epsilon_1^k )^{+}\r)-\frac{1}{2 }\l( \mat{R}\vec{a^k} \r) \cdot
 \l(\mat{R} \vec{\epsilon}_1^k \r) \r)
\end{align*}
with  different flux functions \eqref{eq:Edge_based_cen_flux}--\eqref{eq:split_based_up_flux}.
In \cite{ranocha2018generalised}, different assumptions on $a$ have already been formulated for stability and conservation of the numerical schemes.
First, we consider the general case.
One may recognise the problems which arise by considering variable coefficients in the model problem \eqref{eq:Model_problem}.
Following this, we will formulate analogues assumptions to \cite[Theorem 3.4]{ranocha2018generalised} and proceed with our analysis.

We split the sum in three terms (one for the left physical boundary,
one for the right physical boundary and a sum over the internal element endpoints), and we get
\begin{align*}
 &\quad \sum\limits_{k=1}^K \vec{\epsilon}_1^{k,T} \mat{R}[^T] \mat{B}
\l( \vecfnumk
 \l( (\epsilon_1^k )^{-},(\epsilon_1^k )^{+}\r)-\frac{1}{2 }\l( \mat{R}\vec{a^k} \r) \cdot
 \l(\mat{R} \vec{\epsilon}_1^k \r) \r)   \\
 &=- \Eps_L^1 \l( f^{\mathrm{num},1}_L -\frac{1}{2} a_L^1\Eps_L^1 \r)+
 \sum\limits_{k=2}^K \Bigg( f^{\mathrm{num},k}_L \l( \Eps_R^{k-1}-\Eps_L^{k}  \r)
 \\
  &\quad -\frac{1}{2} \left(a_R^{k-1}\l(\Eps_R^{k-1}\r)^2- a_L^k \l(\Eps_L^{k}\r)^2  \right) \Bigg)
 + \Eps_R^K \l(f^{\mathrm{num},K}_R -\frac{1}{2} a_R^K\Eps^K_R \r).
\end{align*}
We describe with $\Eps_i$ ($i=L,R$) the approximation error $\epsilon_1$, the indices give
 the position in the elements,  $f^{\mathrm{num},k}_L:=f^{\mathrm{num},k} \l(\Eps^{k-1}_R,\Eps^{k}_L \r)$, $ f^{\mathrm{num},1}_L:=f^{\mathrm{num},1} \l(0,\Eps^{1}_L \r)$
 and $f^{\mathrm{num},K}_R:= f^{\mathrm{num},1} \l(\Eps^K_R,0 \r) $. The external states for the physical boundary contributions are zero,
 because $\Ip^N(u)^1=g$ at the left boundary and the external state for $U^1$ is set to $g$.
 The selection of the numerical flux functions \eqref{eq:Edge_based_cen_flux}--\eqref{eq:split_based_up_flux} has an
 influence on the behaviours of the errors and we have to be careful in our study.
 If the interpolation of $a$ is exact and $a$ is continuous over the inter-element boundaries, then the
 influence of the numerical fluxes can be simplified essentially and we are able to analyse the long time error behaviours.
 We will formulate this in detail for the first flux under consideration, the edge based central flux \eqref{eq:Edge_based_cen_flux}.
\begin{itemize}
 \item Edge based central flux $\fnum(u_-, u_+)= a(x)\frac{u_-+u_+}{2}$:
We get for the terms in the sum
\begin{align*}
 &\frac{1}{2}a^k(x_L)\l(\Eps^{k-1}_R+\Eps^k_L \r) \l(\Eps^{k-1}_R-\Eps_L^k \r)-\frac{1}{2} \l( a_R^{k-1}\l( \Eps_R^{k-1 } \r)^2    -a_L^k \l(\Eps^k_L \r)^2 \r)\\
 =&\frac{1}{2}a^k(x_L) \l( \l( \Eps_R^{k-1 } \r)^2    -\l(\Eps^k_L \r)^2 \r)-\frac{1}{2} \l( a_R^{k-1}\l( \Eps_R^{k-1 } \r)^2    -a_L^k \l(\Eps^k_L \r)^2 \r)\\
 =& \frac{1}{2} \l( \Eps^{k-1}_R\r)^2 \l( a^k(x_L)- a_R^{k-1} \r)+\frac{1}{2} \l(\Eps^k_L \r)^2 \l(a_L^k-a(x_L) \r) =0.
\end{align*}
If the interpolation of $a$ is exact and $a$ is continuous, the brackets of $a$ will be zero,
because $a^k(x_L)=a^k_L=a^{k-1}(x_R)=a^{k-1}_R$.
If this is not the case, we get additional terms that can be positive or negative depending
on  brackets.
On the boundaries, one obtains
 \begin{align*}
  &\text{ left:} &-\Eps_L^1 \l( f^{\mathrm{num},1}_L-\frac{1}{2}a_L^1\Eps_L^1 \r)=
  -\Eps_L^1 \l( \frac{a^1(x_L)}{2}\Eps_L^1-\frac{a_L^1}{2} \Eps_L^1 \r)\\
 &  & = \frac{1}{2} \l(\Eps_L^1 \r)^2 \l(a_L^1-a^1(x_L) \r) = 0, \\
&\text{ right:} & \; \Eps^K_R  \l( f^{\mathrm{num},K}_R -\frac{1}{2}a_R^K \Eps^K_R \r)
  = \frac{1}{2} \l(\Eps_R^K \r)^2 \l(a^K(x_R)-a_R^K \r) = 0.
 \end{align*}
\item  \revsec{Using this approach, we get the following results where
the details of the calculation can be found in the
appendix \ref{sec:Appendix}:
\begin{table}[!ht]
\centering
 \begin{tabular}{|l|c|c|c|}
 \hline
 Fluxes & Interiour& Left & Right
 \\ \hline \hline
 Split central & $0$& $0$& $0$
 \\ \hline
 Edge bases upwind & $\frac{1}{2} a_R^{k-1} \jump{E_R^{k-1}}^2$ & $\frac{1}{2} \l(\Eps^1_L \r)a_L^1$ & $\l(\Eps^k_R \r)^2 \l(\frac{a_R^{K}}{2} \r)$
 \\ \hline
 Split upwind & $\frac{1}{2} a_R^{k-1}
 \jump{E_R^{k-1}}^2$ & $\frac{1}{2} \l(\Eps^1_L \r)a_L^1$&
 $\l(\Eps^k_R \r)^2 \l(\frac{a_R^{K}}{2} \r)$ \\ \hline
\end{tabular}
\caption{Error terms of the numerical fluxes}
\label{Table:error_flux}
\end{table}}

\revsec{For the calculation of the split upwind flux, we apply the assumptions of the exactness of the interpolation
and the continuity of $a$.
}
\item Unsplit upwind flux $\fnum(u_-,u_-)=(au)_-$. \\
Unfortunately, for the unsplit numerical fluxes  \eqref{eq:unsplit_based_cen_flux}, \eqref{eq:unsplit_based_up_flux} we are not able
to find such a simplification as above, since the restriction of the product can not be compared to the product of the restriction.
This issue triggers also stability problems, see \cite{ranocha2018generalised} for details. We formulate this now for the unsplit upwind flux as an example. It is:
\begin{align*}
 &(a\Eps)^{k-1}_R \l( \Eps^{k-1}_R-\Eps_L^k \r)-\frac{1}{2} \l( a_R^{k-1}\l( \Eps_R^{k-1 } \r)^2    -a_L^k \l(\Eps^k_L \r)^2 \r)\\
  =& \frac{1}{2}\l(  \left(2(a\Eps)^{k-1}_R \Eps_R^{k-1} -a^{k-1}_R\l(\Eps_R^{k-1}\right)^2\right)-2a^{k-1}_R\Eps_L^k \Eps^{k-1}_R +a_L^k \l(\Eps^k_L\r)^2 \r).
\end{align*}
Because of $(a\Eps)^{k-1}_R\neq a^{k-1}_R\Eps^{k-1}_R$ in general, a further simplification is in this case not possible anymore.
The following error bounds are \textbf{only} valid for the split numerical fluxes.
Nevertheless, we test also the unsplit fluxes in the next section.
\end{itemize}
By comparison, one may recognise that the split upwind flux is equal the edge upwind flux and analogously for the central
fluxes under assumptions. Using central fluxes leads to no additional terms in the inequality \eqref{eq:R_absch_2},
whereas using upwind fluxes does. If the restrictions $a_{L/R}$ to the boundary are positive\footnote{This assumption
is already formulated in \cite[Theorem~3.4]{ranocha2018generalised} to guarantee stability and conservation of the numerical schemes.},
all of these terms are positive.
We reformulate the energy inequality \eqref{eq:R_absch_2} as
\begin{equation}\label{eq:Errorbound_Gauss}
\begin{aligned}
   &\frac{1}{2} \frac{\d}{\d t} ||\epsilon_1||^2_N +
  \underbrace{\frac{\sigma}{2} \l( a^K_R\l(\Eps^K_R  \r)^2 +a_L^1\l(\Eps_L^1\r)^2  \r) +
  \frac{\sigma}{2} \sum\limits_{k=2}^K a^{k-1}_R\l(\jump{\Eps^{k} }\r)^2 }_{BTs }+
 Int_d \\
  \leq&  \hat{\Ep}_G (t,N)
  ||\epsilon_1||_N -\Theta_2,
\end{aligned}
\end{equation}
where $\sigma$ is zero (central flux) or one (upwind flux).
The energy growth energy inequality \eqref{eq:Errorbound_Gauss} is similar to  \eqref{eq:R_absch_2}.
The only difference is the term $\Theta_2$, which will yield a smaller upper bound under
the condition $\Theta_2\geq 0$. We follow the steps of section \eqref{sec:Error_Lobatto}
and get
\begin{equation}\label{eq:eta_def}
  \frac{\d}{\d t} ||\epsilon_1||_N +\underbrace{\frac{BTs +Int_d+\Theta_2}{||\epsilon_1||^2_N}}_{\eta_G(t)} ||\epsilon_1||_N\leq  \hat{\Ep}_G (t,N).
\end{equation}
We have to assume that we can bound the mean value of $\eta_G(t) $ by a positive constant $\delta_G$ from below.
If already  $Int_d +\Theta_2>0$,
this is actually met without restrictions.
\revPP{So, using the central fluxes ($\sigma=0$) does not yield
to problems. Simultaneously, if $BTs+Int_d +\Theta_2$ overall is positive,
the requirement on every $a^k_{L,R}$ to be non-negative can be weaken to make the estimations,
but one should have in mind that  the positivity of $a^k_{L,R}$ is a condition to prove stability.}
This means\footnote{$\delta_0$ from section \ref{sec:Error_Lobatto}, inequality \eqref{eq:Errorbound3}.}
$\overline{\eta_G}(t) \geq \delta_G >\delta_0>0$.
Under the assumption of $u$, the right hand side $\hat{\Ep}_G (t,N) $ is also bounded
in time and one can put $\max\limits_{s\in [0,\infty)}\hat{\Ep}_G (t,N) \leq C_2<\infty$.
Applying this in \eqref{eq:Errorbound2} and integrating over time, the inequality
for the error follows as
\begin{equation}\label{eq:Errorbound4}
 ||\epsilon_1(t)||_N \leq \frac{1-\exp(-\delta_G t)}{\delta_G} C_2.
\end{equation}
Since $\delta_G>\delta_0$, the error using Gauß-Legendre nodes will reach its
asymptotic value faster than the error using a Gauß-Lobatto-Legendre basis. We see this
behaviour in our numerical simulations in the next section.

%% file: 6_Numerics.tex
\section{Numerical Examples}\label{sec:Num}

In this section, we present some numerical experiments using the constructed schemes. We focus on the influence
of the different numerical fluxes on the long time behaviour of the error.
From \cite{kopriva2017error, oeffner2018error}, we know that in case of constant coefficients the choice of the
numerical flux has an essential influence on the error behaviour, especially in the Gauß-Lobatto-Legendre case.

We consider our model problem, the linear advection equation
\begin{equation*}\tag{\ref{eq:Model_problem}}
\begin{aligned}
   \partial_t u(t,x)+\partial_x (a(x)u(t,x))&=0, &t>0,\; x\in (x_L,\;x_R),\\
   u(t,x_L)&=g_L(t), & t\geq 0, \\
   u(0,x)&=u_0(x), & x\in (x_L,\;x_R),
\end{aligned}
\end{equation*}
with smooth speed $a(x)>0$, initial condition $u_0$ and boundary condition $g_L$.
The solution $u$ of the corresponding Cauchy problem can be calculated by the method of characteristics,
see e.g. \cite[Chapter 3]{bressan2000hyperbolic}.
As time integrator, we use the fourth order, ten stage, strong stability
preserving Runge-Kutta method of \cite{ketcheson2008highly} and the time step is
chosen such that the time integration error is negligible. Although the term
``strong-stability preserving'' means the preservation of stability properties of the
explicit Euler method and the explicit Euler method is not stable for our numerical
experiments, this fourth order Runge-Kutta method is strongly stable for linear
equations \cite{ranocha2018L2stability}. All elements are of uniform size.

\subsection{Coefficient \texorpdfstring{$a(x)=x$}{a(x)=x}}\label{subsec:a(x)=x}

In our first experiment, we choose $a(x)=x$ with initial condition $u_0(x) =\sin (12(x-0.1))$.
The interval is $[x_L,x_R] = [0,2\pi]$ and we choose the inflow boundary condition such that
we get the solution
\begin{equation*}
 u(t,x)= \exp(-t) u_0\bigl( x \exp(-t) \bigr).
\end{equation*}
For the coefficient $a(x)=x$, the first derivative of $a$ is strictly positive,
implying $Int_d > 0$.

\begin{figure}[!htp]
\centering
  \begin{subfigure}[b]{0.495\textwidth}
  \centering
    \includegraphics[width=\textwidth]{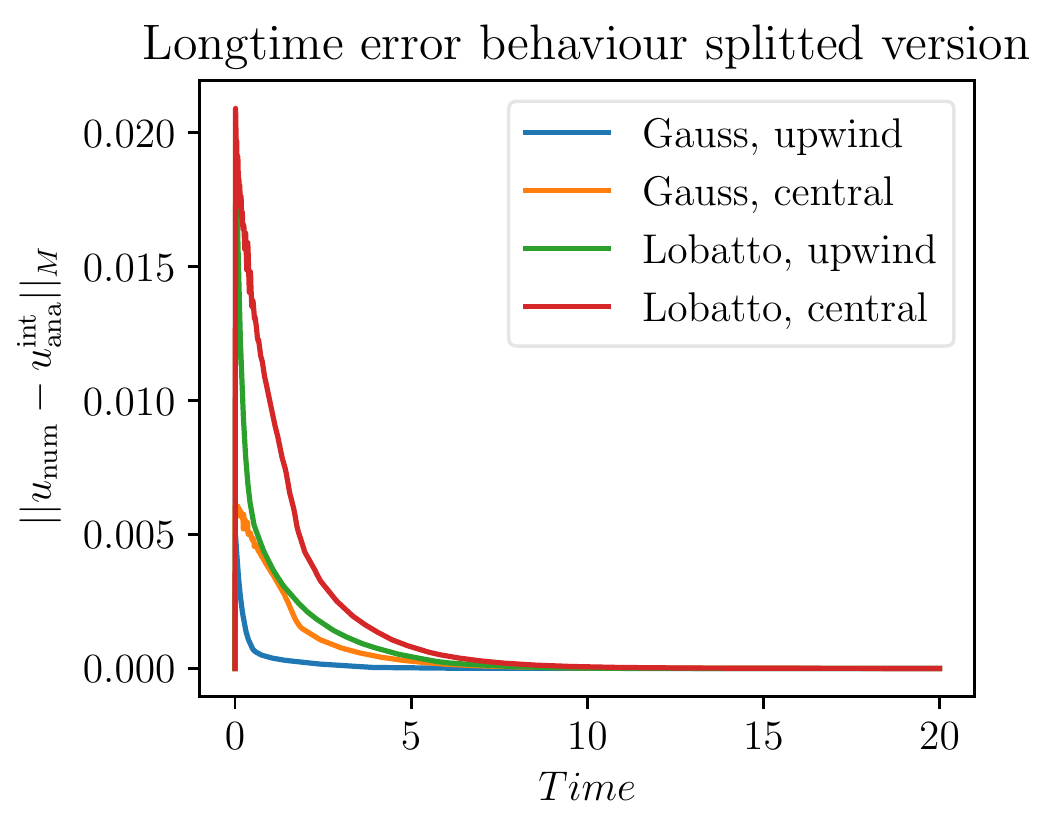}
    \caption{$N=3$, $K=40$, $t=20$.}
  \end{subfigure}%
  ~
  \begin{subfigure}[b]{0.495\textwidth}
  \centering
    \includegraphics[width=\textwidth]{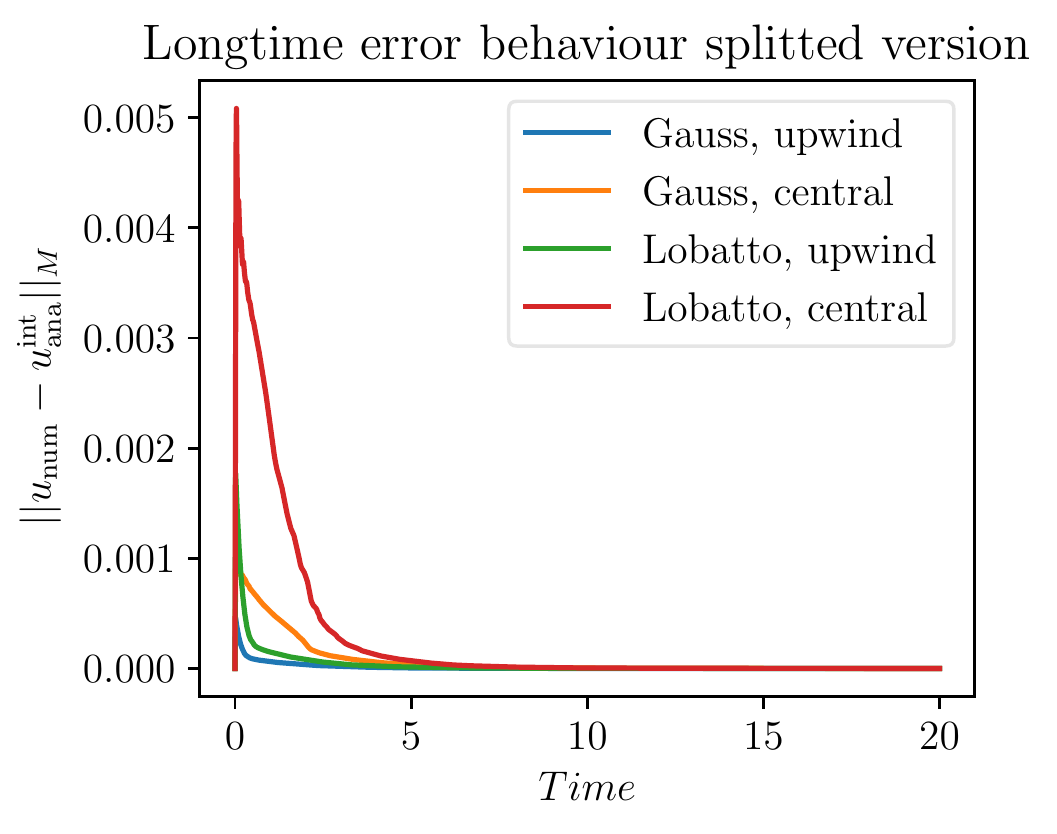}
    \caption{$N=4$, $K=40$, $t=20$.}
  \end{subfigure}\\%

 \begin{subfigure}[b]{0.495\textwidth}
  \centering
    \includegraphics[width=\textwidth]{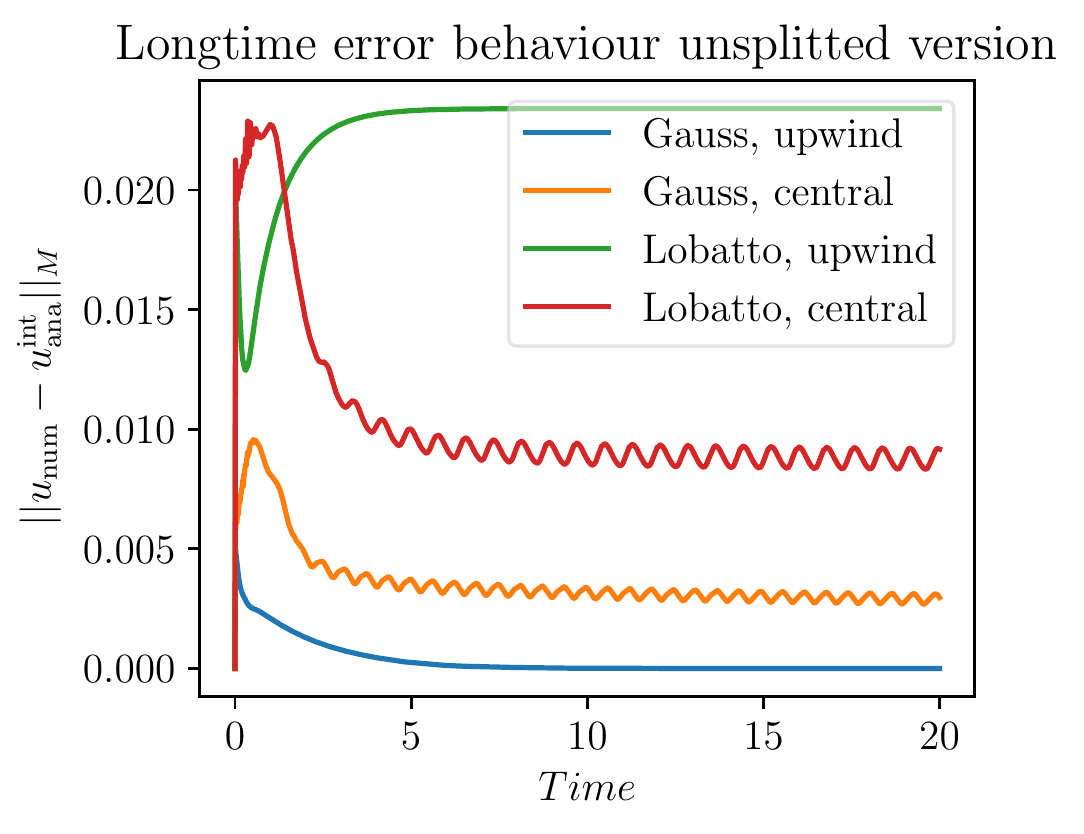}
    \caption{$N=3$, $K=40$, $t=20$.}
  \end{subfigure}%
  ~
  \begin{subfigure}[b]{0.495\textwidth}
  \centering
    \includegraphics[width=\textwidth]{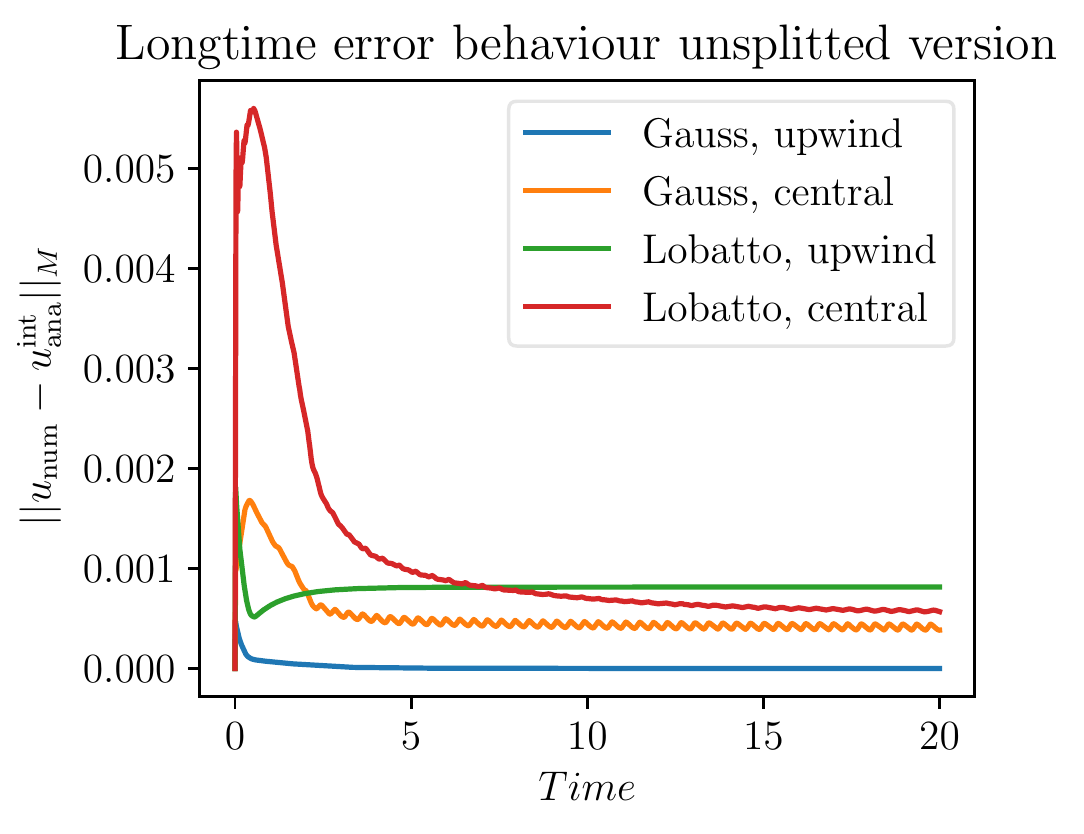}
    \caption{$N=4$, $K=40$, $t=20$.}
  \end{subfigure}%
  \caption{Errors of numerical solutions using split and unsplit discretisations, $a(x) = x$.}
  \label{fig:Splitted_unsplitt_gauss_Lobatto}
\end{figure}

In our first simulation, we use $K=40$ elements and calculate the solutions up to $t=20$ with
\num{200000} time steps. In \autoref{fig:Splitted_unsplitt_gauss_Lobatto}, we plot the long
time error behaviour using polynomial degrees three and four. One recognizes that in all
cases the error remains bounded in time.

In the first row of \autoref{fig:Splitted_unsplitt_gauss_Lobatto}, all terms (surface, flux and volume) are split whereas in the second row
they are not. We see that the error for the split version behaves like in the case of
constant coefficients \cite{kopriva2017error, oeffner2018error}. We mean that the errors using
the upwind fluxes are always lower than the ones using central fluxes and one may recognize
that we have some noisy behaviour using the central fluxes. Using upwind fluxes, the error
reaches its asymptotic value faster than for the central fluxes.

In the second row, the unsplit discretisation is used. We recognize that we lose the predictions
from \cite{kopriva2017error, oeffner2018error} that applying the upwind flux yields a more
accurate solution. The absolute value is also bigger applying the unsplit versions and we have
again the noisy behaviour by applying the central fluxes.

\begin{figure}[!htp]
\centering
  \begin{subfigure}[b]{0.495\textwidth}
    \includegraphics[width=\textwidth]{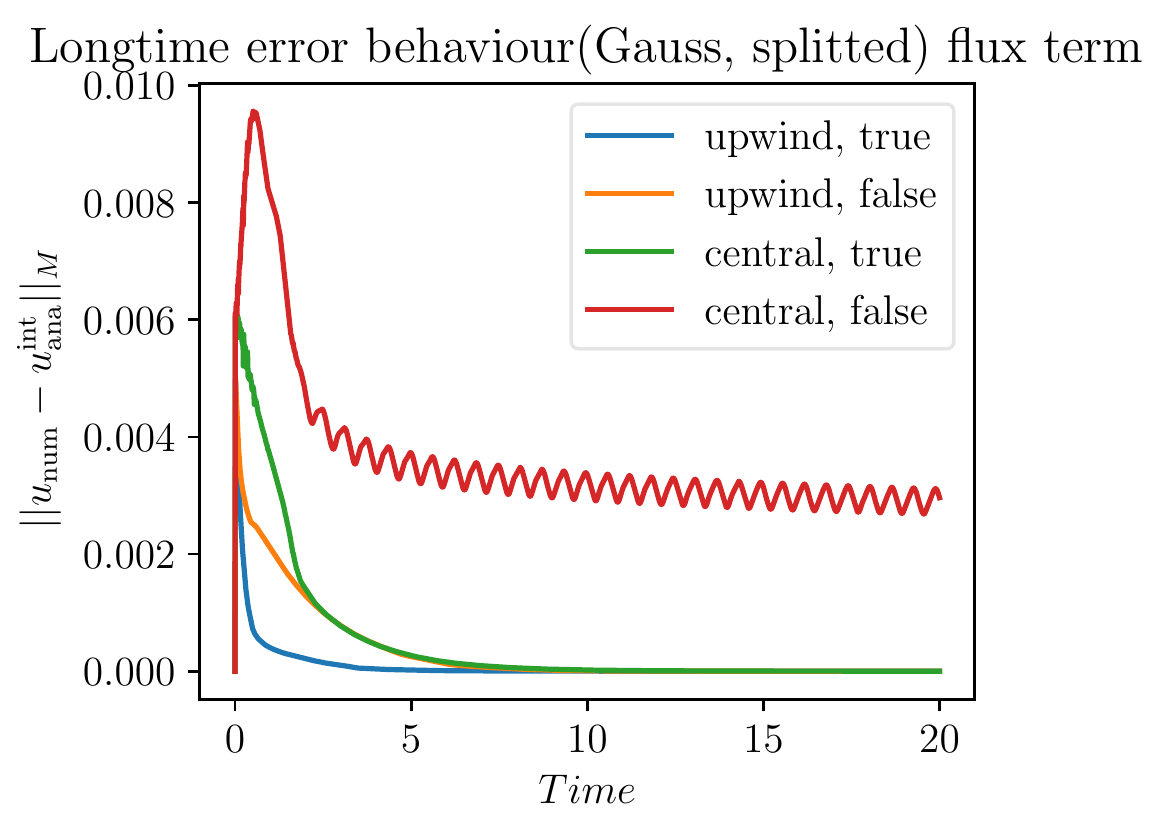}
    \caption{$N=3$, $K=40$, $t=20$.}
  \end{subfigure}%
  ~
  \begin{subfigure}[b]{0.495\textwidth}
    \includegraphics[width=\textwidth]{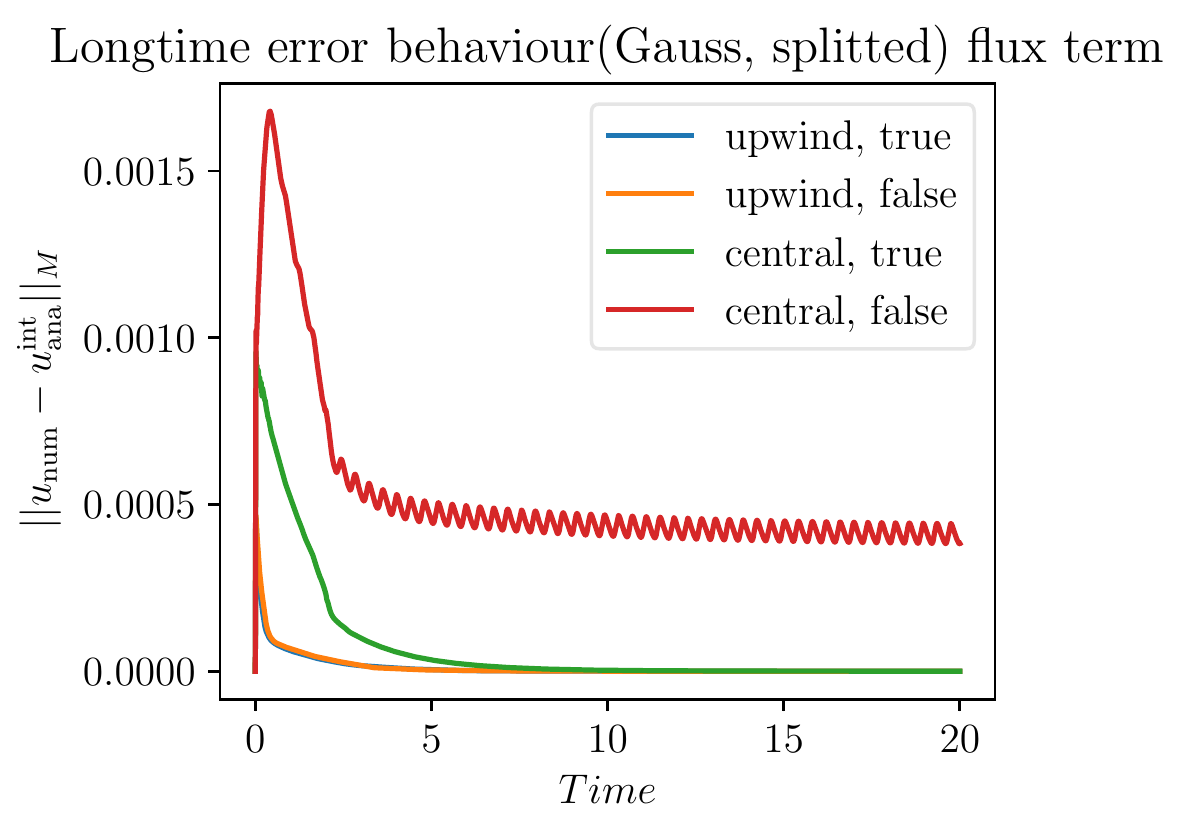}
    \caption{$N=4$, $K=40$, $t=20$.}
  \end{subfigure}%
  \caption{Errors of numerical solutions using the split form and both split (true) and
           unsplit (false) numerical fluxes.}
  \label{fig:flux_split_unsplit_gauss}
\end{figure}

Comparing all four plots, we recognize that the best results are obtained by using Gauß-Legendre nodes and
the split discretisation. Therefore, we have a closer look on this.
In \autoref{fig:flux_split_unsplit_gauss}, we consider only Gauß nodes and compare the split
numerical fluxes and the unsplit numerical fluxes (with split surface and volume terms).
\emph{True} in the legend indicates the split numerical fluxes and
\emph{false} the unsplit ones.
The experiment on the left-hand side demonstrates clearly that the noisy behavior for the central flux
transfers also to the application of Gauß-Legendre nodes if all terms are split. Furthermore, we can
hardly indicate some difference between the usage of split and unsplit upwind fluxes here, whereas
we have a slight different behaviour in the usage of the central fluxes.
The test indicates that the split discretisation (volume/surface and numerical fluxes)
should be preferred, matching our stability analysis.

\subsection{Coefficient \texorpdfstring{$a(x)=x^2$}{a(x)=x²}}\label{subsec:a(x)=x2}

In our second experiment, we choose $a(x)=x^2$ with initial condition
$u_0(x) = \cos\left( \frac{\pi x}{2} \right)$. The interval is $[x_L,x_R] = [0.1,1]$
and we choose the inflow boundary condition according to the solution
\begin{equation}\label{sol_a_x_x}
  u(t,x)= \frac{u_0\bigl(x/(1+tx)\bigr)}{(1+tx)^2}.
\end{equation}

\begin{figure}[!htp]
\centering
  \begin{subfigure}[b]{0.495\textwidth}
    \includegraphics[width=\textwidth]{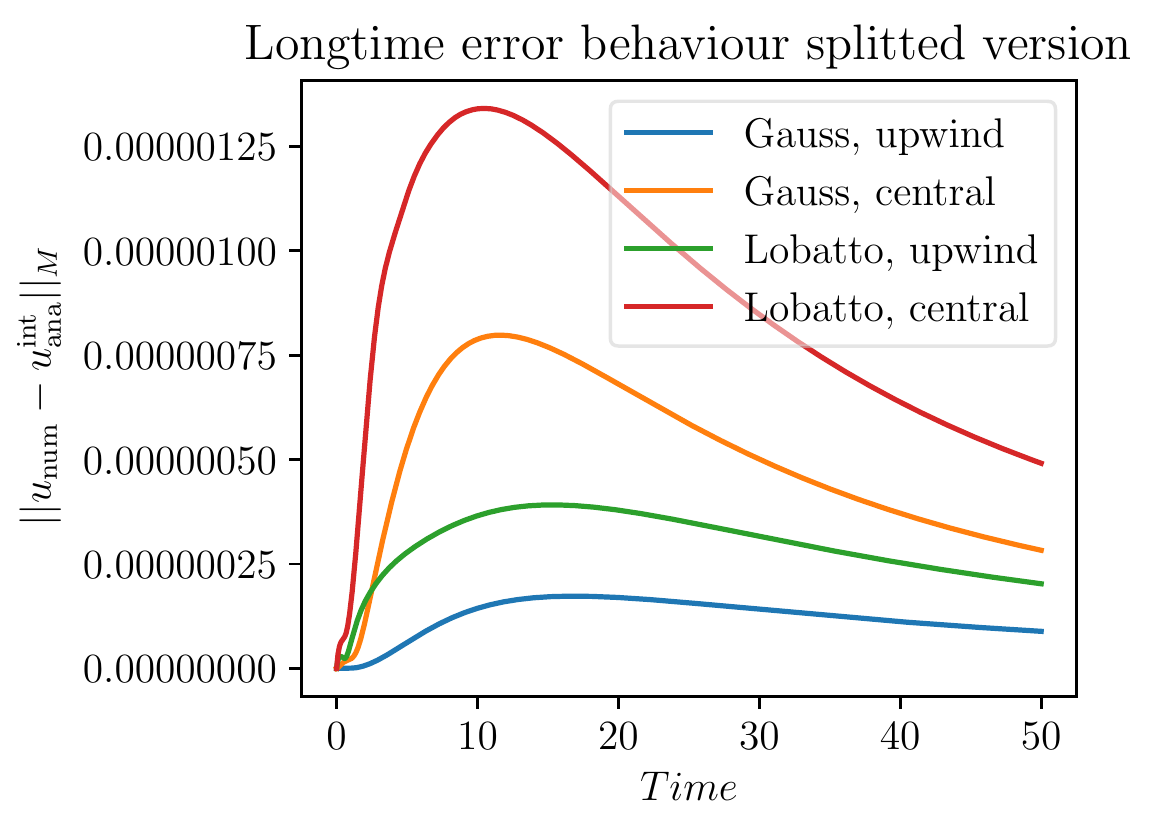}
    \caption{$N=3$, $K=40$, final time $t = 50$.}
  \end{subfigure}%
  ~
  \begin{subfigure}[b]{0.495\textwidth}
    \includegraphics[width=\textwidth]{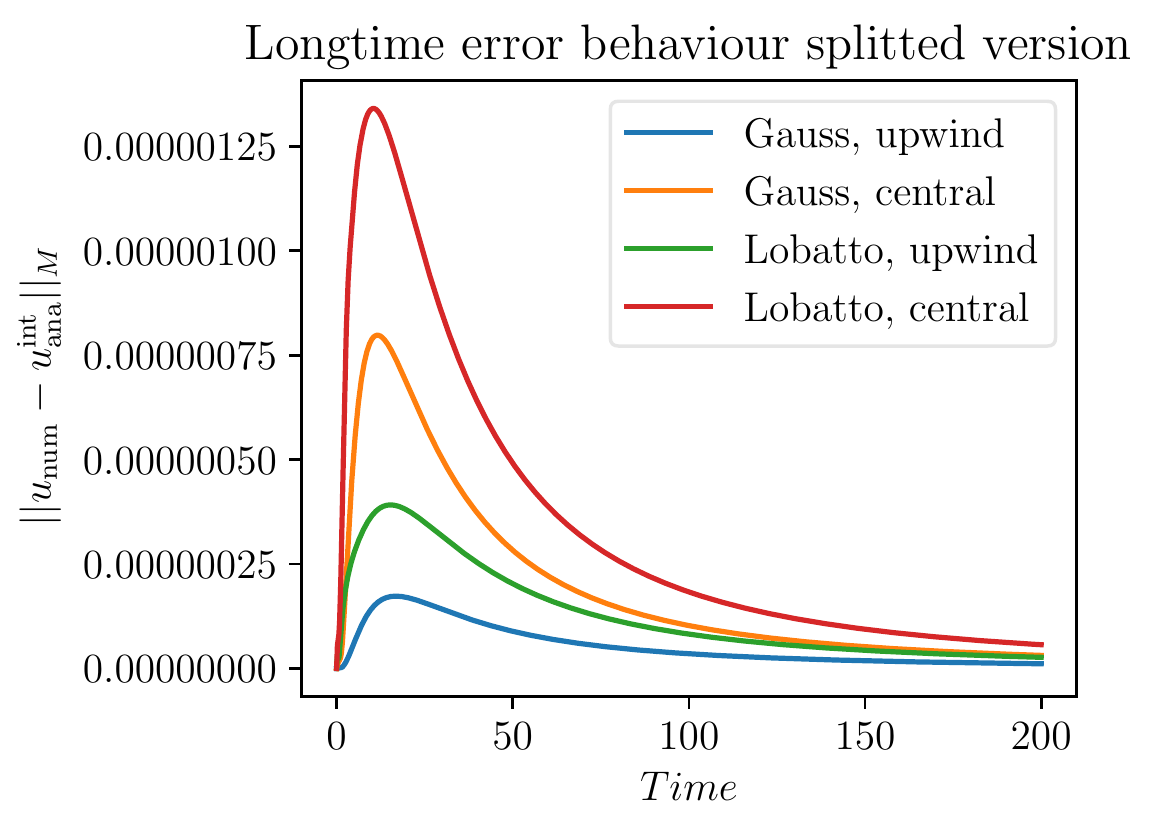}
    \caption{$N=3$, $K=40$, final time $t = 200$.}
  \end{subfigure}%
  \caption{Errors of discretisations for $a(x) = x^2$, $[x_L,x_R] = [0.1,1]$,
           and $u_0(x) = \cos\left( \frac{\pi x}{2} \right)$.}
  \label{fig:flux_splitted_unsplit_gauss2}
\end{figure}

\begin{figure}
\centering
  \includegraphics[width=0.5\textwidth]{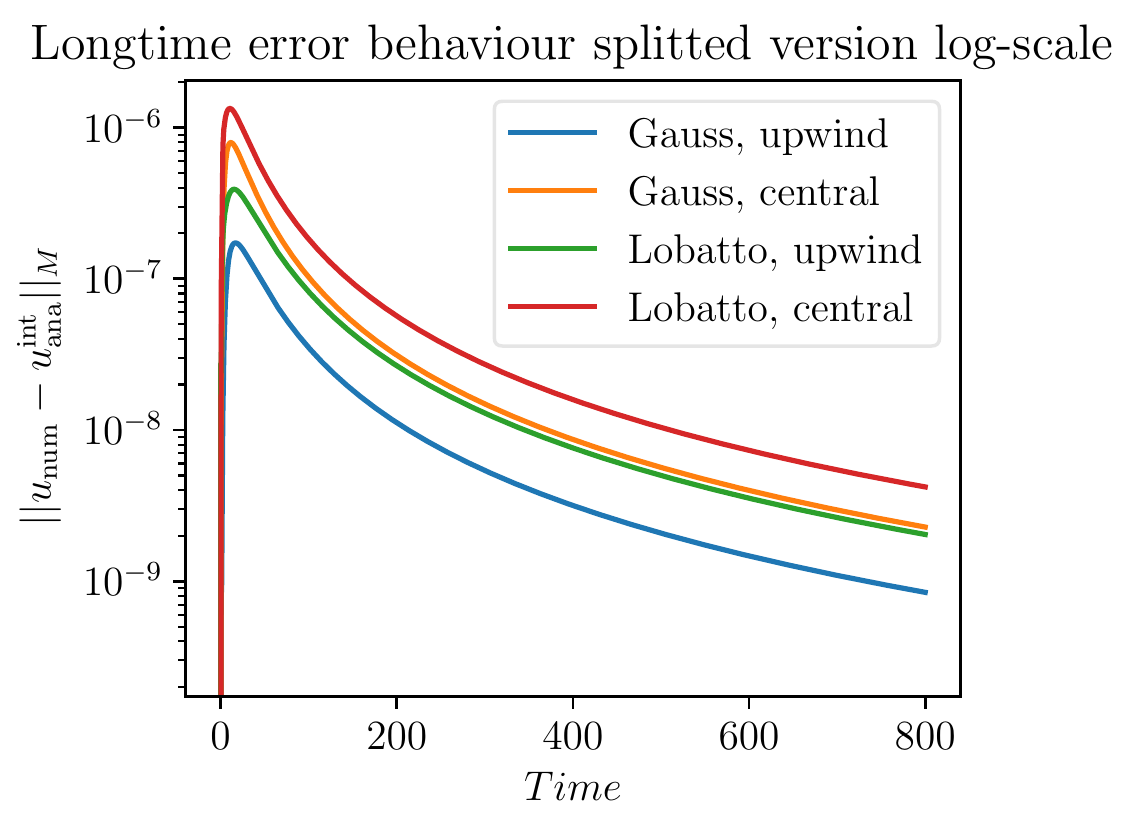}
  \caption{Error behaviour on a logarithmic scale for $a(x) =  x^2$, $t=800$.}
  \label{fig:flux_splitted_unsplit_gauss2_log}
\end{figure}

In our simulation shown in \autoref{fig:flux_splitted_unsplit_gauss2},
we apply different numbers of time steps up to $t = 200$.
First, we recognize that all errors are bounded in time,
but different from the first case we do not have any noisy behavior
of the central fluxes, at least we can not identify some.
Simultaneously, the unsplit central flux error  with Lobatto nodes
increases at first rapidly before it finally tends to its
asymptotic value. In all cases, the errors are small
but we get always the \emph{best} results by applying Gauß-Legendre nodes.
Nevertheless, it takes a lot of time for the errors to reach the asymptotic values.
Even at time $t=800$, the asymptotic is still not reached, cf. \autoref{fig:flux_splitted_unsplit_gauss2_log}.

\begin{figure}[!htp]
\centering
  \begin{subfigure}[b]{0.324\textwidth}
    \includegraphics[width=\textwidth]{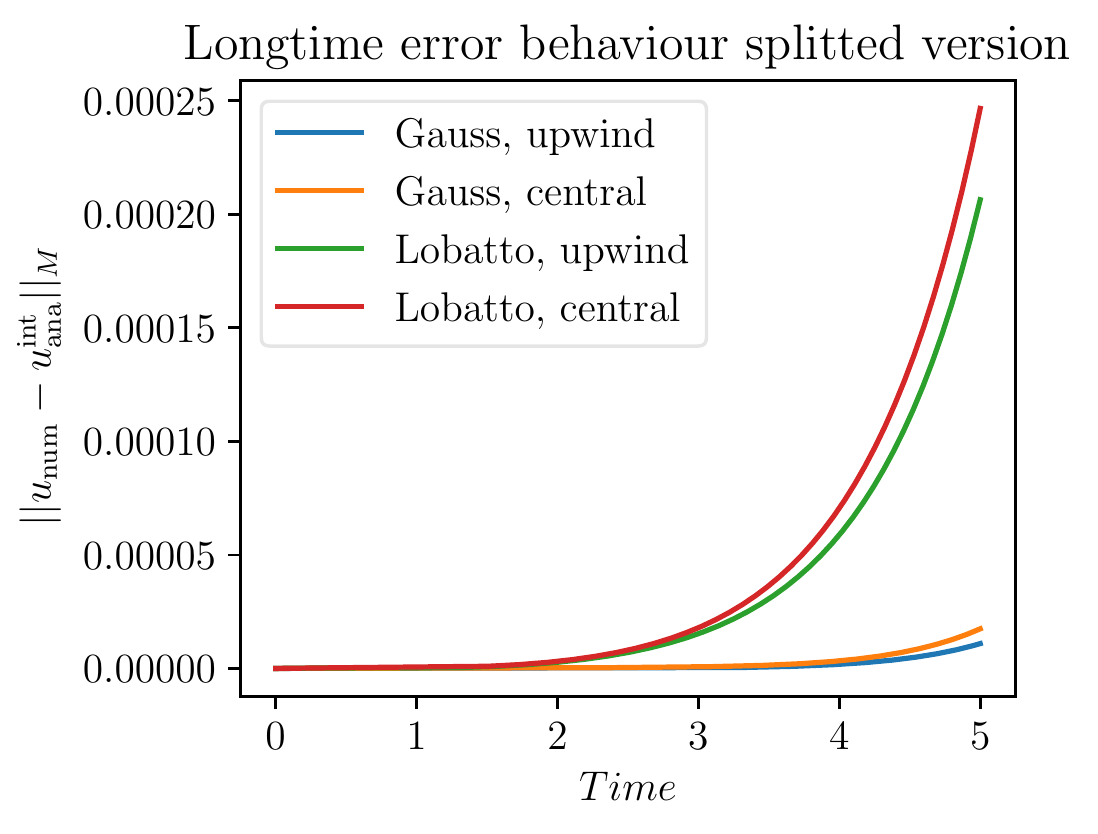}
    \caption{$N=3$, $K=40$, $t=5$.}
  \end{subfigure}%
  ~
  \begin{subfigure}[b]{0.324\textwidth}
    \includegraphics[width=\textwidth]{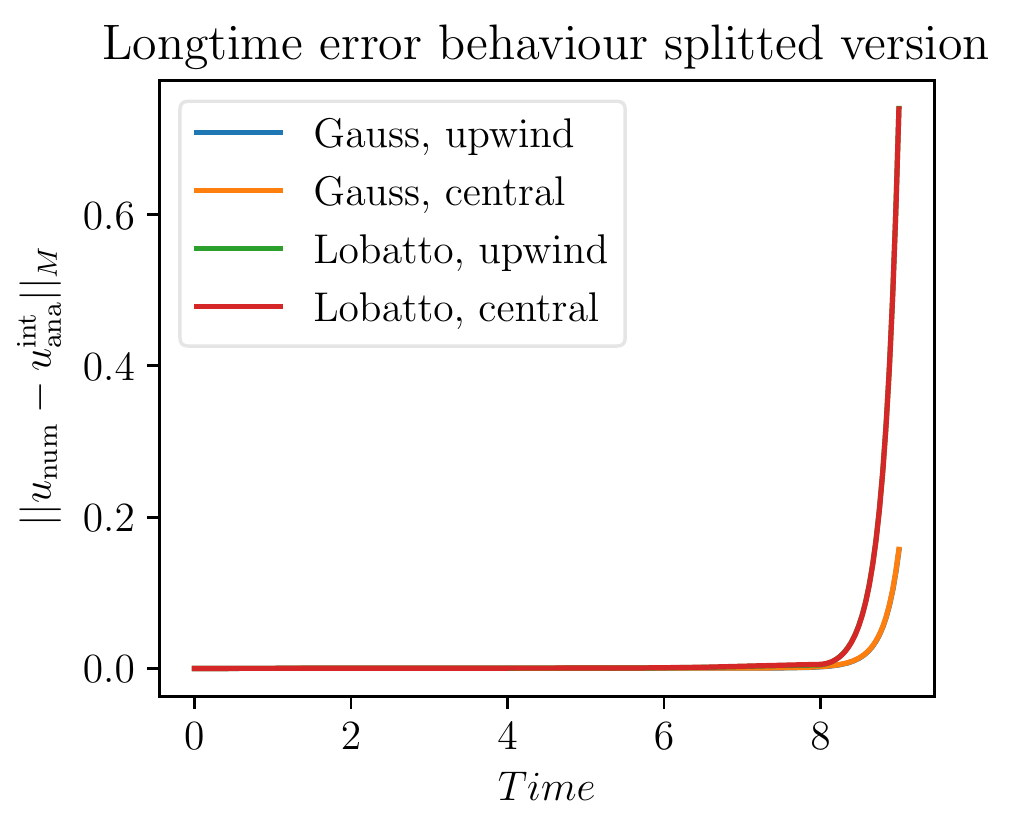}
    \caption{$N=3$, $K=40$, $t=9$.}
  \end{subfigure}%
    ~
  \begin{subfigure}[b]{0.324\textwidth}
    \includegraphics[width=\textwidth]{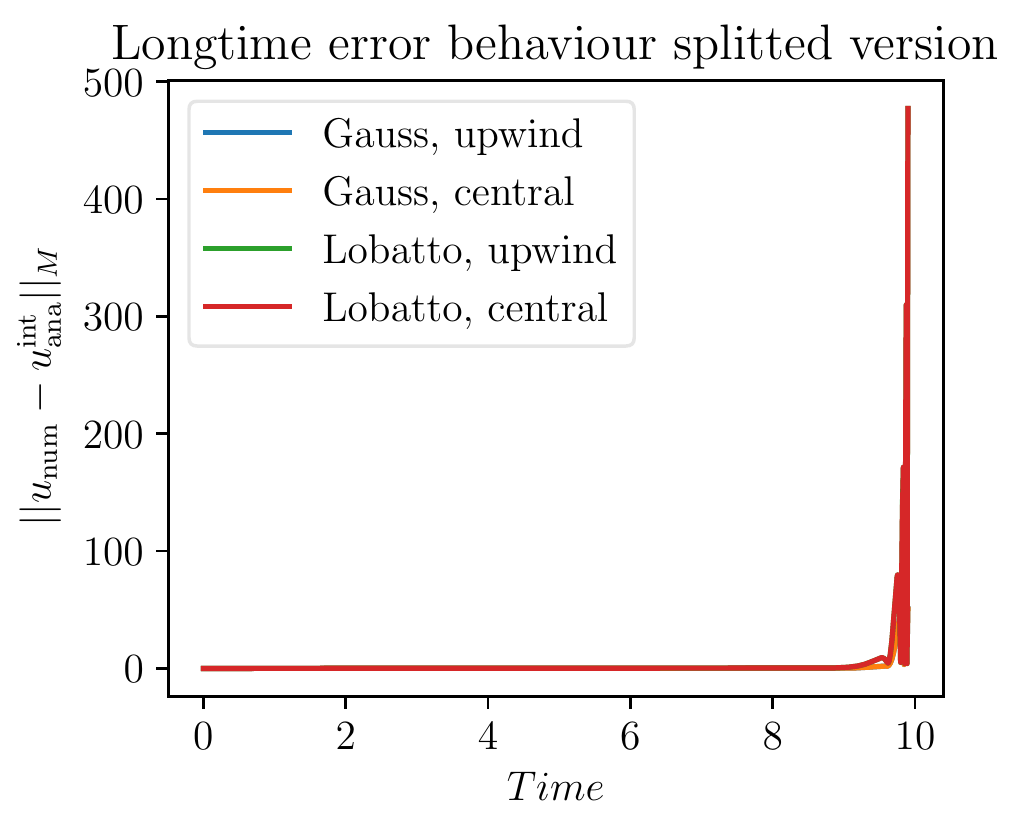}
    \caption{$N=3$, $K=40$, $t=9.9$.}
  \end{subfigure}
  \caption{Errors of discretisations for $a(x) = x^2$, $[x_L,x_R] = [-0.1,1]$,
           and $u_0(x) = \cos\left( \frac{\pi x}{2} \right)$.}
  \label{fig:flux_splitted_unsplit_gauss3}
\end{figure}

In the first simulation, the interval has been chosen as $[x_L,x_R] = [0.1,1]$ to guarantee the
positivity of the derivative of $a$ and also of its interpolation. Now, we change the interval
to $[x_L,x_R] = [-0.1,1]$, resulting in two major issues. First, the first derivative of $a$
is not strictly
positive anymore and the solution develops a pole at time $t=10$. Here, the solution is also
not uniformly bounded in its Sobolev norm and our error bounds \eqref{eq:Errorbound3}
and \eqref{eq:Errorbound4} do not hold. Nevertheless, the error behaviour can be investigated.
Using only the split discretisation for different times, we see in
\autoref{fig:flux_splitted_unsplit_gauss3}
that the errors increase and will increase further. They are unbounded.
Simultaneously, we also recognize
that the errors using Gauß-Legendre nodes still increase slower
due to the fact that the methods using
these nodes are more accurate.

\begin{figure}[!htp]
\centering
  \begin{subfigure}[b]{0.324\textwidth}
    \includegraphics[width=\textwidth]{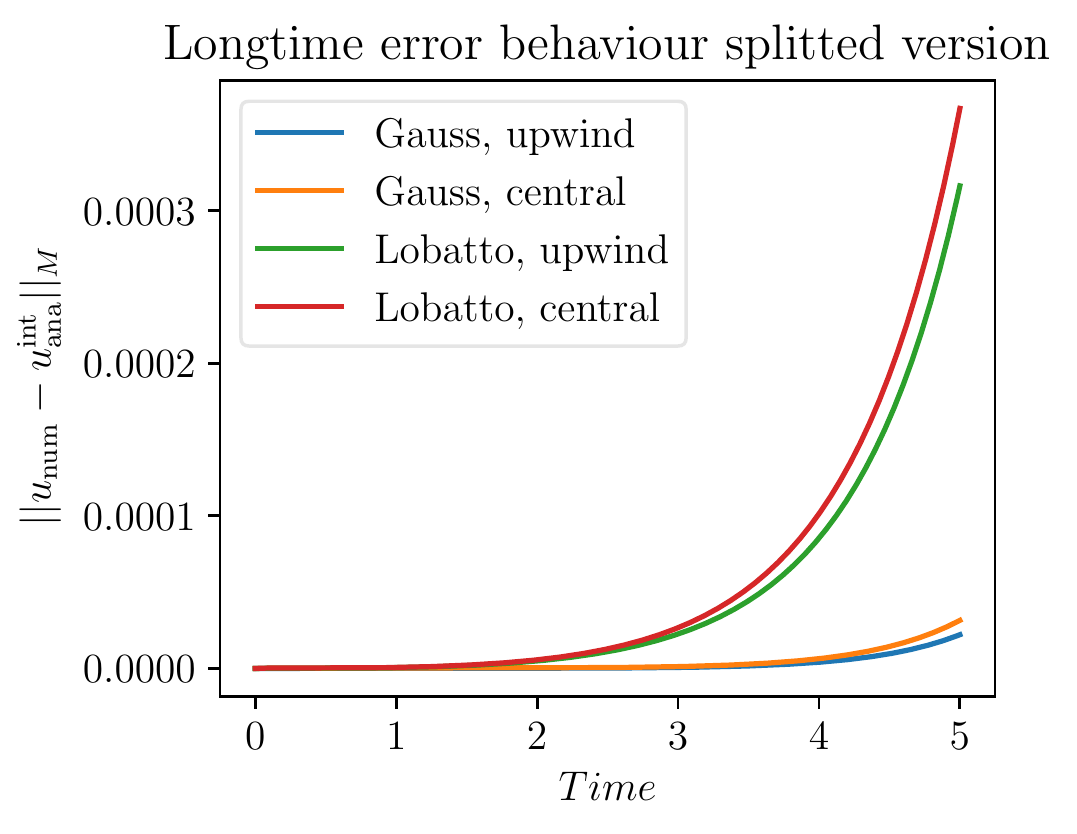}
    \caption{$N=3$, $K=40$, $t=5$.}
  \end{subfigure}%
  ~
  \begin{subfigure}[b]{0.324\textwidth}
    \includegraphics[width=\textwidth]{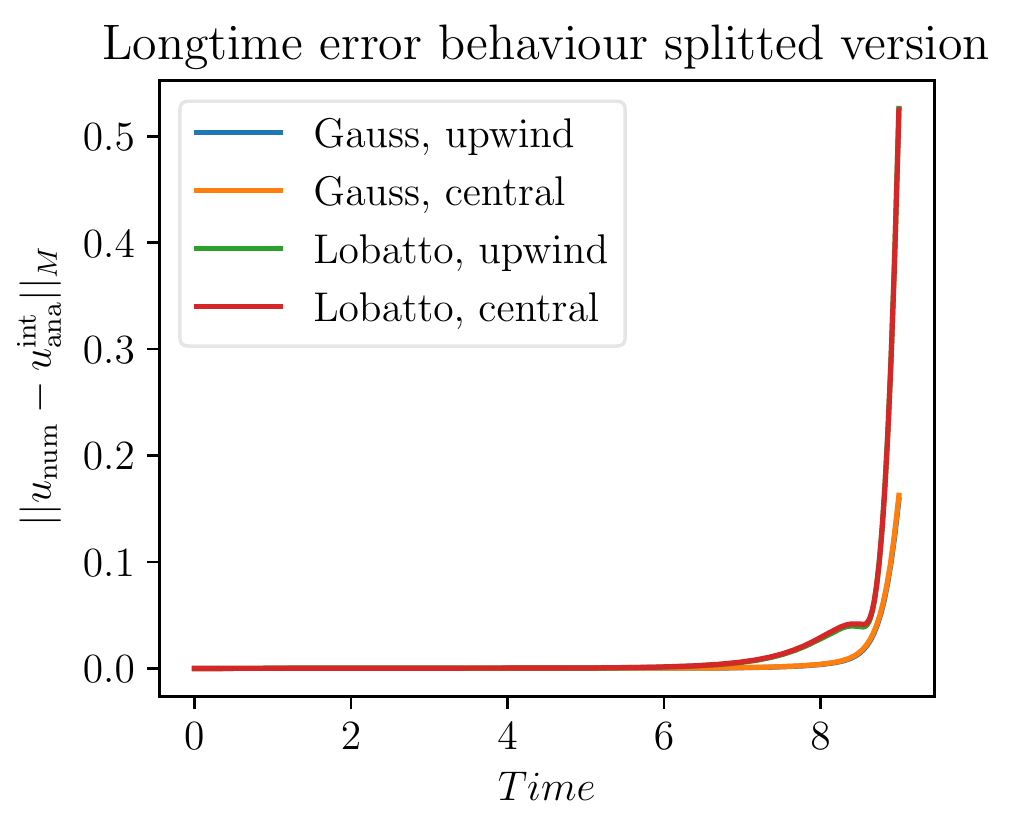}
    \caption{$N=3$, $K=40$, $t=9$.}
  \end{subfigure}%
    ~
  \begin{subfigure}[b]{0.324\textwidth}
    \includegraphics[width=\textwidth]{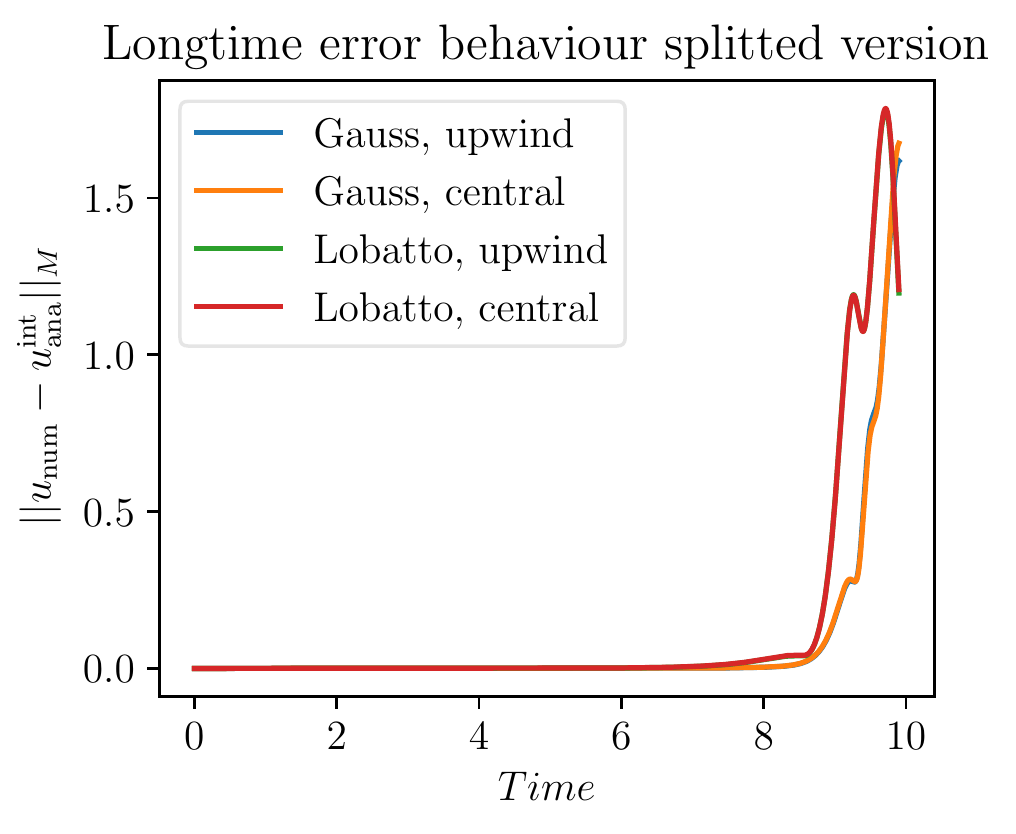}
    \caption{$N=3$, $K=40$, $t=9.9$.}
  \end{subfigure}
  \caption{Errors of discretisations for $a(x) = x^2$, $[x_L,x_R] = [-0.1,1]$, and $u_0(x) = \exp(-x^4)$.}
  \label{fig:flux_splitted_unsplit_gauss3_3}
\end{figure}

Furthermore, by changing the initial condition to $u_0(x)=\exp(-x^4)$ instead of
$u_0(x) =\cos\left( \frac{\pi x}{2} \right)$, we are able to avoid the pole in the solution
\eqref{sol_a_x_x} since the exponential function will tend fast enough to zero compared to
$(1+tx)^2$ and we can extend the solution. Nevertheless we get further problems here.
If we have a look on the error behaviour in \autoref{fig:flux_splitted_unsplit_gauss3_3},
we see that we get a similar increase of the errors like in \autoref{fig:flux_splitted_unsplit_gauss3},
but they are much smaller. Nevertheless they are still unbounded, but why do we have this behaviour?
The analytical solution is for fixed times bounded, nevertheless we demand
as one assumption right at the beginning at equation \eqref{eq:Model_problem} the solution to be
\textbf{uniformly} bounded in time. However, this is not the case anymore. This demonstrates
again how essential this assumption is.

The same issue arises if we are investigate $a(x)=\cosh(x)+1$ as in \cite{ranocha2018generalised}.
Therefore, we skip this case here.

\subsection{Coefficient \texorpdfstring{$a(x) = 1-x^2$}{a(x) = 1-x²}}\label{subsec:test_1_x_x}

Here, we choose the coefficient $a(x) = 1 - x^2$. The solution of the Cauchy problem is
\begin{equation}\label{solutions_ax_1_x}
  u(t,x)
  =
  \frac{u_0\bigl( (-x \cosh(t) + \sinh(t)) / (x \sinh(t) - \cosh(t)) \bigr)}
       {(\cosh(t) - x \sinh(t))^2}.
\end{equation}
Using the domain $[x_L,x_R] = [-1, 0.9]$ and the initial condition $u_0(x) = \sin(\pi x)$,
the solution remains bounded but $a'(x) < 0$ for $x > 0$.
\begin{figure}
\centering
  \includegraphics[width=0.5\textwidth]{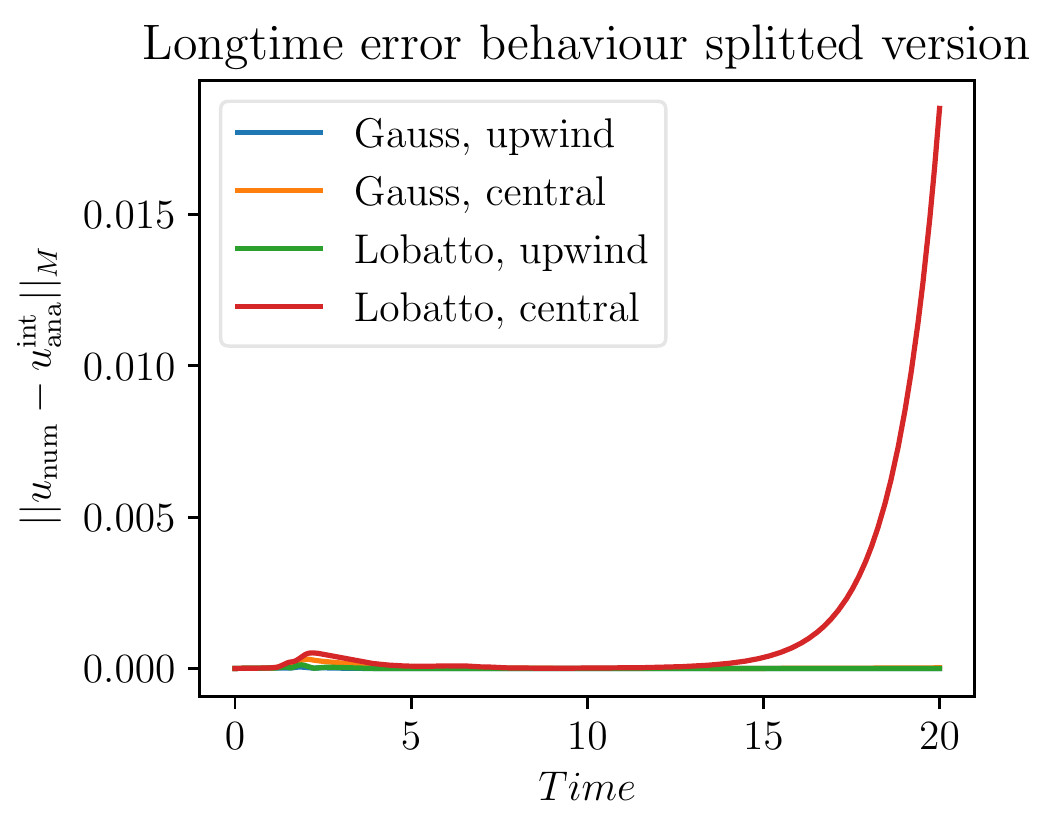}
  \caption{Error behaviour for $a(x) = 1 - x^2$, $t=20$.}
  \label{fig:flux_splitted_no_convergece}
\end{figure}
If we investigate now the long time error behaviour, we get a huge  increase of the errors
if we apply the central fluxes, cf. \autoref{fig:flux_splitted_no_convergece}. This matches
perfectly our theoretical investigations in \autoref{sec:Error_Lobatto} and
\autoref{sec:Error_Gauss}, cf. Remark~\ref{remark_int}.
We explain the reasons again in detail in the next test case and a physical interpretation
and illustration is given afterwards.

\subsection{Coefficient \texorpdfstring{$a(x) = \cos(x)$}{a(x) = cos(x)}}\label{subsec:a(x)=cos(x)}

Here, we choose $a(x) = \cos(x)$ and $u_0(x) = \sin(5 x)$. The solution of the Cauchy problem is
\begin{equation}\label{eq:Test_a_x_cos_x}
\begin{aligned}
  u(t,x)& = u_0\bigl( x_0(t,x) \bigr) \frac{\cos \bigl( x_0(t,x) \bigr) }{\cos (x)},
  \\
  x_0(t,x)
  &=
  - 2 \arctan\bigl( \tanh\bigl( t/2 - \artanh( \tan(x/2) ) \bigr) \bigr).
\end{aligned}
\end{equation}
We can find an interval for our solution \eqref{eq:Test_a_x_cos_x} so that $a'(x)\leq 0$
and $u(t,x)$ does not blow up, e.g. $[x_L, x_R] = [0.1, \pi/3]$. The solution remains
bounded but $a'(x) < 0$.

 \begin{figure}[!htp]
\centering
  \begin{subfigure}[b]{0.332\textwidth}
    \includegraphics[width=\textwidth]{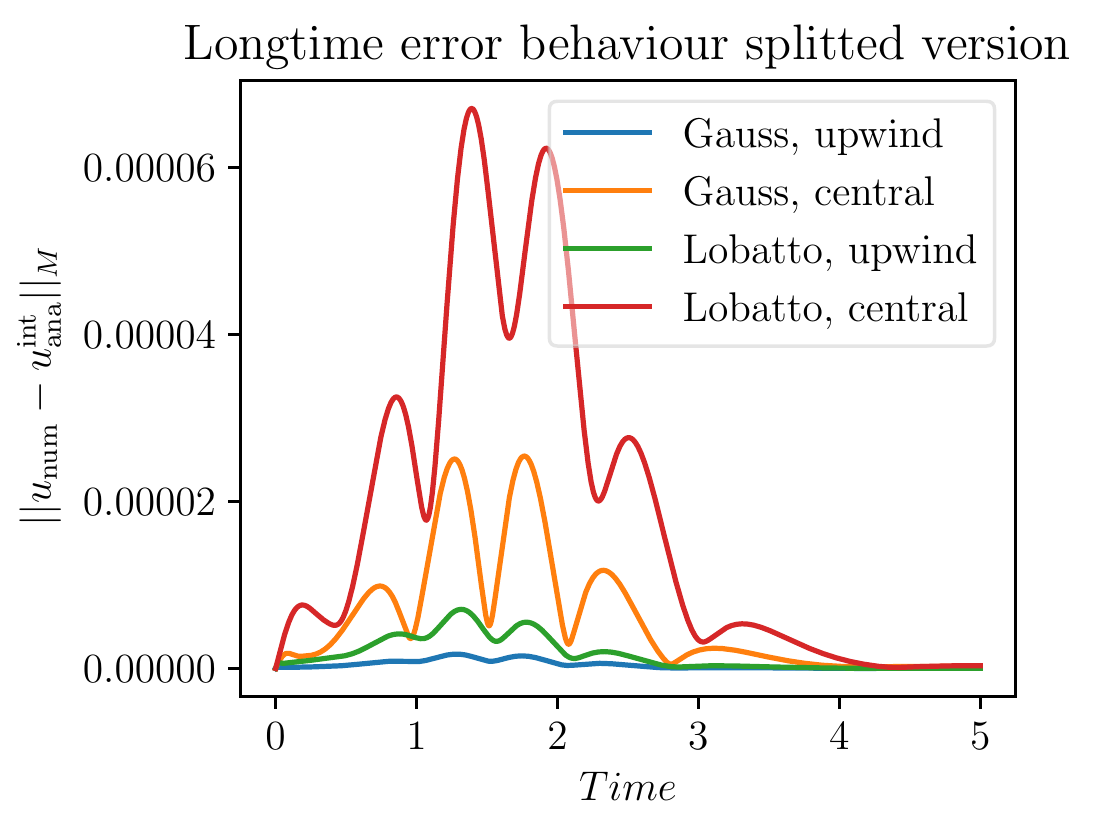}
    \caption{$N=3$, $K=30$, $t=5$.}
  \end{subfigure}%
   ~
  \begin{subfigure}[b]{0.332\textwidth}
    \includegraphics[width=\textwidth]{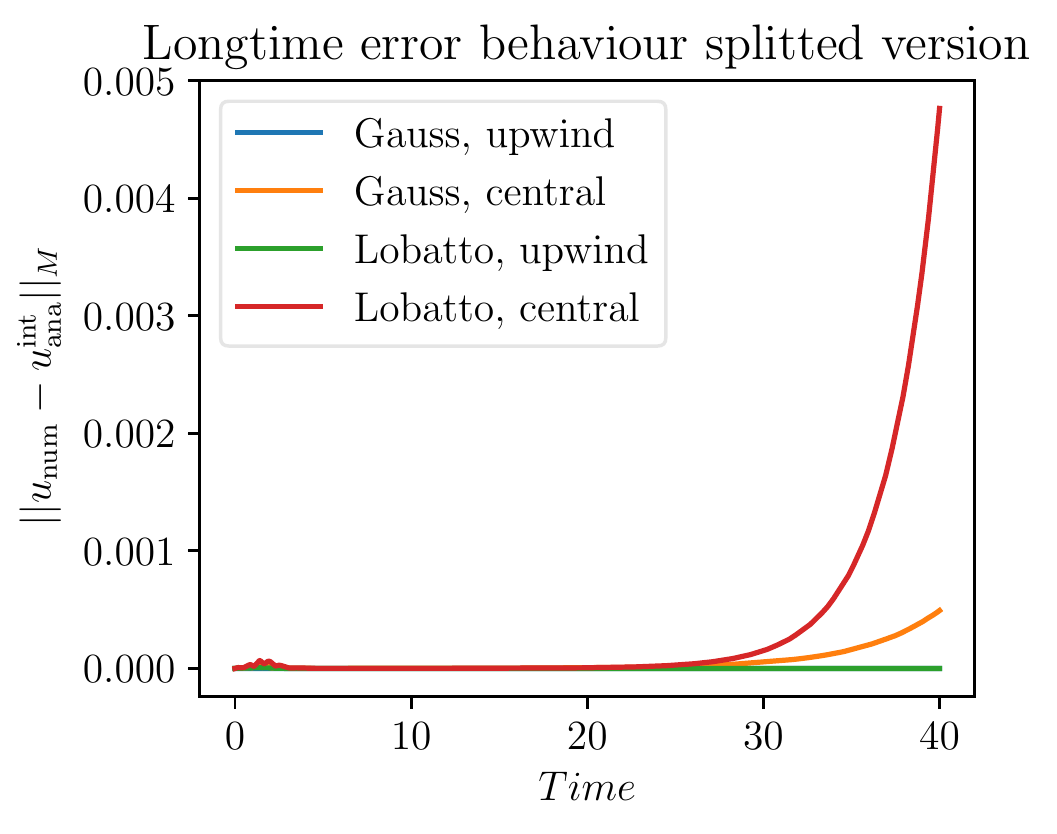}
    \caption{$N=3$, $K=30$, $t=40$.}
  \end{subfigure}%
  ~
  \begin{subfigure}[b]{0.332\textwidth}
    \includegraphics[width=\textwidth]{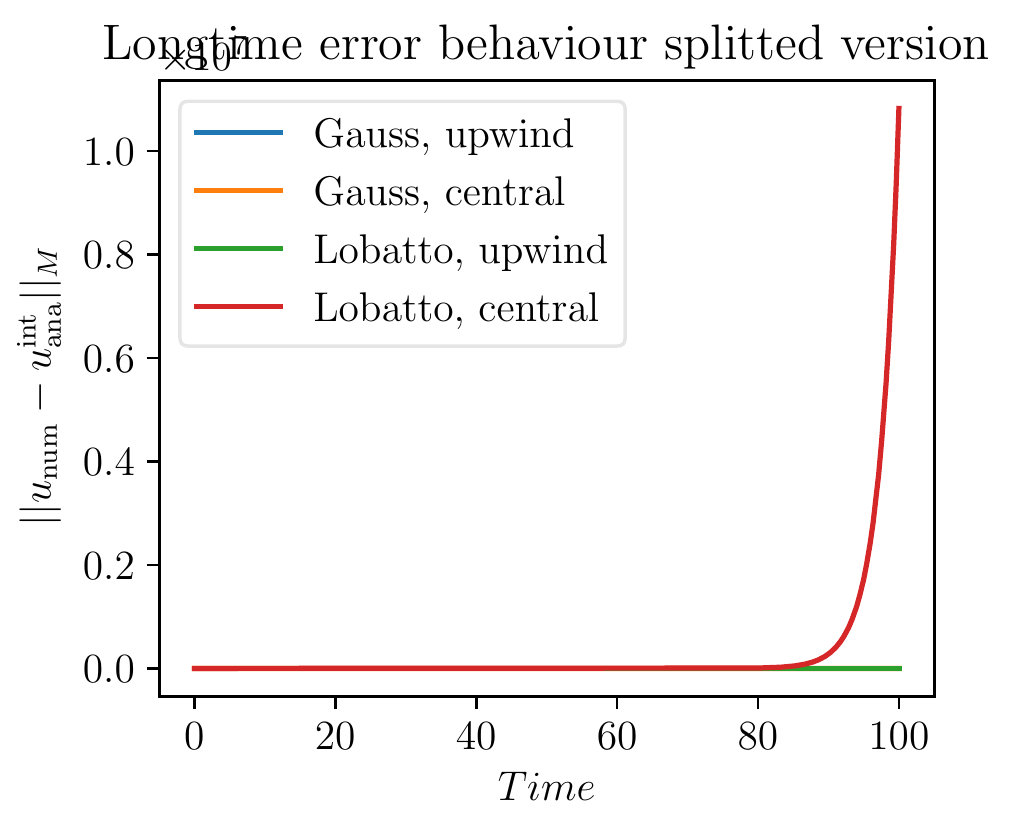}
    \caption{$N=3$, $K=30$, $t=100$.}
  \end{subfigure}%
  \caption{For $a(x) = \cos(x)$, the slope of $a$ is negative. The errors of numerical solutions
           using the central flux increase, whereas the upwind flux results in bounded errors.}
  \label{fig:flux_splitted_unsplit_gauss4}
\end{figure}

In \autoref{fig:flux_splitted_unsplit_gauss4}, we see the behaviour of the error for different times.
First, one may suppose that the error remains bounded in time, but this is not the case as can be seen
stepping further in time.
Using the central fluxes ($\sigma=0$), the $BTs$ terms are zero and we do not find an $\eta$
which is bounded with a positive constant from below away from zero.
One may recognize also that for Gauß-Legendre nodes, the error increases much slower
(second picture).
Surely, one reason for this is the smaller error in the Gauß-Legendre case. Furthermore, also
the term $\Theta$ may have a positive impact of the error behaviour.

\begin{figure}[!htp]
\centering
  \begin{subfigure}[b]{0.49\textwidth}
    \includegraphics[width=\textwidth]{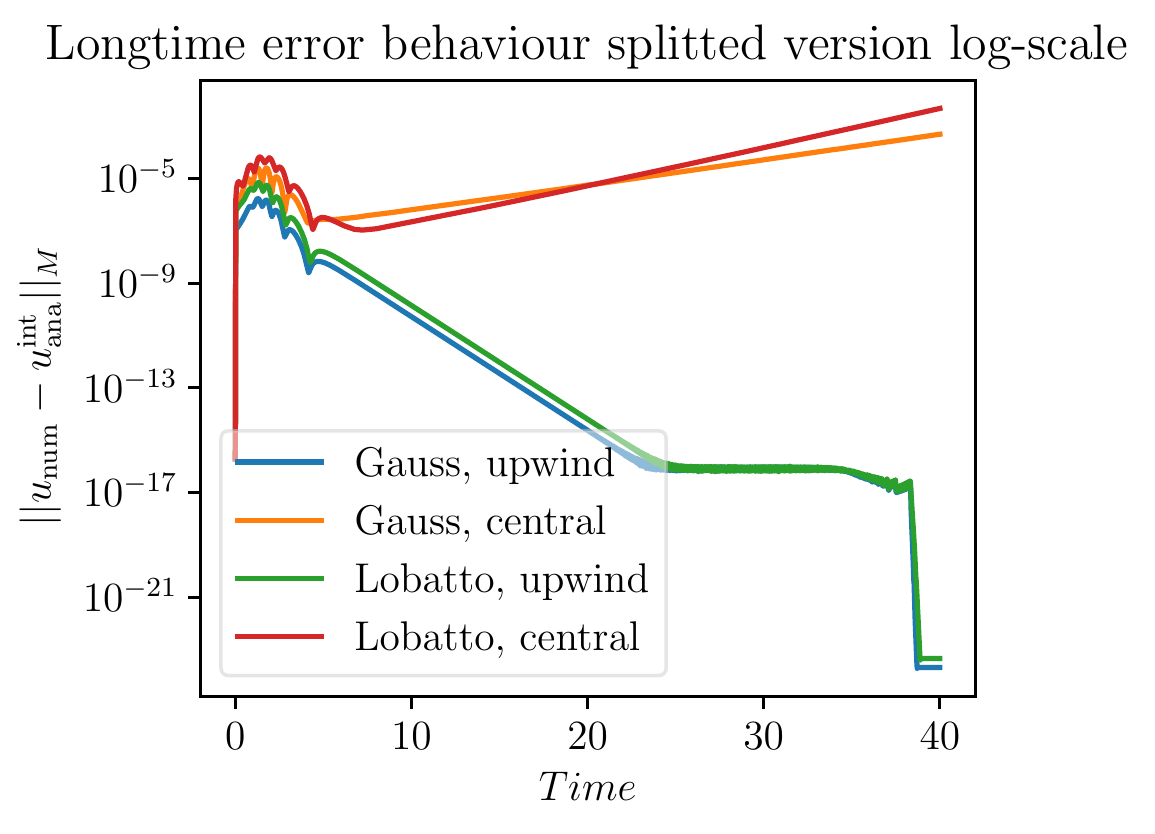}
    \caption{$N=3$, $K=30$, $t=40$.}
  \end{subfigure}%
  ~
  \begin{subfigure}[b]{0.49\textwidth}
    \includegraphics[width=\textwidth]{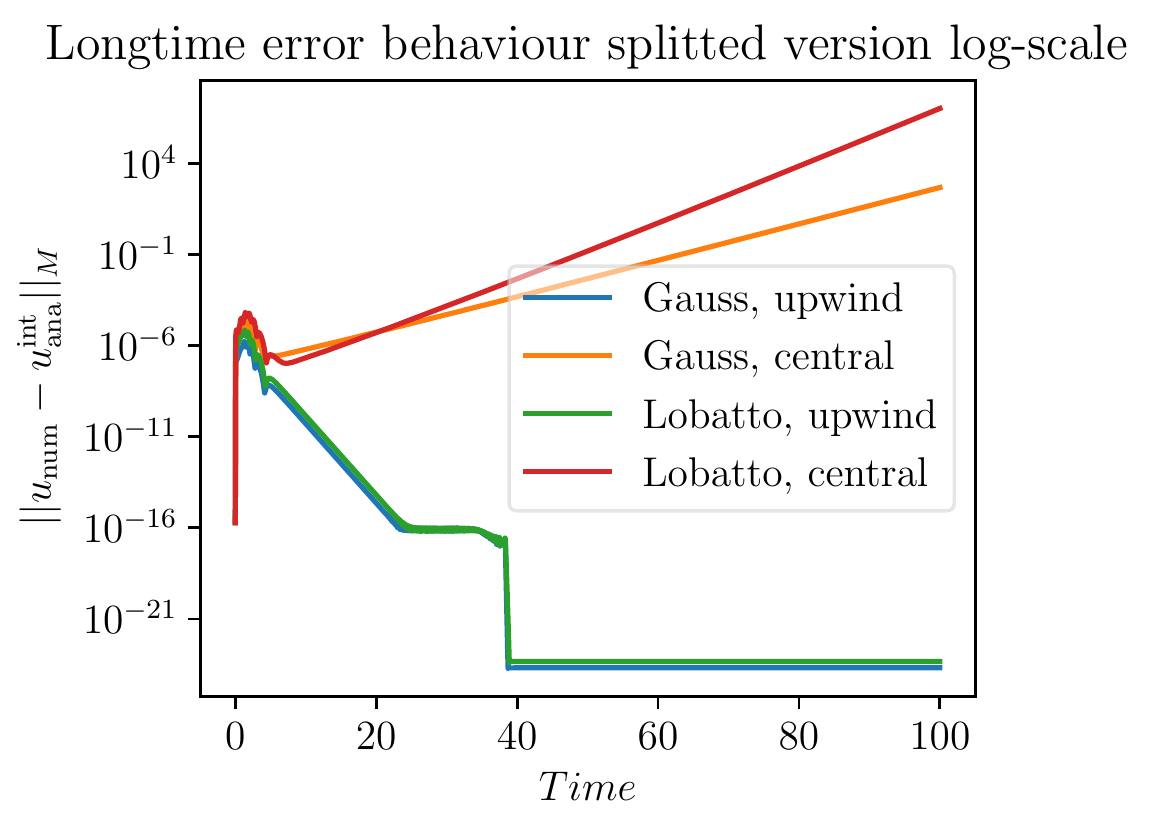}
    \caption{$N=3$, $K=30$, $t=100$.}
  \end{subfigure}%
  \caption{Error of numerical solutions for $a(x) = \cos(x)$ in logarithmic scale.}
  \label{fig:flux_splitted_unsplit_gauss4_log}
\end{figure}

However, this example demonstrates well that the condition $a'(x)>0$ is essential for the
boundedness of the error, also in the test case of section~\ref{subsec:test_1_x_x}. One can
rescue \eqref{eq:Errorbound3} and \eqref{eq:Errorbound4} by applying the upwind flux like it can be seen
in this test case and especially in \autoref{fig:flux_splitted_unsplit_gauss4_log}.

\begin{figure}
  \centering
  \captionsetup{skip=0pt}
  \includegraphics[width=\textwidth]{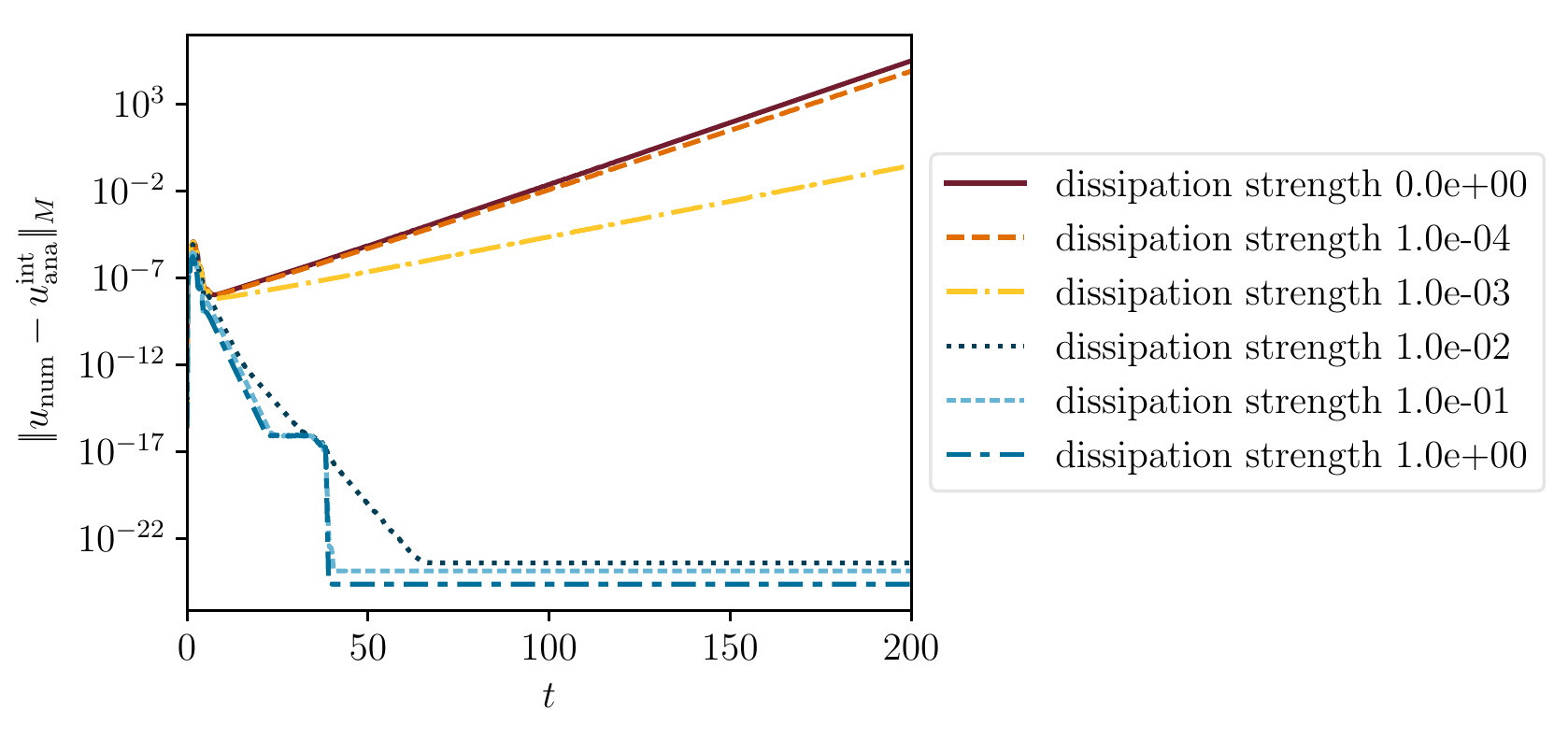}
  \caption{Errors of numerical solutions using SBP finite difference schemes.}
  \label{fig:SBP_FD__a_cos__order_8__N_100__split_true}
\end{figure}

By applying an SBP-SAT finite difference scheme with one block, the internal terms $BTs$
do not exist. Using the SBP difference operator of \cite{mattsson2004summation} with interior
order of accuracy eight, the split form, and \num{100} nodes for this problem, the error
is unbounded, as can be seen in \autoref{fig:SBP_FD__a_cos__order_8__N_100__split_true}.
However, if the high-order artificial dissipation operator of \cite{mattsson2004stable} is
applied additionally, the error remains bounded.

Comparing Figures~\ref{fig:flux_splitted_unsplit_gauss4_log} and
\ref{fig:SBP_FD__a_cos__order_8__N_100__split_true} demonstrates that stabilisation
induced by upwind fluxes or artificial dissipation operators is crucial and comparable.
Furthermore, Gauss-Legendre nodes not including boundary points provide some stabilisation.

\subsection{A first analytical study}\label{subsec:analytical}
\revsec{
As can be seen in Figures \ref{fig:flux_splitted_unsplit_gauss4_log}
and \ref{fig:SBP_FD__a_cos__order_8__N_100__split_true}, if $a'(x)$ is not positive
the long time errors show different behaviours depending on the dissipation which is
added to the scheme by numerical fluxes or artificial dissipation terms.
Here, we give a short rough analysis on this topic under what conditions we can guarantee boundedness.
A more detail analysis should follow in future research with more validations. \\
We are starting considering  $\eta_G(t)$ from \eqref{eq:eta_def}.
It is
\begin{equation}\label{eq:eta_analysis}
 \eta(t):=\frac{BTs +Int_d+\Theta_2}{||\epsilon_1||^2_N}
\end{equation}
with $Int_d:= \frac{1}{2} \sum_{k=1}^K \frac{\Delta x_k}{2} \l(
\partial_x \vec{a}^k, \mat{\epsilon}_1^k
  \vec{\epsilon}_1^k \r)_N$.
A sufficient condition for the mean of $\eta(t)$  to be positive is that every value of $\eta(t)$  is positive.
Therefore, we require
\begin{equation*}
 \frac{BTs +Int_d+\Theta_2}{||\epsilon_1||^2_N} >0.
\end{equation*}
If the derivative of $a$ is negative, we can reformulate the inequality above as
\begin{equation*}
 \left( BTs +\Theta_2 \right)\frac{1}{||\epsilon_1||^2_N} >
 \frac{1}{2||\epsilon_1||^2_N} \left( \sum_{k=1}^K \frac{\Delta x_k}{2}
 \l(   \left|\partial_x \vec{a}^k \right|,   \mat{\epsilon}_1^k
  \vec{\epsilon}_1^k  \r)_N\right),
\end{equation*}
and even strengthen our assumptions by requiring
\begin{equation}\label{eq:condtion_on_a}
\begin{split}
  &\left( BTs +\Theta_2 \right)\frac{1}{||\epsilon_1||^2_N} >
 \frac{\max_{x\in (x_0,x_K)} \left|\partial_x \vec{a} \right| }{2||\epsilon_1||^2_N}
 \left( \sum_{k=1}^K \frac{\Delta x_k}{2} \l(    \vec{1}, \mat{\epsilon}_1^k
  \vec{\epsilon}_1^k \r)_N \right)\\
 \text{ or }&\left( BTs +\Theta_2 \right)\frac{1}{||\epsilon_1||^2_N} >
 \frac{\max_{x\in (x_0,x_K)} \left|a'(x) \right| }{2} .
\end{split}
\end{equation}
From \eqref{eq:condtion_on_a} we realize the $BTs$-terms are responsible to guarantee
that this sufficient condition is fulfilled. In case of a central numerical flux, $BTs\equiv0$
and we have to add additional dissipation to the scheme as it is done in the SBP-SAT schemes in
Figure \ref{fig:SBP_FD__a_cos__order_8__N_100__split_true}. However, also the dependence of the error is
important and we may also realize that in case of using Gau\ss-Legendre we rather get the condition
\eqref{eq:condtion_on_a} fulfilled. However, this estimation is rough and should be improved in
further research.
}
\subsection{Physical Interpretation and Illustration}
\label{sec:phys-interpretation}

In order to understand some results better, a physical interpretation of the
advection equation can be used. This serves also as illustration and explains
the rational behind some of the choices regarding for example the numerical
experiments.

The advection equation $\partial_t u + \partial_x (a u) = 0$ with non-negative
velocity $a(x)$ is a conservation law with varying coefficients. Thus, the total
mass $\int u$ is conserved and $u$ is transported from left to right due to
$a(x) \geq 0$. In order to compute analytical solutions of the Cauchy problem,
the method of characteristics can be used, cf. \cite[Chapter 3]{bressan2000hyperbolic}.
\begin{itemize}
  \item
  Solve the ODE $x'(t) = a\bigl( x(t) \bigr)$, $x(0) = x_0$, for $x(t) = x(t; x_0)$.
  Compute also the inverse function $x_0 = x_0(t; x)$.

  \item
  Solve the ODE $z'(t) = - a'\bigl( x(t;x_0) \bigr) z(t)$, $z(0) = z_0$, for
  $z(t) = z(t; z_0, x_0)$.

  \item
  Set $z_0 = z_0(x_0) = u_0(x_0) = u_0\bigl( x_0(t;x) \bigr)$ and obtain the analytical
  solution $u(t,x) = z(t; z_0, x_0) = z(t; z_0\bigl(x_0(t;x)), x_0(t;x)\bigr)$.
\end{itemize}
In the second step, if $a' > 0$, the absolute value of $z(t)$ decreases. Contrary,
if $a' < 0$, the absolute value of $z(t)$ increases. This corresponds directly
to the physical interpretation as transport problem. Since $u$ is conserved and
transported with velocity $a(x)$, there is a loss of $u$ if $a' > 0$, since
there is less new mass coming from the left than going to the right. Similarly,
$a' < 0$ yields an increase of $u$, since more mass is coming from the left than
transported to the right. This explains also the critical role of $a'(x)$.
If $a' < 0$, there can be blow-up phenomena in the solution $u$, resulting in
possibly finite life spans and increasing energies and errors of numerical
solutions. If $a' > 0$, this cannot happen.

If one wants to investigate a situation with $a'(x) > 0$ in some parts and
$a'(x) < 0$ in other parts of the domain, there are basically two possibilities.
Firstly, there can be a local minimum of $a(x)$, e.g. for $a(x) = x^2$.
In this case, there can be a blow-up of the solution $u$,
since more mass is coming from the left than transported to the right at this
minimum. However, this blow-up phenomenon caused by the varying transport velocity
$a(x)$ can be balanced by the initial condition $u_0$. If there is simply not
enough mass on the left, than the higher transport speed there can not cause a
blow-up of the solution $u$. This explains our choice of the intervals and the
initial conditions for these cases.

Secondly, there can be a local maximum of $a(x)$, e.g. for $a(x) = 1 - x^2$ or
$a(x) = \cos(x)$. Now, there is no blow-up at the critical point, since more
mass is transported to the right. However, both examples have stagnation points
with $a(x) = 0$. At such points, there will be a  blow-up of the solution, since
mass is coming from the left but not transported to the right. In order to
avoid this phenomenon of the Cauchy problem, the interval can be chosen adequately,
i.e. bounded away at the right from the point with $a(x) = 0$. Then, the blow-up
of the solution of the Cauchy problem does not cause any problems for the corresponding
solution of the initial value problem. This explains our choices of the domains
for these cases.

%% file: 7_Generalisation.tex
\section{Possible Generalisation and Examples}\label{sec:Generalisation}
\revsec{
As has been demonstrated hitherto, the error of numerical solutions of scalar
hyperbolic conservation laws with varying coefficients does not necessarily remain
bounded in finite domains, contrary to the expectation for linear systems with
constant coefficients. Here, some further remarks concerning generalisations of
this result are given.
}
{\color{black}
\subsection{Linearized Euler Equations}

We start by considering the theory for the linearized Euler Equations which are one of the  most --- if not the most ---
investigated system in computational fluid dynamics.
The one-dimensional compressible Euler equations in conservation form are
\begin{equation}\label{eq:Euler_equation}
\partial_t \bU+ \partial_x \bF(\bU) =0,
\end{equation}
where $ \bU$ is the state vector of the conserved quantities and $\bF$ is the flux.
Thus,
\begin{equation}\label{Eq:Euler}
 \bU=\begin{pmatrix}
      \rho \\
       m\\
       E
     \end{pmatrix},
     \quad
     \bF(\bU)=\begin{pmatrix}
               m\\
               \rho u^2 +p\\
               u(E+p)
              \end{pmatrix},
\end{equation}
where $\rho$ is the mass density, $m=\rho u$ is the momentum, $E$ is the total energy,
$u$ is the velocity and $p$ is the pressure related to
$\bU$ by the equation of state $p=(\gamma-1)(E-\rho \frac{u^2}{2})$ using $\gamma$ for
the specific heat capacities.
We can rewrite \eqref{eq:Euler_equation} as
\begin{equation*}
 \partial_t \bU+ \bf{A}(\bU) \partial_x \bU =0,
\end{equation*}
where $\bf{A}=\partial_{\bU} \bF$ is the Jacobian matrix which has only real eigenvalues and
can be diagonalized by the matrix $\bf{R}$ of eigenvectors.
Indeed, $\bf{A}=\bf{R}\Lambda\bf{R}^{-1}$, where $\Lambda=\diag{\lambda_1,\lambda_2,\lambda_3}
=\diag{u+c,u,u-c}$. Here, $c$ is the speed of sound which satisfies
\begin{equation}
  c(\rho)^2 = p'(\rho) > 0.
\end{equation}
As mentioned before, a lot of investigations of \eqref{eq:Euler_equation} can be
found in the literature where also different linearization techniques were used
depending on the numerical schemes
\cite{roe1981approximate, fey1998multidimensional, steger1981flux, van1997flux}.
Here, we will focus on this topic and the problems which can appear.
This yields us to some outlook for future research.

%

\subsubsection*{Linerization around a Smooth Solution --- An Outlook}

We are not considering  the full system  \eqref{Eq:Euler}  for a smooth solution,
but the truncated/simplified/shortened version \cite{gustafsson2013time}
\begin{equation}\label{eq:Short_Euler}
\begin{aligned}
  \partial_t \rho + \partial_x (\rho u) &= 0,
  \\
  \partial_t u + u\partial_x u + \frac{1}{\rho} \partial_x p(\rho) &= 0,
\end{aligned}
\end{equation}
to explain the problem. Using a Taylor series approach for the linearization around
a smooth solution $(\hrho,\; \hu)$ yields a linear system with variable coefficients
of the form
\begin{equation*}
  \partial_t \vect{\rho \\ u}
  +
   \begin{pmatrix}
     \hu & \hrho
    \\
     \frac{c(\hrho)^2}{\hrho} & \hu
   \end{pmatrix}
   \partial_x \vect{\rho \\ u}
  + C  \vect{\rho \\ \hu}
  =
  0,
\end{equation*}
where $C$ depends on $(\hrho,\; \hu)$ and their derivatives  such that $C=0$ if $\hrho$ and $\hu$ are constant.
This system can be symmetrized using $\rho_S := \frac{c(\hrho)}{\hrho} \rho$, resulting in
\begin{equation}
  \partial_t \vect{\rho_S \\ u}
  +
  \begin{pmatrix}
    \hu & c(\hrho)
    \\
    c(\hrho) & \hu
  \end{pmatrix}
  \partial_x \vect{\rho_S \\ u}
  + \tilde C \vect{\rho_S \\ u}
  =
  0,
\end{equation}
where $\tilde C$ depends on $(\hrho,\; \hu)$ and their derivatives such that $\tilde C = 0$
if $\hrho$ and $\hu$ are constant. If we have constant coefficients, this investigation belongs
to the case which was already studied in \cite{oeffner2018error, nordstrom2007error, kopriva2017error}
and the error remains bounded under the conditions give there.
Otherwise, all entries of $\tilde C$ are non-trivial and the equations cannot be decoupled.
Already for symmetric systems, we get further problems depending on the estimation of the
energy growth, as described in section \ref{sec:Mulitdimensional_systems}.
The investigation of the error behaviour for this problem is not straightforward and
should be considered in more  detail in future work.

\begin{remark}
As mentioned above, there are different techniques for linearizing
the Euler equations. They  depend on the numerical schemes which are used/constructed
for these system.
Here, we  only mention the approach by Roe \cite{roe1981approximate} about   flux difference splitting
or the flux vector splitting in \cite{steger1981flux}. The linearization is used in the construction
of the numerical schemes in some sense. To follow their ideas together with our analysis about the long time error behaviour
is  an alternative ansatz and will also be considered in future research.
\end{remark}

\subsection{Multidimensional Systems} \label{sec:Mulitdimensional_systems}

We consider the linear magnetic induction equation
\begin{equation}
\begin{aligned}
  \partial_t B(t,x) &= \nabla \times \bigl( u(t,x) \times B(t,x) \bigr),
  && t \in (0,50), x \in (0,1)^3,
  \\
  B(0,x) = u(t,x)
  &=
  \begin{pmatrix}
    \sin(\pi x) \cos(\pi y) \cos(\pi z) \\
    \cos(\pi x) \sin(\pi y) \cos(\pi z) \\
    - 2 \cos(\pi x) \cos(\pi y) \sin(\pi z)
  \end{pmatrix},
  && x \in [0,1]^3,
\end{aligned}
\end{equation}
supplemented with the divergence constraint $\div B(t,x) = 0$, cf. \cite{mishra2010stability,
koley2009higher}. This specific example is taken from \cite{ranocha2018numerical}.
Here, $B$ is the magnetic field and $u$ the particle velocity.
Since $u$ vanishes at the boundary of the domain, no boundary condition is specified.
In order to get a symmetric hyperbolic system, the nonconservative
source term $-u \div B$ is added to the right hand side, resulting in an energy
estimate if a splitting is used as described in the references listed above.
There are several discrete forms of the equation allowing an energy estimate
\cite{ranocha2018numerical}. Using the terminology introduced there, the most
obvious one uses the same split form as applied at the continuous level and is
called (product, central, split). Another choice described there is
(central, central, central). The implementations of \cite{ranocha2018induction}
are used in the following.

Applying both discretisations, SBP FD operators of interior order of accuracy $4$,
and $40^3$ nodes to discretize the domain yields the results visualized in
\autoref{fig:SBP_FD__IEQ}. As can be seen there, the form (product, central, split)
results in an exponential growth of both the energy and the error while the other
form yields a bounded error. Adding artificial dissipation does not change the
result significantly.

\begin{figure}
  \centering
  \captionsetup{skip=2pt}
  \begin{subfigure}{0.8\textwidth}
    \includegraphics[width=\textwidth]{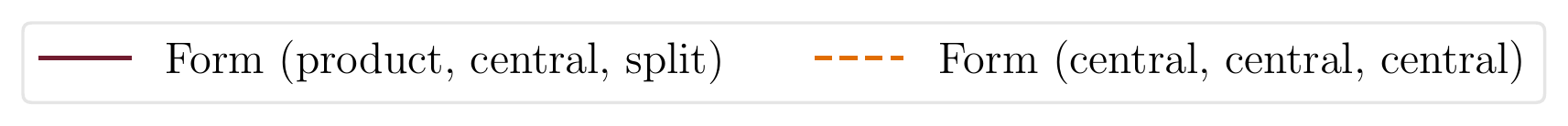}
  \end{subfigure}%
  \\
  \begin{subfigure}{0.495\textwidth}
    \captionsetup{skip=0pt}
    \includegraphics[width=\textwidth]{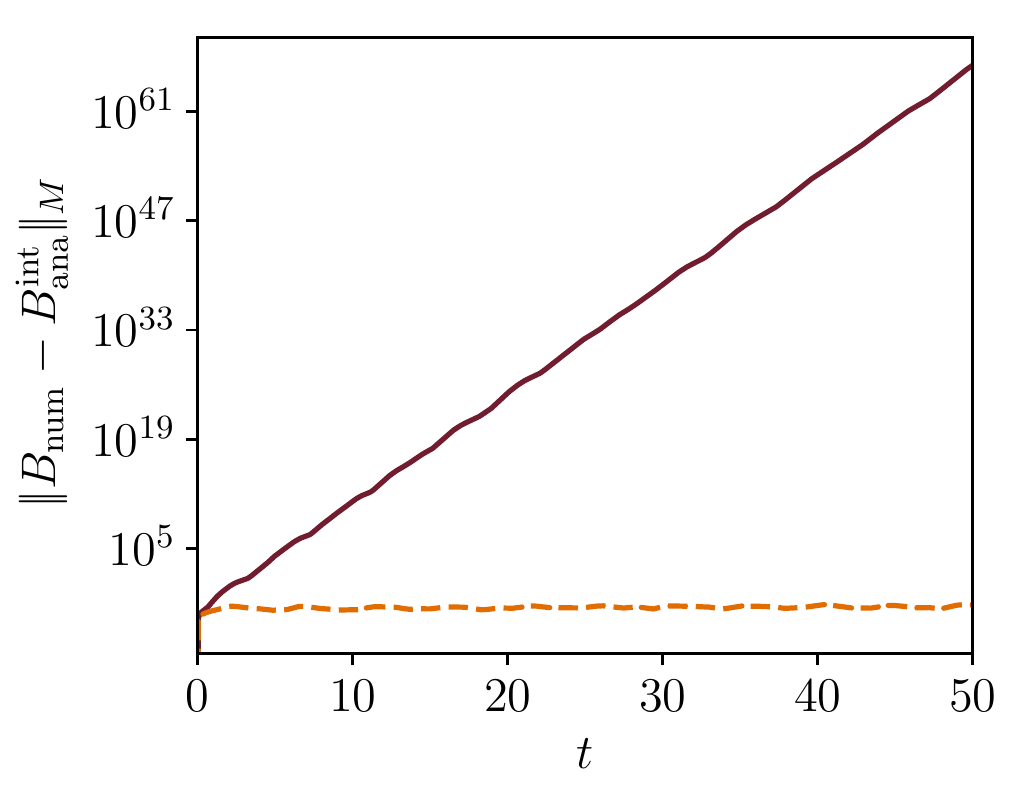}
    \caption{Error.}
  \end{subfigure}%
  ~
  \begin{subfigure}{0.495\textwidth}
    \captionsetup{skip=0pt}
    \includegraphics[width=\textwidth]{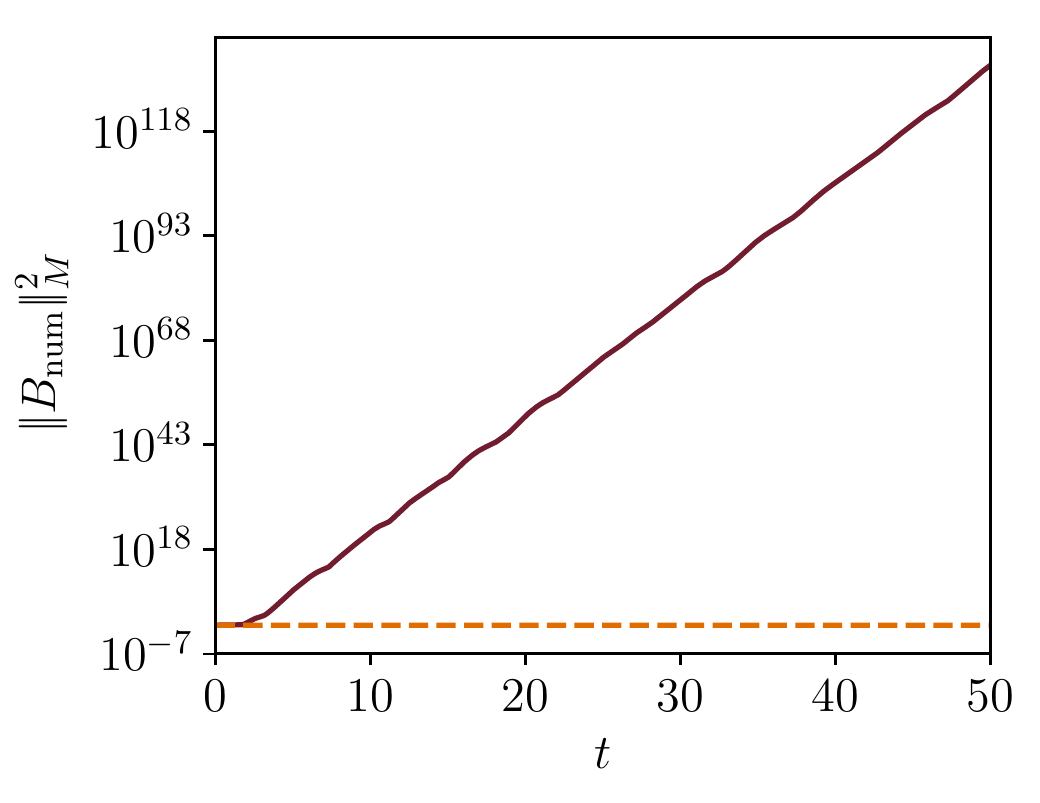}
    \caption{Energy.}
  \end{subfigure}%
  \caption{Errors and energies of numerical solutions of the induction equation.}
  \label{fig:SBP_FD__IEQ}
\end{figure}

These results are in accordance with the energy estimates (an exponential growth
is allowed as worst case estimate) and the investigations in this article. The
main complications for (product, central, split) are presumably a combination of
\begin{itemize}
  \item
  The velocity $u$ vanishes at the boundary and errors cannot be transported
  out of the domain; instead, they accumulate.

  \item
  While the analytical solution has a bounded energy, the worst case estimate
  allows an exponential growth.

  \item
  The analytical solution is a steady state which is not necessarily represented
  exactly by the discretisation.
\end{itemize}

This shows that severe problems can be expected for general symmetric hyperbolic
systems with varying coefficients in multiple space dimensions.
}

%% file: 8_Summary.tex
\section{Summary and Discussion}
\label{sec:summary}

In this article, we have conducted an analysis of the long-time behaviour of the error
of numerical solutions to the linear advection equation with variable coefficients
in bounded domains. Using flux reconstruction schemes/\revPP{discontinuous Galerkin methods} with summation-by-parts operators,
we provide a detailed analysis of the influence of both the choice of the numerical flux
and the polynomial basis. If boundary conditions are imposed in a provably stable way using
numerical fluxes, the error can be bounded uniformly in time, depending on the variable
coefficient $a(x)$ and the numerical fluxes at the interior boundaries. However, there
can be also an unbounded growth of the error if certain conditions are not satisfied.

Firstly, if the varying coefficient $a(x)$ behaves nicely, inducing a decay of the analytical
solution, the long time behaviour of the numerical error is comparable to the case of
constant coefficients. The application of upwind fluxes at interior boundaries results
in a smaller asymptotic value of the error and this value is also attained faster. Using
Gauß-Legendre nodes results in smaller errors compared to Gauß-Lobatto-Legendre nodes.

However, if the varying coefficient $a(x)$ induces a possible growth or blow-up of the analytical
solution, the situation is totally different. Of course, if the solutions blows up in
finite time, so does the error. This behaviour is not possible for constant coefficients.
Moreover, there can still be some problems, even if the solution does not blow up.
Indeed, the variable coefficients can trigger a growth of the error that has to be balanced
by additional stabilisation such as upwind numerical fluxes compared to central ones or
artificial dissipation, e.g. in finite difference methods. We have explained this behaviour
and have presented several numerical examples, where upwind numerical fluxes or artificial
dissipation result in uniformly bounded errors while the errors increase without bound if
central numerical fluxes or no additional dissipation operators are applied.

\revsec{Finally, in the last section we have extended our analysis of the long time error behaviour
to systems. Here, several problems emerge and we have given an outlook for further research topics in this context focussing
on coupled symmetric systems with variable coefficients such as the linearized Euler
or magnetic induction equations. As can be seen there, further problems can arise
for general symmetric hyperbolic systems in multiple space dimensions, even if energy
stable discretizations are used.}

%% file: Appendix.tex
\section*{Appendix}\label{sec:Appendix}

\subsection*{Technical explanation of the investiagtion in section \ref{sec:Error_Lobatto}}
\revsec{We presented the ideas how to reach  \eqref{eq:Delta} from \eqref{eq:continuous_error}.
Applying the interpolation operator together with discrete norms results in\footnote{
Since $\phi \in \P^N$ and if $a\equiv 1$, the volume term is
\begin{equation*}\label{eq:interpolation_foot}
 \est{ \Ip^N(u^k), \partial_\xi \phi^k}=
 \l(  \vec{\Ip^N(u^k)}, \partial_\xi \vec{\phi}^{k,T} \r)_N =\vec{\phi}^k \mat{D}[^T] \mat{M} \vec{\Ip^N(u^k)}
\end{equation*}
and also the terms \eqref{eq:interpolation2}--\eqref{eq:interpolation4} simplify \revPP{and can be brought together}, see
inter alia \cite{kopriva2017error} for details.
}
\begin{align*}
\est{\partial_t \Ip^N(u^k),\phi^k  }=&
\l( \partial_t \vec{\Ip^N(u^k)}, \vec{\phi}^k \r)_N \nonumber \\ &
+\l\{ \est{ \partial_t \Ip^N(u^k), \phi^k }-
 \l( \partial_t \vec{\Ip^N(u^k)}, \vec{\phi}^k \r)_N  \r\}, \tag{\ref{eq:interpolation}}\\
\frac{1}{2}\est{a^k \Ip^N(u^k), \partial_\xi \phi^k}  &=
\frac{1}{2} \l(  \mat{a^k} \vec{\Ip^N(u^k)} , \partial_\xi \vec{\phi}^k\r)_N \nonumber \\ &
+\frac{1}{2} \l\{ \est{a^k \Ip^N(u^k), \partial_\xi \phi^k} -
\l( \mat{a^k} \vec{\Ip^N(u^k)}  , \partial_\xi \vec{\phi}^k\r)_N  \r\},
\stepcounter{equation}\tag{\theequation}\label{eq:interpolation2}\\
\frac{1}{2}\est{a^k   \partial_\xi \Ip^N(u^k),\phi^k}  &=
\frac{1}{2} \l(  \mat{a^k}  \partial_\xi \vec{\Ip^N(u^k) } , \vec{\phi}^k  \r)_N\nonumber \\ &
+\frac{1}{2} \l\{ \est{a^k \partial_\xi \Ip^N(u^k) , \phi^k  } -
\l(  \mat{a^k} \partial_\xi \vec{\Ip^N(u^k)} , \vec{\phi}^k  \r)_N \r\},
\stepcounter{equation}\tag{\theequation}\label{eq:interpolation3}\\
\frac{1}{2}\est{ \Ip^N(u^k) \partial_\xi a^k,  \phi^k }  &=
\frac{1}{2} \l(  \mat{\Ip^N(u^k)} \partial_\xi \vec{a}^k ,\vec{\phi}^k  \r)_N\nonumber \\ &
+\frac{1}{2} \l\{ \est{ \Ip^N(u^k)  \partial_\xi a^k,\phi^k }  -
 \l(  \mat{\Ip^N(u^k)} \partial_\xi \vec{a}^k ,\vec{\phi}^k  \r)_N \r\}.
\stepcounter{equation}\tag{\theequation}\label{eq:interpolation4}
\end{align*}
It is well known \cite[Section~5.4.3]{canuto2006spectral} that the integration error arising from the use of Gauß quadrature (Gauß-Legendre and Gauß-Lobatto-Legendre) decays spectrally
fast. Indeed, for all $\phi \in \P^N$ and $m\geq 1$,
\begin{equation*}
 \l|\est{u,\phi}-(\vec{u},\vec{\phi})_N \r|\leq C N^{-m}|u|_{H^{m,N-1}(-1,1)}||\phi||_{\L^2(-1,1)},
\end{equation*}
where $C$ is a constant independent of $m$ and $u$.
The curly brackets of \eqref{eq:interpolation}, \eqref{eq:interpolation2}--\eqref{eq:interpolation4} have to
be reformulated. Using
\begin{multline}\label{Soblev_approach}
  \est{ \partial_t \Ip^N(u^k), \phi^k }- \l( \partial_t \vec{\Ip^N(u^k)}, \vec{\phi}^k \r)_N
  =
 \est{ \underbrace{ \partial _t \l( \Ip^N(u^k)-P^m_{N-1} \l( \Ip^N(u^k)\r)\r)}_{=:Q(u^k)}, \phi^k }\\
  - \l( \partial_t \l( \vec{\Ip^N(u^k)}-\vec{P^m_{N-1} \l(\Ip^N(u^k)\r) }  \r), \vec{\phi}^k   \r)_N ,
\end{multline}
 where $P^m_{N-1}$ is the orthogonal projection of $u$ onto $\P^{N-1}$ using the inner product of $H^m(e^k)$, gives a
 new formulation for \eqref{eq:interpolation}. \revPP{The projection operator is defined by
 the classical truncated Fourier series $P^{N-1}u=\sum_{k=0}^{N-1} \hat{u}_k \varPhi_k$ up to order $N-1$ where  Sobolev type orthogonal polynomials
 $\{\varPhi_k \} $ are used
 as basis functions in the Hilbert space $H^m(e^k)$.
 The coefficients are calculated using the inner product of $H^m(e^k)$ given by
 \begin{equation*}
  \est{u,v}_m=\sum_{k=0}^m \int_{e_k} \frac{\d^k u}{\d x^k}(x)\frac{\d^k v}{\d x^k}(x) \d x.
 \end{equation*}
For more details about the projection operator and about approximation results, we strongly recommend \cite[Section 5]{canuto2006spectral}
and also \cite{bernardi1989properties, bernardi1992polynomial}.
 }
 An analogous approach as \eqref{Soblev_approach}
 leads to terms with $Q_1$ for \eqref{eq:interpolation2}, $Q_2$ for \eqref{eq:interpolation3}
 and $Q_3$ for \eqref{eq:interpolation4}.
The $Q_j$ measure the projection error of a polynomial of degree $N$ to a polynomial of degree $N-1$.
Since $u$ and $a$ are bounded, also these
values have to be bounded. This values can be introduced and finally one obtains \eqref{eq:Delta}.} \\
\revsec{
Later, in this section the error of the fluxes hase to be calulated.
We obtain for the left and right boundary:
\begin{align*}
 &\text{left:} &-\Eps_L^1 \l( f^{\mathrm{num},1}_L -
 \frac{1}{2}  a_L^1\Eps_L^1 \r)=-\Eps_L^1  \left( \left( a_L^1\frac{0+\Eps_L^1 }{2}-
 \sigma a_L^1 \frac{\Eps_L^1 }{2} \right)-\frac{a_L^1 \Eps_L^1 }{2} \right) \\
 & &= \frac{\sigma a_L^1}{2}\l(\Eps_L^1\r)^2,
 \\
 &\text{right:}&\Eps_R^K \l(f^{\mathrm{num},K}_R -
 \frac{1}{2} a_R^K\Eps^K_R \r)= \Eps^K_R  \l( \l(  a_R^K\frac{0+\Eps^K_R }{2} +
 \frac{1}{2} \sigma  a_R^K \Eps^K_R  \r)-\frac{\Eps^K_R  a_R^K}{2}\r) \\
 & &= \frac{\sigma  a_R^K}{2} \l(\Eps^K_R  \r)^2.
\end{align*}
}

\subsection*{Technical steps of the development in section \ref{sec:Error_Gauss}}
\revsec{
Here, we are presenting the main steps to reach \eqref{eq:R_absch_2}}.
\revsec{
\begin{align*}
 &\frac{\Delta x_k}{2} \est{\partial_t \Ip^N(u^k), \phi^k }
 +\frac{1}{2} \Bigg( a^k \Ip^N(u^k) \phi^k \Bigg|_{-1}^1 -
 \est{a^k \Ip^N(u^k),\partial_\xi \phi^k  }
 + \est{\partial_\xi \Ip^N(u^k), a^k \phi^k } \\
 &+ \revPP{\est{ \partial_\xi a^k , \phi^k \Ip^N(u^k)} }  \Bigg)
 = -\frac{\Delta x_k}{2} \est{\partial_t \epsilon_p^k, \phi^k }-\frac{1}{2} \l(a^k \epsilon_p^k \phi^k \Bigg|_{-1}^1  \r) +\frac{1}{2} \est{a^k \epsilon_p^k, \partial_\xi \phi^k }
  \\ &-\frac{1}{2} \est{\partial_\xi \epsilon_p^k, a^k \phi^k } -
  \frac{1}{2} \est{ \phi^k \epsilon_p^k, \partial_\xi a^k }.
\end{align*}
Integration-by-parts yields
\begin{equation*}
 -\frac{1}{2} \l(a^k \epsilon_p^k \phi^k \Bigg|_{-1}^1 - \est{a^k \epsilon_p^k, \partial_\xi \phi^k } \r)= -\frac{1}{2} \est{ \partial_\xi (a^k \epsilon_p^k), \phi^k }.
\end{equation*}
With \eqref{eq:interpolation},\eqref{eq:interpolation2}-- \eqref{eq:interpolation4}, one obtains
\begin{equation*} \tag{\ref{eq:Error_Gauss_2}}
 \begin{aligned}
   &\quad \frac{\Delta x_k}{2} \l( \partial_t \vec{\Ip^N (u^k)}, \vec{\phi}^k \r)_N +\vec{\phi}^{k,T} \mat{R}[^T] \mat{B} \l( \vecfnumk \l(\Ip^N(u^k)^-, \Ip^N(u^k)^+\r)-
 \frac{1}{2} \l(\mat{R}\vec{a^k} \r)  \cdot
 \l(\mat{R} \vec{u} \r) \r)
 \\&\quad
 + \underbrace{\l( \frac{1}{2}  a^k \Ip^N(u^k)\phi^k\Bigg|_{-1}^1 -
 \vec{\phi}^{k,T} \mat{R}[^T] \mat{B} \l( \vecfnumk \l(\Ip^N(u^k)^-, \Ip^N(u^k)^+\r)-
 \frac{1}{2} \l(\mat{R}\vec{a^k} \r) \cdot
 \l(\mat{R} \vec{u} \r) \r) \r)}_{=:\epsilon_2^k(a^k)}
 \\&\quad
 - \frac{1}{2} \l(  \mat{a^k} \vec{\Ip^N(u^k)} ,  \partial_\xi \vec{\phi}^k \r)_N
  +\frac{1}{2} \l( \partial_\xi \vec{\Ip^N(u^k)}, \mat{a^k} \vec{\phi}^k    \r)_N
  +\frac{1}{2} \l( \partial_\xi \vec{a}^k , \mat{\Ip^N(u^k)}  \vec{\phi}^k \r)_N\\
 &= \frac{\Delta x_k}{2} \est{ \hat{T}^k(u^k) , \phi^k} +\frac{\Delta x_k}{4} \est{Q_1(u^k), \partial_x \phi^k }
 \\&\quad
 +\frac{\Delta x_k}{4}  \Big\{
  \l( \vec{Q(u^k)}, \vec{\phi}^k \r)_N -
   \l( \vec{Q_1(u^k)}, \partial_x \vec{\phi}^k \r)_N +
   \revPP{ \l(  \vec{Q_2(u^k)}, \mat{a}^k \vec{\phi}^k \r)_N } \\ & \quad +
     \l( \partial_x \vec{a}^k, \mat{Q_3(u^k)}\vec{\phi}^k \r)_N \Big \}
 \end{aligned}
\end{equation*}
 with
 \begin{equation*}
 \begin{aligned}
     \hat{T}^(u^k):=& -\Bigg\{\partial_t \epsilon_p^k+\frac{1}{2} \l( \partial_x \l( a^k \epsilon_p^k\r)
     +\epsilon_p^k \partial_x a^k \partial_x \epsilon_p^k \r)\\
 &+\frac{1}{2} \l( Q(u^k) +a^kQ_2(u^k) +Q_3(u^k) \partial_x a^k \r)
 \Bigg\}.
 \end{aligned}
 \end{equation*}}
\revPP{We transposed every term in \eqref{eq:semi_disc_ver_Gauss}  and subtracted
it from equation \eqref{eq:Error_Gauss_2}.}\revsec{
Using $\epsilon_1^k=\Ip^N(u^k)-U^k$ yields
\begin{align*}
 &\frac{\Delta x_k}{2} \l( \partial_t \vec{\epsilon}_1^k, \vec{\phi}^k \r) +
 \vec{\phi}^{k,T} \mat{R}[^T] \mat{B} \l( \vecfnumk \l((\epsilon_1^k)^-, (\epsilon_1^k)^+\r)-
 \frac{1}{2} \l(\mat{R}\vec{a^k} \r) \cdot
 \l(\mat{R} \epsilon_1^k \r) \r)\\
 & +\epsilon_2^k(a^k)- \frac{1}{2} \l(  \mat{a^k} \vec{\epsilon}_1^k ,  \partial_\xi \vec{\phi}^k \r)_N
  +\frac{1}{2} \l( \partial_\xi  \vec{\epsilon}_1^k , \mat{a^k} \vec{\phi}^k    \r)_N
  +\frac{1}{2} \l( \partial_\xi \vec{a}^k , \mat{\epsilon}_1^k   \vec{\phi}^k \r)_N\\
  =& \frac{\Delta x_k}{2} \est{ \hat{T}^k(u^k) , \phi^k} +\frac{\Delta x_k}{4} \est{Q_1(u^k), \partial_x \phi^k }
  +\frac{\Delta x_k}{4}  \Big\{
  \l( \vec{Q(u^k)}, \vec{\phi}^k \r)_N - \\
 &
   \l( \vec{Q_1(u^k)}, \partial_x \vec{\phi}^k \r)_N +
  \revPP{ \l( \vec{Q_2(u^k)}, \mat{a}^k \vec{\phi}^k \r)_N }+
     \l( \partial_x \vec{a}^k, \mat{Q_3(u^k)}\vec{\phi}^k \r)_N \Big \} .
\end{align*}
Putting $\phi^k=\epsilon_1^k$ results in the energy equation similar to \eqref{eq:energy1}:
\begin{align*}
 &\frac{\Delta x_k}{4} \frac{\d }{\d t}||\epsilon_1^k||_N^2 +
 \vec{\epsilon}_1^{k,T} \mat{R}[^T] \mat{B} \l( \vecfnumk \l((\epsilon_1^k)^-, (\epsilon_1^k)^+\r)-
 \frac{1}{2} \l(\mat{R}\vec{a^k} \r) \cdot
 \l(\mat{R} \vec{\epsilon}_1^k \r) \r)\\
 &
 +\epsilon_2^k(a^k)- \frac{1}{2} \l( \mat{a^k} \vec{\epsilon}_1^k , \partial_\xi \vec{\epsilon}_1^k \r)_N
  +\frac{1}{2} \l( \partial_\xi  \vec{\epsilon}_1^k , \mat{a^k} \vec{\epsilon}_1^k    \r)_N
  +\frac{1}{2} \l( \partial_\xi \vec{a}^k , \mat{\epsilon}_1^k   \vec{\epsilon}_1^k \r)_N\\
  =& \frac{\Delta x_k}{2} \est{ \hat{T}^k(u^k) , \epsilon_1^k} +\frac{\Delta x_k}{4} \est{Q_1(u^k), \partial_x \epsilon_1^k } \\
 &+\underbrace{\frac{\Delta x_k}{4}  \Big\{
  \l( \vec{Q(u^k)}, \vec{\epsilon}_1^k \r)_N -
   \l( \vec{Q_1(u^k)}, \partial_x \vec{\epsilon}_1^k \r)_N +
  \revPP{  \l(  \vec{Q_2(u^k)}, \mat{a}^k\vec{\epsilon}_1^k \r)_N }+
     \l( \partial_x \vec{a}^k, \mat{Q_3(u^k)}\vec{\epsilon}_1^k \r)_N \Big \} }_{\hat{Q}^k}.
\end{align*}
Together with \eqref{eq:Deleting_a_terms}, one obtains
\begin{multline}\label{eq:Gauss_error}
  \frac{\Delta x_k}{4} \frac{\d }{\d t}||\epsilon_1^k||_N^2 +
 \vec{\epsilon}_1^{k,T} \mat{R}[^T] \mat{B} \l( \vecfnumk \l((\epsilon_1^k)^-, (\epsilon_1^k)^+\r)-
 \frac{1}{2} \l(\mat{R}\vec{a^k} \r)  \cdot
 \l(\mat{R} \vec{\epsilon}_1^k \r) \r)
 \\+\epsilon_2^k(a^k)
  +\frac{1}{2} \l( \partial_\xi \vec{a}^k , \mat{\epsilon}_1^k   \vec{\epsilon}_1^k \r)_N
  = \frac{\Delta x_k}{2} \est{ \hat{T}^k(u^k) , \epsilon_1^k} +\frac{\Delta x_k}{4} \est{Q_1(u^k), \partial_x \epsilon_1^k }
 +\hat{Q}^k.
\end{multline}
Summing this up over all elements results in
\begin{equation*}
 \begin{aligned}
 &\quad \frac{1}{2}   \frac{\d}{\d t} \sum_{k=1}^K \frac{\Delta x_k}{2}||\epsilon_1^k||_N^2+\sum_{k=1}^K
 \vec{\epsilon}_1^{k,T} \mat{R}[^T] \mat{B}
 \Bigg(\vecfnumk
 \l( (\epsilon_1^k )^{-},(\epsilon_1^k )^{+}\r)  \\&\quad-
 \frac{1}{2 } \l(\mat{R}\vec{a^k} \r)
 \l(\mat{R} \vec{\epsilon}_1^k \r)\Bigg)
 +\frac{1}{2} \sum_{k=1}^K \frac{\Delta x_k}{2}
 \revPP{\l( \partial_x \vec{a}^k , \mat{\epsilon}_1^k  \vec{\epsilon}_1^k  \r)_N }
  +\sum_{k=1}^K\frac{\Delta x_k}{4}
     \epsilon_2^k(a^k)
  \\&=
  \sum_{k=1}^K\frac{\Delta x_k}{2} \est{\hat{T}^k(u^k), \epsilon_1^k}
  + \sum_{k=1}^K\frac{\Delta x_k}{4}  \est{Q_1(u^k), \partial_x \epsilon_1^k}
  + \sum_{k=1}^K  \hat{Q}^k .
 \end{aligned}
\end{equation*}
Applying the same approach like in equations \eqref{eq:T_est}--\eqref{Eq:est_Q_1}
and the fact that $\epsilon_1 \in \P^N$, it is$
||\partial_x \vec{\epsilon}_1^k||_N^2 \leq c_1 N^2   ||\vec{\epsilon}_1^k||_N^2$
and we get finally \eqref{eq:R_absch_2}.}\\

\subsubsection*{Calculating the fluxes from Table \ref{Table:error_flux}}
\revsec{
\begin{itemize}
\item Split central flux $\fnum(u_-, u_+)= \frac{a_-u_-+a_+u_+}{2}$:
One obtains
\begin{align*}
 &\frac{1}{2}\l(a_R^{k-1}\Eps^{k-1}_R+a_L^{k}\Eps^k_L \r) \l(\Eps^{k-1}_R-\Eps_L^k \r)-\frac{1}{2} \l( a_R^{k-1}\l( \Eps_R^{k-1 } \r)^2    -a_L^k \l(\Eps^k_L \r)^2 \r)\\
 =& \frac{1}{2} \l(a_R^{k-1} \left(\Eps_R^k\right)^2-a_R^{k-1}\Eps_L^k \Eps^{k-1}_R + a_L^k \Eps^k_L \Eps^{k-1}_R-a_L^k \l( \Eps_L^k \r)^2 \r) \\
 &-\frac{1}{2} \l( a_R^{k-1}\l( \Eps_R^{k-1 } \r)^2    -a_L^k \l(\Eps^k_L \r)^2 \r)
 = \frac{1}{2} \Eps_L^k \Eps_R^{k-1} \l(a_L^k-a_R^{k-1} \r) = 0
\end{align*}
and
\begin{align*}
 &\text{ left:}& -\Eps_L^1 \l( f^{\mathrm{num},1}_L-\frac{1}{2}a_L^1\Eps_L^1 \r)=-\Eps_L^1 \l( \frac{a^1_L}{2}\Eps_L^1-\frac{a_L^1}{2} \Eps_L^1 \r)=0, \\
 &\text{ right:}& \Eps^K_R  \l( f^{\mathrm{num},K}_R -\frac{1}{2}a_R^K \Eps^K_R \r)=\frac{1}{2}
 \l(\Eps_R^K \r)^2 \l(a^K_R-a_R^K \r) =0 .
\end{align*}
\item Edge based upwind flux $\fnum(u_-,u_-)=a(x)u_-$:
It is
\begin{align*}
 &a^{k-1}(x_R)\Eps^{k-1}_R \l( \Eps^{k-1}_R-\Eps_L^k \r)-\frac{1}{2} \l( a_R^{k-1}\l( \Eps_R^{k-1 } \r)^2    -a_L^k \l(\Eps^k_L \r)^2 \r)\\
  =& a^{k-1}(x_R) \left(\Eps_R^{k-1}\right)^2-a^{k-1}(x_R)\Eps_L^k \Eps^{k-1}_R
 -\frac{1}{2} \l( a_R^{k-1}\l( \Eps_R^{k-1 } \r)^2 -a_L^k \l(\Eps^k_L \r)^2 \r) \\
 =&\l( \Eps_R^{k-1 } \r)^2 \l( a^{k-1}(x_R) -\frac{1}{2} a_R^{k-1} \r) +\frac{\Eps^k_L}{2} \l(a_L^k \Eps_L^k-a^{k-1}_R \Eps_R^{k-1} \r)
 \\ =& \frac{1}{2} a_R^{k-1} \jump{E_R^{k-1}}^2
\end{align*}
and
\begin{align*}
 &\text{ left:}& -\Eps_L^1 \l( f^{\mathrm{num},1}_L-\frac{1}{2}a_L^1\Eps_L^1 \r)&=\frac{1}{2} \l(\Eps^1_L \r)a_L^1, \\
 &\text{ right:}& \Eps^K_R  \l( f^{\mathrm{num},K}_R -\frac{1}{2}a_R^K \Eps^K_R \r)&=\l(\Eps^k_R \r)^2 \l(a^K(x_R)-\frac{a_R^{K}}{2} \r)
 & & = \l(\Eps^k_R \r)^2 \l(\frac{a_R^{K}}{2} \r).
\end{align*}
\item Split upwind flux $\fnum(u_-,u_-)=a_-u_-$:
It is
\begin{align*}
 &a^{k-1}_R\Eps^{k-1}_R \l( \Eps^{k-1}_R-\Eps_L^k \r)-\frac{1}{2} \l( a_R^{k-1}\l( \Eps_R^{k-1 } \r)^2    -a_L^k \l(\Eps^k_L \r)^2 \r)\\
  =& \frac{1}{2}\l( a^{k-1}_R\l(\Eps_R^{k-1}\right)^2-2a^{k-1}_R\Eps_L^k \Eps^{k-1}_R +a_L^k \l(\Eps^k_L\r)^2 \r)
  = \frac{1}{2} a_R^{k-1} \jump{E_R^{k-1}}^2,
\end{align*}
where we used in the last step the assumption about the exactness of the interpolation and the continuity of $a$.
At the boundaries we get
\begin{align*}
 &\text{ left:}& \frac{a_L^1 }{2}  \l(\Eps^1_L\r)^2,\\
 &\text{ right:}&  \frac{a_R^K }{2} \l(\Eps^K_R\r)^2.
\end{align*}
\end{itemize}
}